\documentclass[11pt, epsfig]{article}
\usepackage{epsfig, amsmath, amssymb, amsthm}
\def\theequation{\arabic{section}.\arabic{equation}}\def\theequation{\arabic{section}.\arabic{equation}}
\renewcommand{\theequation}{\thesection.\arabic{equation}}

\def\p{\partial}
\def\Re{{\rm I}\!{\rm R}}
\def\Na{{\rm I}\!{\rm N}}
\def\lbl{\label}
\def\be{\begin{equation}}
\def\ee{\end{equation}}
\def\lbl{\label}

\def\mG{\mathcal{G}}
\def\mI{\mathcal{I}}
                                                                        
\title{Multimodal oscillations in systems with strong contraction}
\author{ 
Georgi Medvedev and Yun Yoo \thanks{
Department of Mathematics, Drexel University, 3141 Chestnut Street,
Philadelphia, PA 19104, {\it medvedev@drexel.edu}
 }
 }

\begin{document}
\maketitle

\begin{abstract}
One- and two-parameter families of flows in $\Re^3$
near an Andronov-Hopf bifurcation (AHB) are investigated in this work. 
We identify conditions on the 
global vector field, which yield a rich family of multimodal orbits 
passing close to a weakly unstable
saddle-focus and  perform a detailed asymptotic analysis of the
trajectories in the vicinity of the saddle-focus. Our analysis covers both
cases of sub- and supercritical AHB. For the supercritical case, we find
that the periodic orbits born from the AHB are bimodal when viewed in the
frame of coordinates generated by the linearization about the bifurcating
equilibrium. If the AHB is subcritical, it
is accompanied by the  appearance of multimodal orbits, which consist
of long series of nearly harmonic oscillations separated by large
amplitude spikes. We analyze the dependence of the interspike intervals
(which can be extremely long) on the control parameters. In particular,
we show that the interspike intervals grow logarithmically as the boundary
between regions of sub- and supercritical AHB is approached in the parameter
space. We also identify a window of complex and possibly chaotic oscillations
near the boundary between the regions of sub- and supercritical AHB
and explain the mechanism generating these oscillations. This work is motivated 
by the numerical results for a finite-dimensional approximation of a free 
boundary problem modeling solid fuel combustion.
\end{abstract}

\section{Introduction}

Systems of differential equations in both finite- and infinite-dimensional settings
close to an AHB have been subject to intense research 
due to their dynamical complexity and importance in applications. The latter
range from models in fluid dynamics \cite{MarsdenMcCracken} to those in the life sciences, 
in particular, in computational neuroscience 
\cite{BER, DK05, DRSE, GW, KR83, MC, RE89, ROWK}.
When the proximity to the AHB 
coincides with certain global properties of the vector field,
it may result in a very complex dynamics \cite{BS, DENG90, DENG95, GW, Hirschberg_Knobloch}.
The formation of the Smale horseshoes in systems 
with a homoclinic orbit to a saddle-focus equilibrium  
provides one of the most representative examples of this type 
\cite{AAIS}.
Canard explosion in relaxation systems affords another example \cite{CDD, KS01b}. 
Recent studies of relaxation systems,
motivated mainly by applications in the life sciences, have revealed that
the proximity to an AHB has a significant impact on the system
dynamics. It manifests itself as a family of multimodal periodic solution
that are composed of large-amplitude relaxation oscillations 
(generated by the global structure of the vector field)
and small-amplitude nearly harmonic oscillations (generated by the vector field near the
equilibrium close to the AHB)
\cite{DRSE, KOP, MC, MSLG, PSS, VV} (see Figure \ref{f.1}). 
These families of solutions possess rich bifurcation structure.
A remarkable example of an infinite-dimensional system close to the AHB
has been recently studied by Frankel and Roytburd 
\cite{FR05a, FR05b, FR03, FR95, FR94}. They derived and systematically studied 
a model of solid fuel combustion 
in the form of a free boundary problem for a $1D$ heat equation with
nonlinear conditions imposed at the free boundary modeling the interface between 
solid fuel mixture and a solid product.
The investigations of this model revealed a wealth of spatial-temporal patterns
ranging from a uniform front propagation to periodic and aperiodic front oscillations. 
The transitions between different dynamical regimes involve a variety of nontrivial bifurcation
phenomena including period-doubling cascades, period-adding sequences, 
and windows of chaotic dynamics. To elucidate the mechanisms responsible for 
different dynamical regimes and transitions between them, Frankel and Roytburd
employed pseudo-spectral techniques to derive a finite-dimensional approximation
for the interface dynamics in the free boundary problem \cite{FR95}. 
As shown in \cite{FKRT, FR95},
a system of three ordinary differential equations captured the essential features of 
the bifurcation structure of the infinite-dimensional problem.
The numerical bifurcation analysis of the finite-dimensional approximation revealed a rich
family of multimodal periodic solutions similar to those reported in the context of relaxation
systems near the AHB \cite{FKRT}. The bifurcation diagrams presented in \cite{FKRT}
and in \cite{MC} share a striking similarity, despite the absence of any apparent common 
structures in the underlying models (except to the proximity to the AHB). 
In particular, in both models, topologically distinct multimodal periodic solutions
are located on isolas, closed curves in the parameter space.
The methods of analysis of the mixed-mode solutions in \cite{KOP, MC, MSLG, WESCH} used
in an essential way the relaxation structure present in these problems. These approaches
can not be applied directly to analyzing the model in \cite{FKRT}, because it is not a priori clear 
what creates the separation of the time scales in this model, in spite of the evident fast-slow
character of the numerical solutions. This is partly due to the spectral method, which
was used to derive the system of equations in \cite{FKRT}: while it has captured well
the finite-dimensional attractor of the interface dynamics, it has disguised 
the structure of the physical model. One of the goals of the present paper is to identify 
the structure responsible for the generation of the multimodal oscillations 
in a finite-dimensional model for the interface dynamics and 
to relate it to those studied in the context of relaxation oscillations.

The family of flows in \cite{FKRT} includes 
in a natural way two types of the AHBs. Depending on the parameter values,
the equilibrium of the system of ordinary differential equations in \cite{FKRT}  
undergoes either a sub- or a supercritical AHB. A similar
situation is encountered in certain neuronal models (see, e.g., \cite{DK05,OP}).
In either case, the global multimodal periodic solutions are created after
the AHB. However, 
in the case of a supercritical bifurcation, they are preceded by a series of 
period-doubling bifurcations of small amplitude limit cycles, arising from
the AHB. 
On the other hand, in the subcritical case, the AHB gives rise to 
multimodal solutions, whose lengths and time intervals between successive large
amplitude oscillations can be very long. In the present paper,
we perform a detailed asymptotic analysis of the trajectories in a class of systems motivated
by the problem in \cite{FKRT}. Our analysis includes both cases of the sub- and supercritical AHBs. 
We also investigate the dynamical regimes arising near the border between the regions of sub- and
supercritical AHB. This region in the parameter space contains a number of nontrivial oscillatory
patterns including multimodal trajectories with substantial time intervals between successive
spikes, irregular, and possibly chaotic oscillations, as well as a family of periodic orbits 
undergoing a cascade of period-doubling bifurcations. Our analysis shows that these dynamical
patterns and the order in which they appear under the variation of the control parameters are
independent on the details of the model, but are characteristic to the transition from sub-
to supercritical AHB.

The outline of the paper is as follows. After introducing the model and rewriting it in the
normal coordinates, we present a set of the numerical experiments to be explained in the remainder of the paper.
Then we state our results for each of the following cases: supercritical AHB,
subcritical AHB, and the transition layer between the regions of sub- and supercritical AHB. In Section 3, we analyze
the local behavior of trajectories near a weakly unstable saddle-focus. 
The local expansions used in this section are similar to those used in \cite{CMP} for analyzing 
the AHB using the method of averaging. However, rather than establishing existence of periodic solutions, the goal
of the present section is to approximate the trajectories near the saddle-focus. For this,
after rescaling the variables and recasting the system into cylindrical coordinates, we reduce
the system dynamics to a $2D$ slow manifold. By integrating the leading order approximation of the
reduced system, we obtain necessary information about the local behavior of trajectories. The results
of this section are summarized in Theorem 3.1. In Section 4, we  study oscillatory patterns generated by the
class of systems under investigation. It is divided into three subsections devoted to the oscillations
triggered by the supercritical AHB (\S 4.1), subcritical AHB (\S 4.2), and those found in the transition region
between sub- and supercritical AHB (\S 4.3). In the supercritical case, we show that the oscillations just
after the AHB are already bimodal. Generically, one needs to use two harmonics to describe the
limit cycle born at the supercritical AHB in $\Re^3$. We give a geometric interpretation of this
effect, which we call the frequency doubling, due to the fact that the second frequency is twice as large
as that predicted by the Hopf Bifurcation Theorem \cite{MarsdenMcCracken}. 
We also compute the curvature and the torsion of the periodic orbit as a curve
in $\Re^3$. The latter is useful for the geometric explanation of the frequency doubling. In Subsection 4.2,
we study the subcritical case. We show that under two general assumptions on the global vector field,
the presence of the return mechanism ({\bf G1}) and the strong contraction property ({\bf G2}), the
subcritical AHB results in sustained multimodal oscillations. Even though the oscillations may not be periodic
(our assumptions do not warrant periodicity), the time intervals between consecutive spikes of the resultant motion comply to
the uniform bounds given in Theorem 4.1. This subsection also contains the definitions of the multimodal oscillations
and the precise formulation of the assumptions on the global vector field. The proof of Theorem 4.1
is relegated to Section 5. Subsection 4.3 describes the
transition between sub- and supercritical AHB. This transition contains a distinct bifurcation scenario, when
a small positive real part of the complex conjugate pair of eigenvalues is positive and fixed and the first 
Lyapunov coefficient changes sign. The numerical simulations show that this results 
in the formation of the chaotic attractor followed
by the reverse period-doubling cascade. To explain this bifurcation sequence, we derive a $1D$ first return map.
The asymptotic analysis in this subsection is complemented by numerical extension of the map to the region
inaccessible by local asymptotic expansions. The first return map obtained by the combination of the analytic
and numerical techniques reveals the principal traits of the transition from sub- to supercritical AHB.
Finally, the discussion of the results of the present paper and their relation to the previous work is
given in Section 6.  

\section{The model and numerical results}
In the present section, we formulate the model and present a set of
numerical results, which motivated our study. The following system of three ODEs was derived in 
\cite{FKRT}, as a finite-dimensional
approximation for the interface dynamics in a free boundary problem modeling solid fuel 
combustion:
\begin{eqnarray} \lbl{6.1}
\dot v_1 &=& \frac{3(v_3+v_2-v_1)-\nu k(v_1)-v_1^2}{ \nu k^\prime(v_1)},\\ \lbl{6.2}
\dot v_2 &=& v_3-v_1,\\ \lbl{6.3}
\dot v_3 &=& 9(v_1-v_3)-6v_2+\nu(v_1+1)k(v_1)+2v_1^2.
\end{eqnarray}
Here, $v_1(t)$ approximates the velocity of the interface between the solid fuel and the
burnt material. Functions $v_2(t)$ and $v_3(t)$ are
the first two coefficients in the series expansion for the spatial temperature profile
with respect to the basis of Chebyshev-Laguerre polynomials. Equations (\ref{6.1})-(\ref{6.3})
are obtained by projecting the original infinite-dimensional problem onto a finite-dimensional
function space and using the method of collocations for determining the unknown coefficients.
An important ingredient of the model, nonlinear kinetic function $k(v)$ is given by
\be\lbl{6.4}
k(v)={(1-v)^{p}-(1-v)^{-1}\over p+1}.
\ee
\begin{figure}
\begin{center}
{\bf a}\epsfig{figure=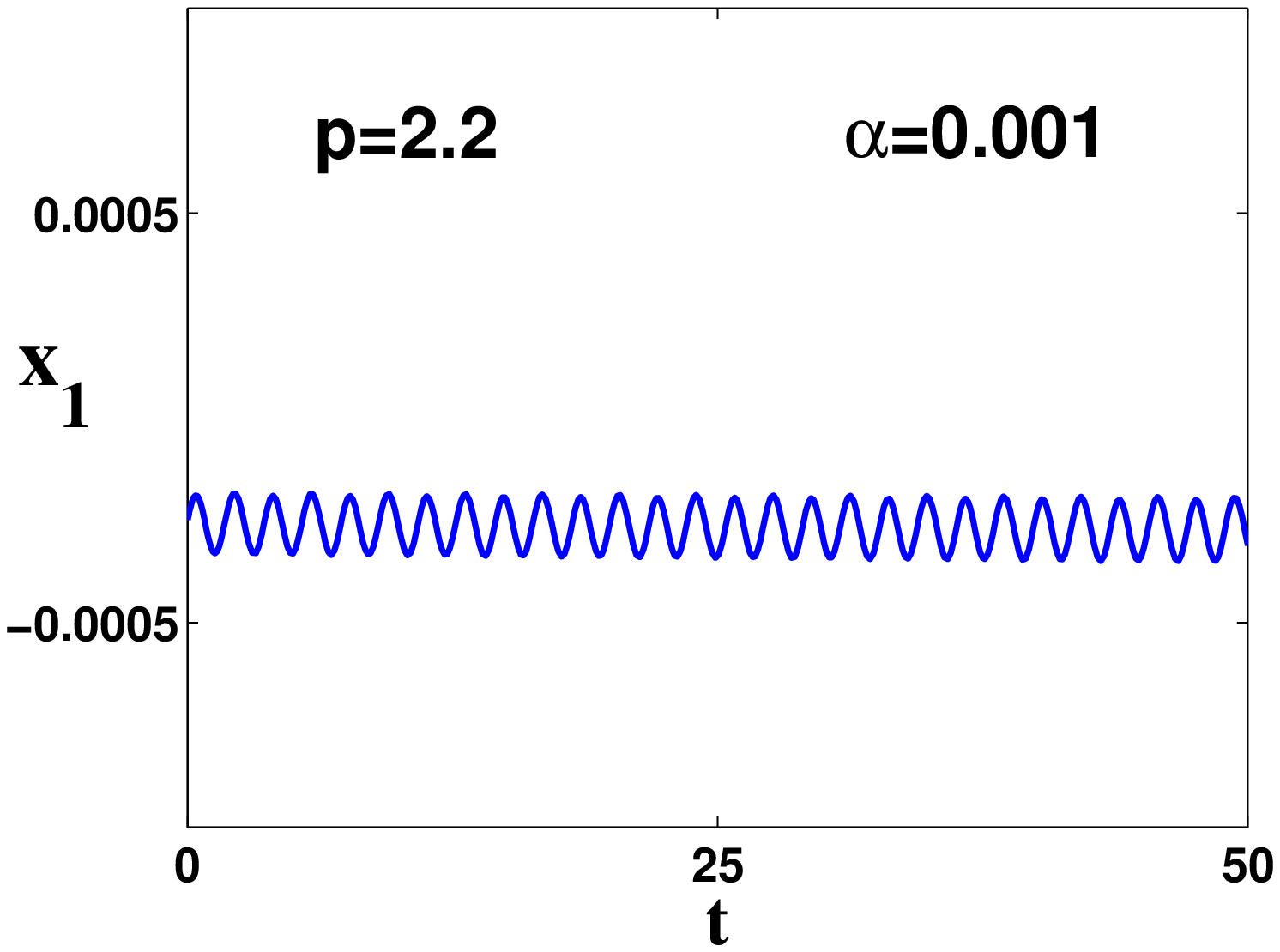, height=1.8in, width=2.5in, angle=0} 
{\bf b}\epsfig{figure=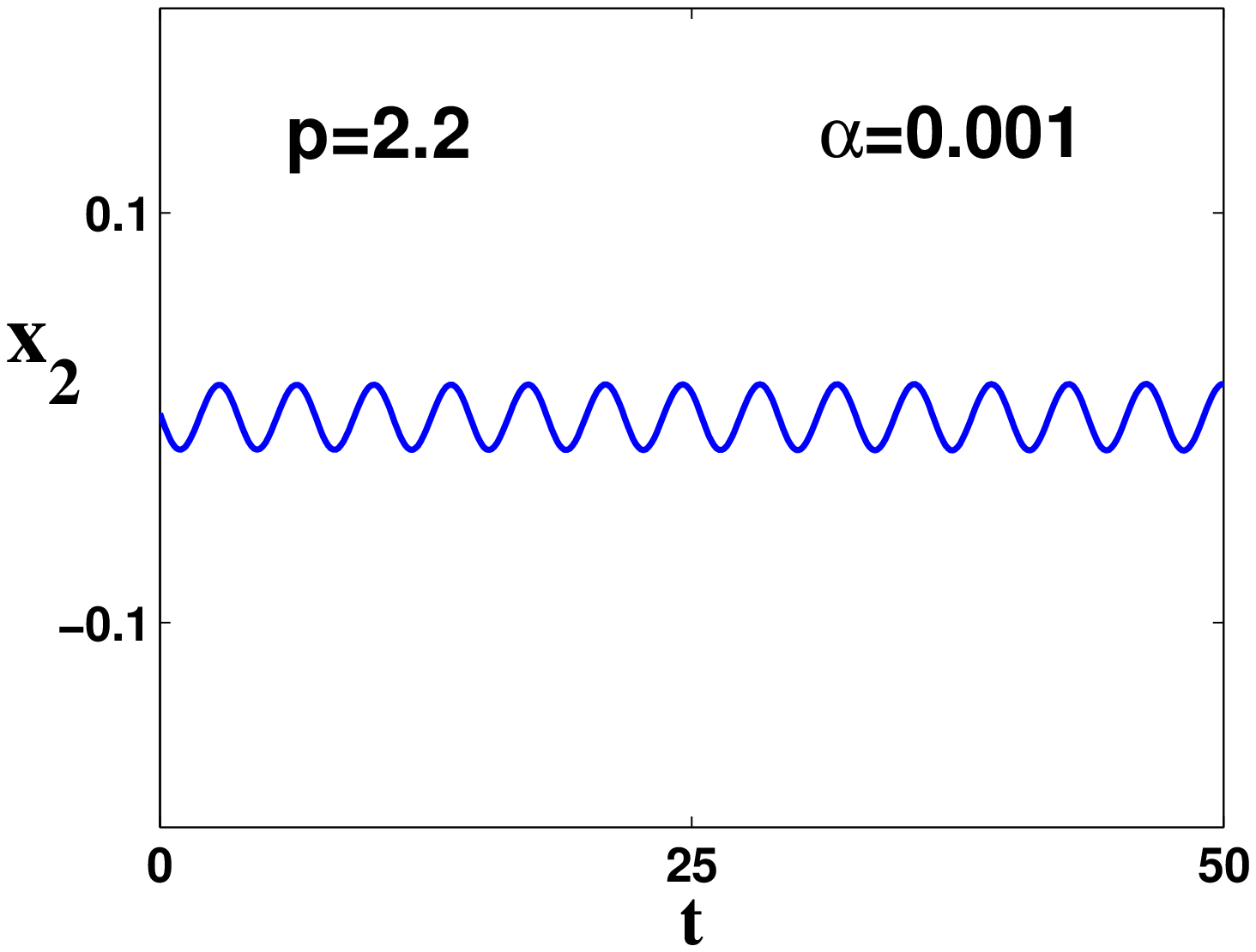, height=1.8in, width=2.5in, angle=0}\\
{\bf c}\epsfig{figure=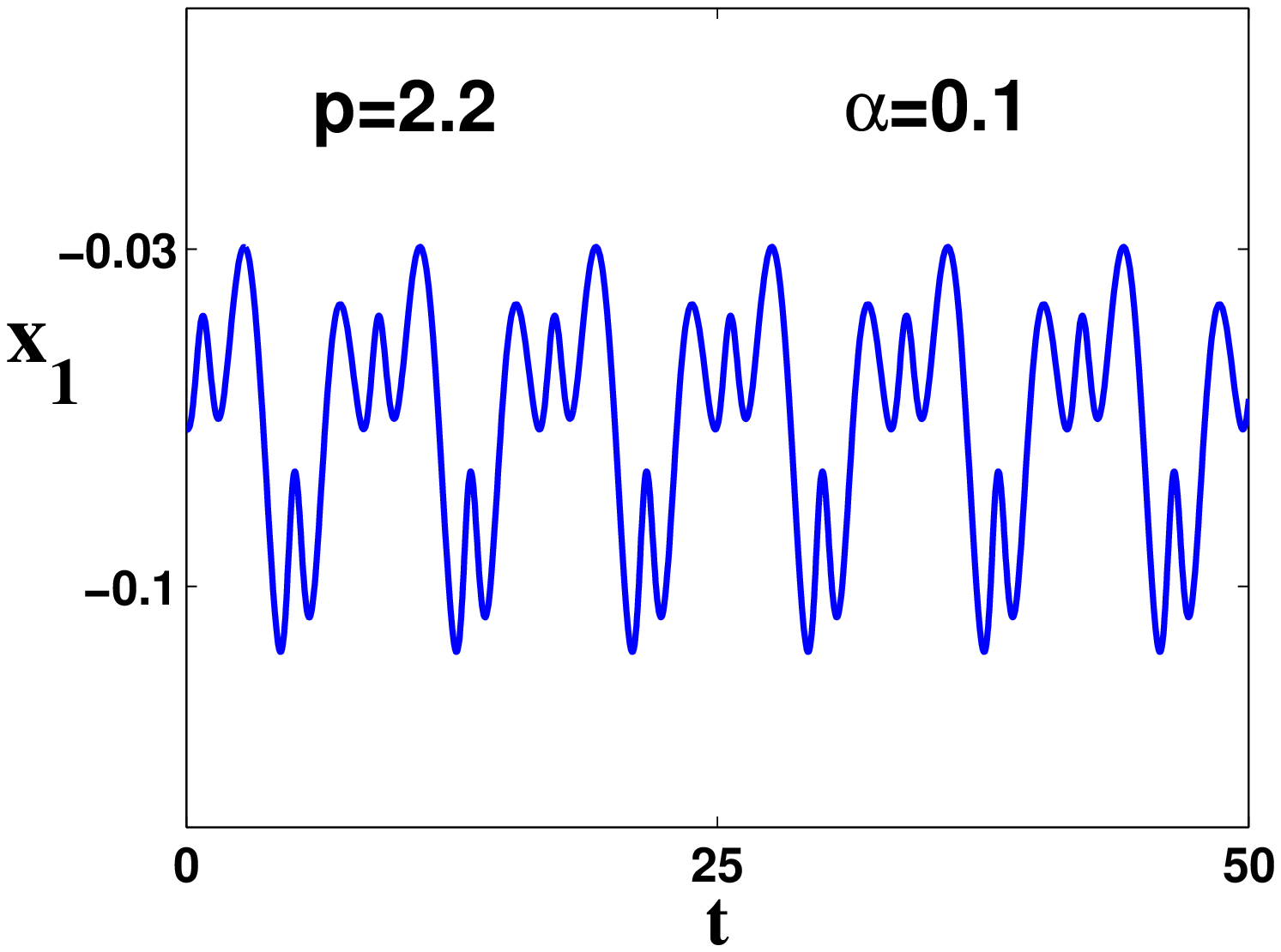, height=1.8in, width=2.5in, angle=0}
{\bf d}\epsfig{figure=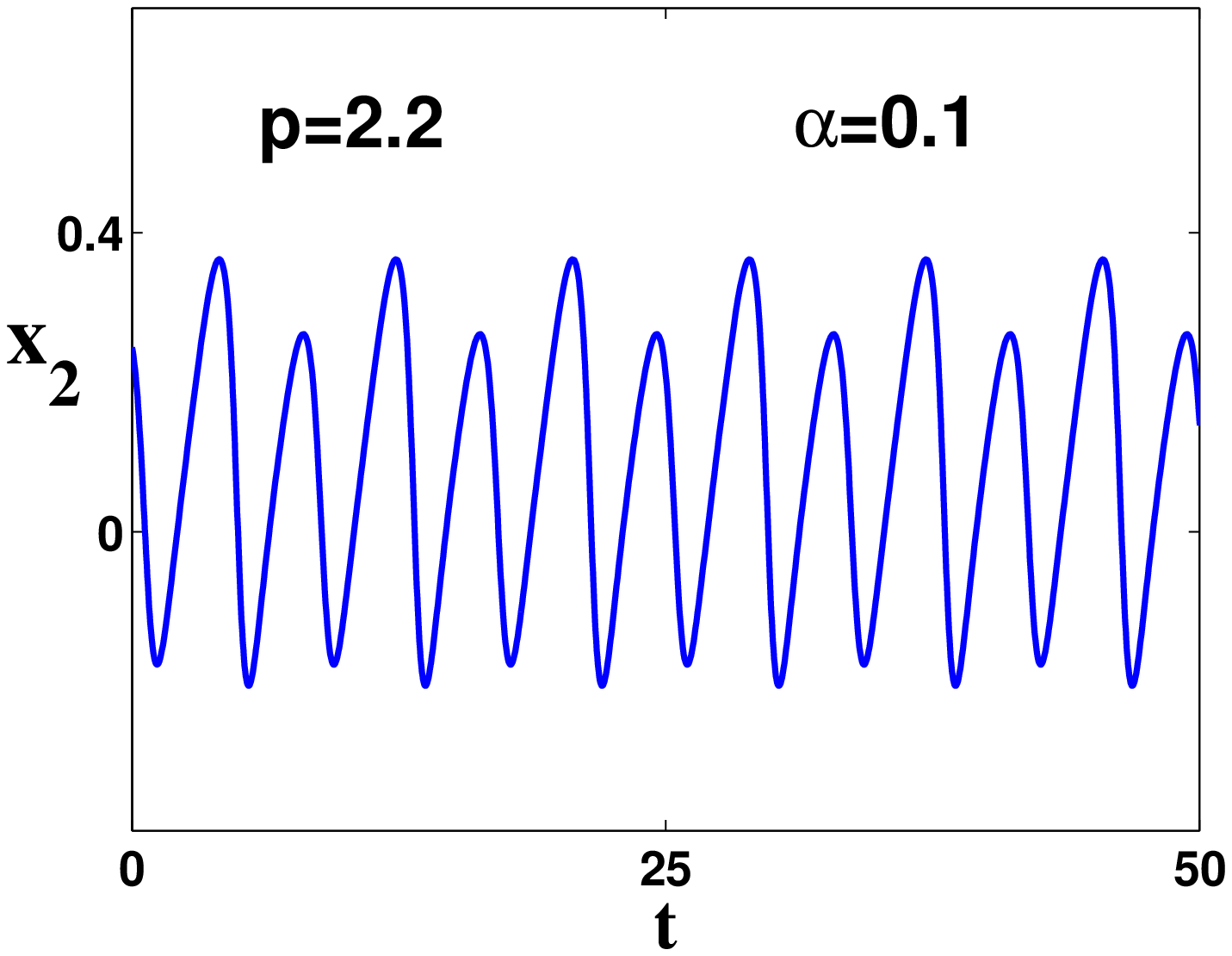, height=1.8in, width=2.5in, angle=0}\\
{\bf e}\epsfig{figure=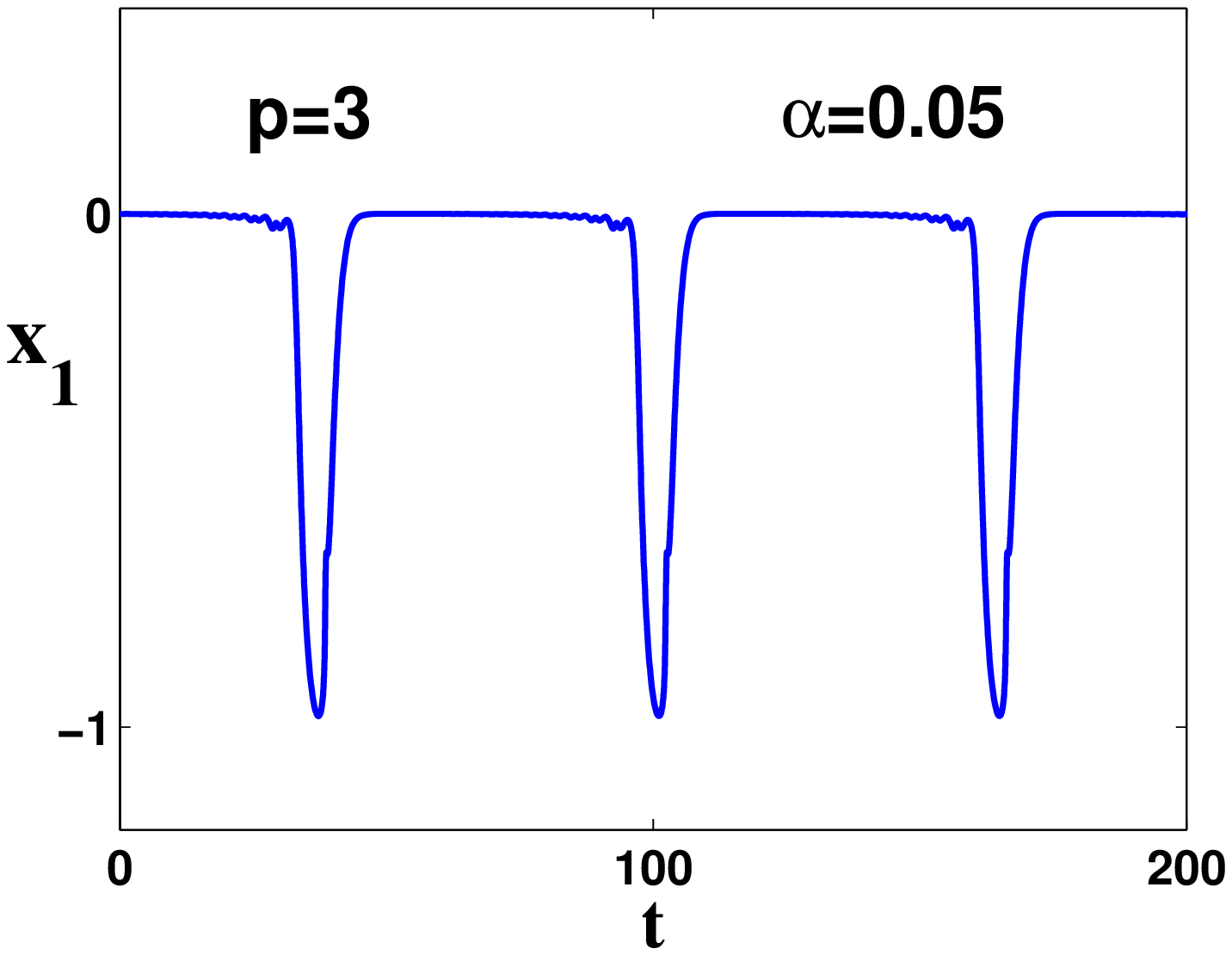, height=1.8in, width=2.5in, angle=0}
{\bf f}\epsfig{figure=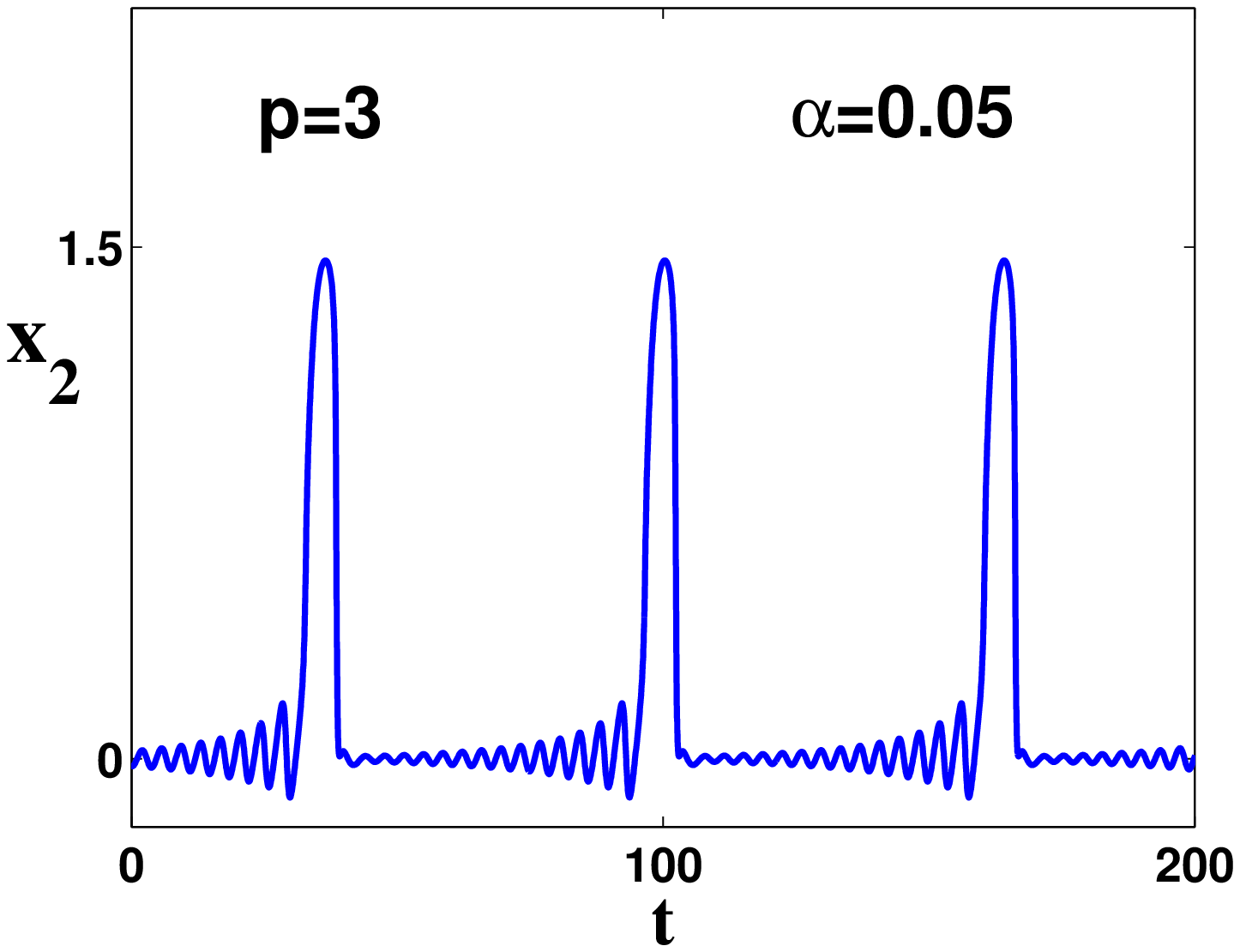, height=1.8in, width=2.5in, angle=0}\\
{\bf g}\epsfig{figure=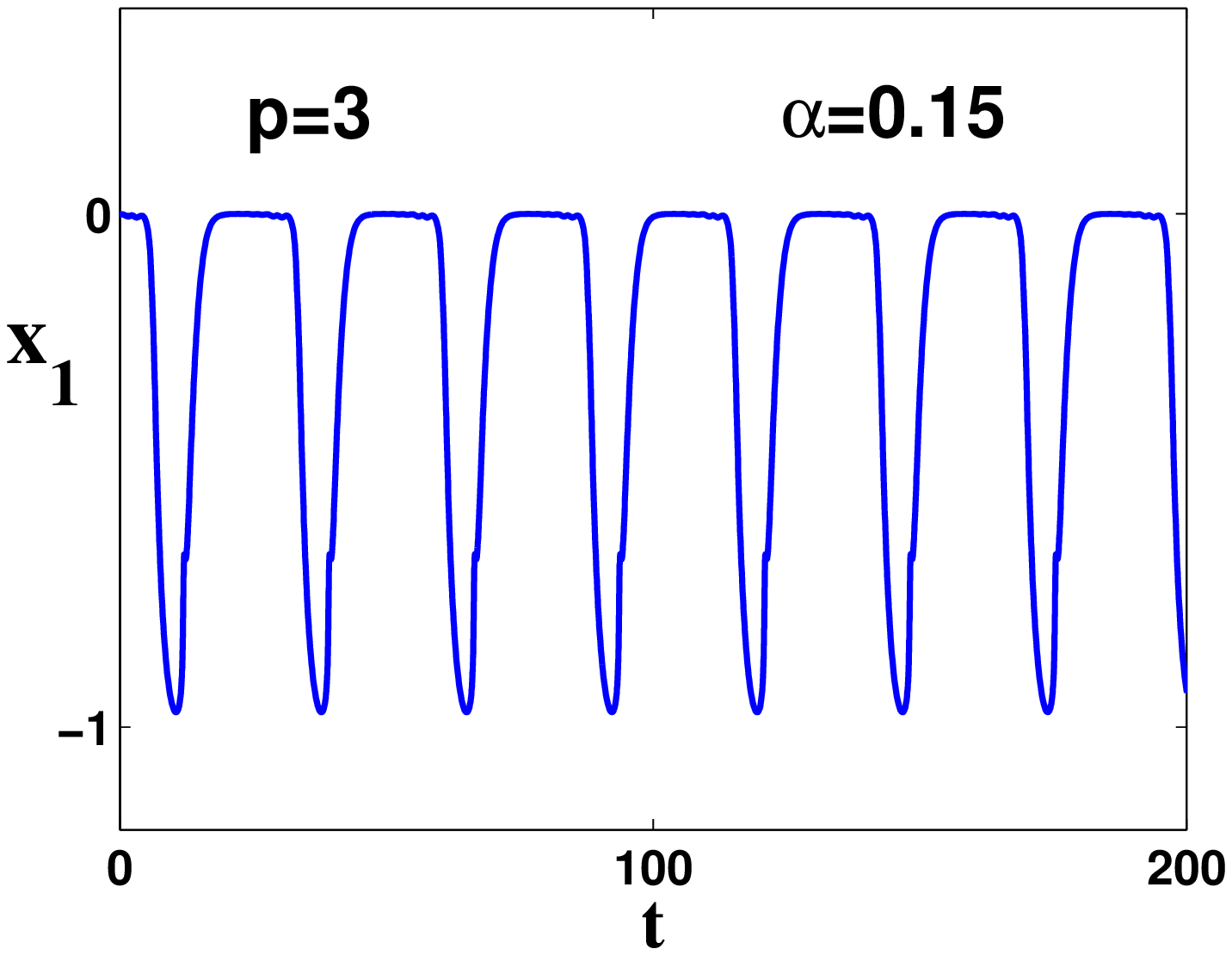, height=1.8in, width=2.5in, angle=0}
{\bf h}\epsfig{figure=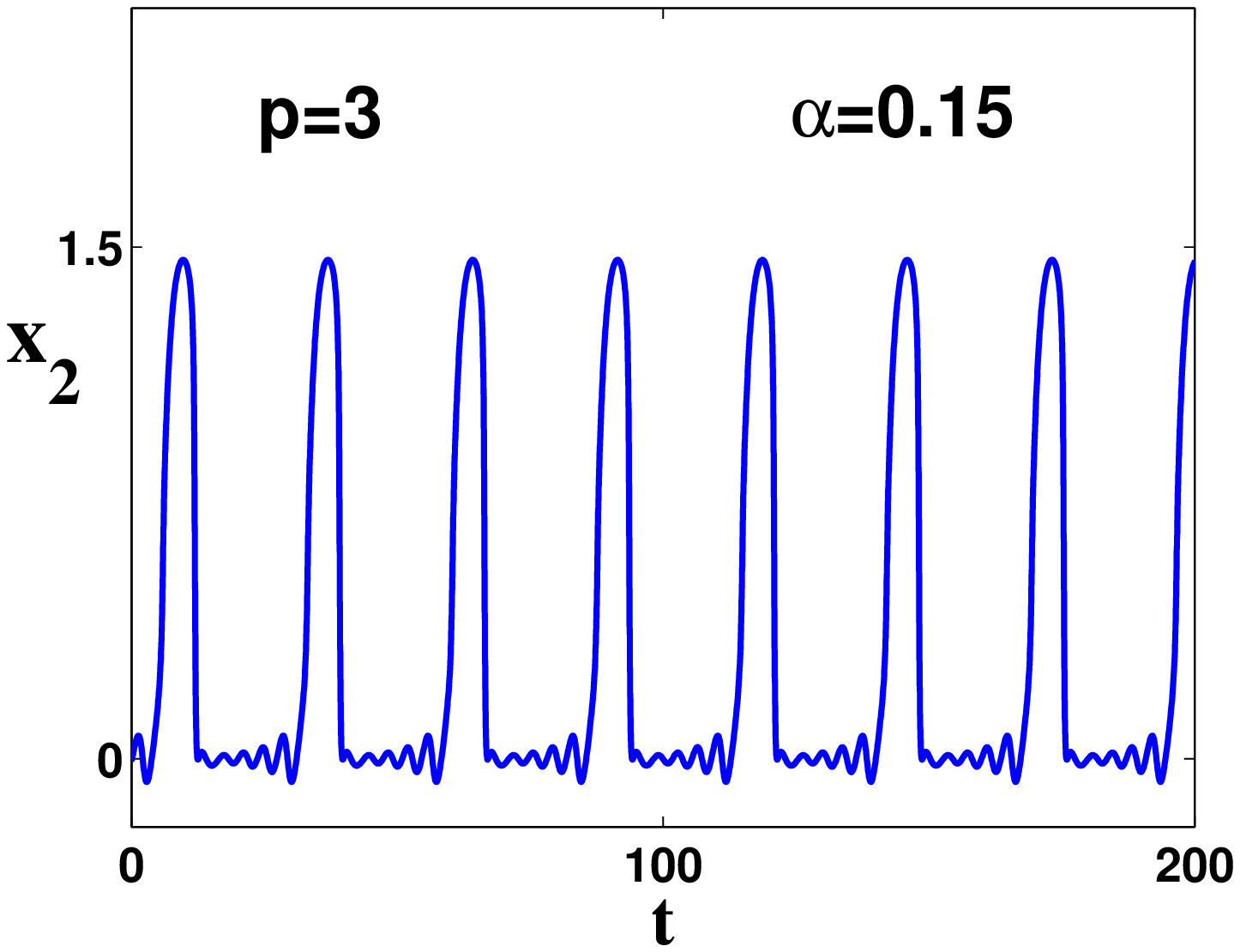, height=1.8in, width=2.5in, angle=0}
\end{center}
\caption{Stable periodic solutions of system of differential equations (\ref{6.5}). Plots
in a-d show solutions near a supercritical AHB ($p=2.2$) and those in e-h
correspond to the subcritical case ($p=3$). The time series of $x_1$ are presented in the
left column and those of $x_2$ are given in the right column. The time series for 
$x_3$ and $x_2$ are qualitatively similar. 
}\label{f.1}
\end{figure} 
It reflects the dependence of the velocity of propagation on the temperature in the front. This 
relation is very complex and is not completely understood at present. To account for a range of 
possible kinetic mechanisms, the model includes two control parameters: $\nu$ and $p$.
Both the interface velocity and the temperature profile are calculated in the uniformly
moving frame of reference. Therefore, system of equations (\ref{6.1})-(\ref{6.3}) describes
the deviations of the interface dynamics from that of a front traveling with constant
speed. In particular, periodic and aperiodic oscillations generated by (\ref{6.1})-(\ref{6.3})
correspond to complex spatio-temporal patterns in the infinite-dimensional model. The numerical
results presented in \cite{FKRT} show a remarkable similarity between the complex patterns 
generated by the infinite- and finite-dimensional models 
and between the scenarios for transitions between different regimes in both models.
For more information about the derivation of the free-boundary problem for solid fuel combustion
and its finite-dimensional approximation, we refer the reader to \cite{FKRT} and bibliography
therein.

A simple inspection of (\ref{6.1})-(\ref{6.3}) shows that it has an equilibrium at the origin $O=(0,0,0)^T$ for all
values of $\nu$ and $p>0$. We linearize (\ref{6.1})-(\ref{6.3}) about the equilibrium at the origin
$$
\dot v = A(\nu)v+ \dots,\quad\mbox{where}\quad
A(\nu) =\left( \begin{array}{ccc}
\ \frac{3-\nu}{\nu} &  \frac {-3}{\nu}  &  \frac {-3}{\nu} \\
-1 & 0 & 1 \\
9-\nu & -6  &  -9 \end{array} \right),
\quad v=\left(v_1, v_2, v_3\right)^T.
$$
Near $\nu_{AH}={1\over 3}$, Jacobian matrix $A(\nu)$ has a negative eigenvalue $-1$ and a pair of complex
conjugate eigenvalues $\lambda=a(\nu) + ib(\nu)$ and $\bar\lambda$. The latter crosses the imaginary
axis transversally at $\nu=\nu_{AH}$:
$$
a\left(\nu_{AH}\right)=0\qquad\mbox{and}\qquad a^\prime\left(\nu_{AH}\right)\ne 0.
$$
Therefore, the equilibrium of (\ref{6.1})-(\ref{6.3}) undergoes an AHB. In the  neighborhood of $\nu_{AH}$,
$\alpha=a(\nu)$ defines a smooth invertible function. We shall use $\alpha$ as a new control parameter.
After a linear coordinate transformation,
system of equations (\ref{6.1})-(\ref{6.3}) has the following form:  
\be\lbl{6.5} 
\dot x = 
\left(\begin{array}{ccc}
                      -1 & 0 & 0\\ 
                       0 & \alpha & -b(\alpha)\\ 
                       0 & b(\alpha) &  \alpha
      \end{array}\right) x + h(x,\alpha),
\quad x=\left(x_1, x_2, x_3\right)^T, \quad h=\left(h_1, h_2, h_3\right)^T.
\ee
Here, by $b(\alpha):=b\left(a^{-1}(\alpha)\right)$ we denote the imaginary part of $\lambda$ as a function of
the new control parameter. Nonlinear function $h:\Re^3\times \Re\to \Re^3$  is a smooth function 
such that $h\left(0,\alpha\right)=0$ and ${\p h (0,\alpha)\over \p x}=0$
for values of $\alpha$ close to $0$. More precisely, $h$ stands for a family of functions parametrized by $p$
(see (\ref{6.4})). To keep the notation simple, we omit the dependence of $h$ on the second control parameter $p$.
\begin{figure}
\begin{center}
{\bf a}\epsfig{figure=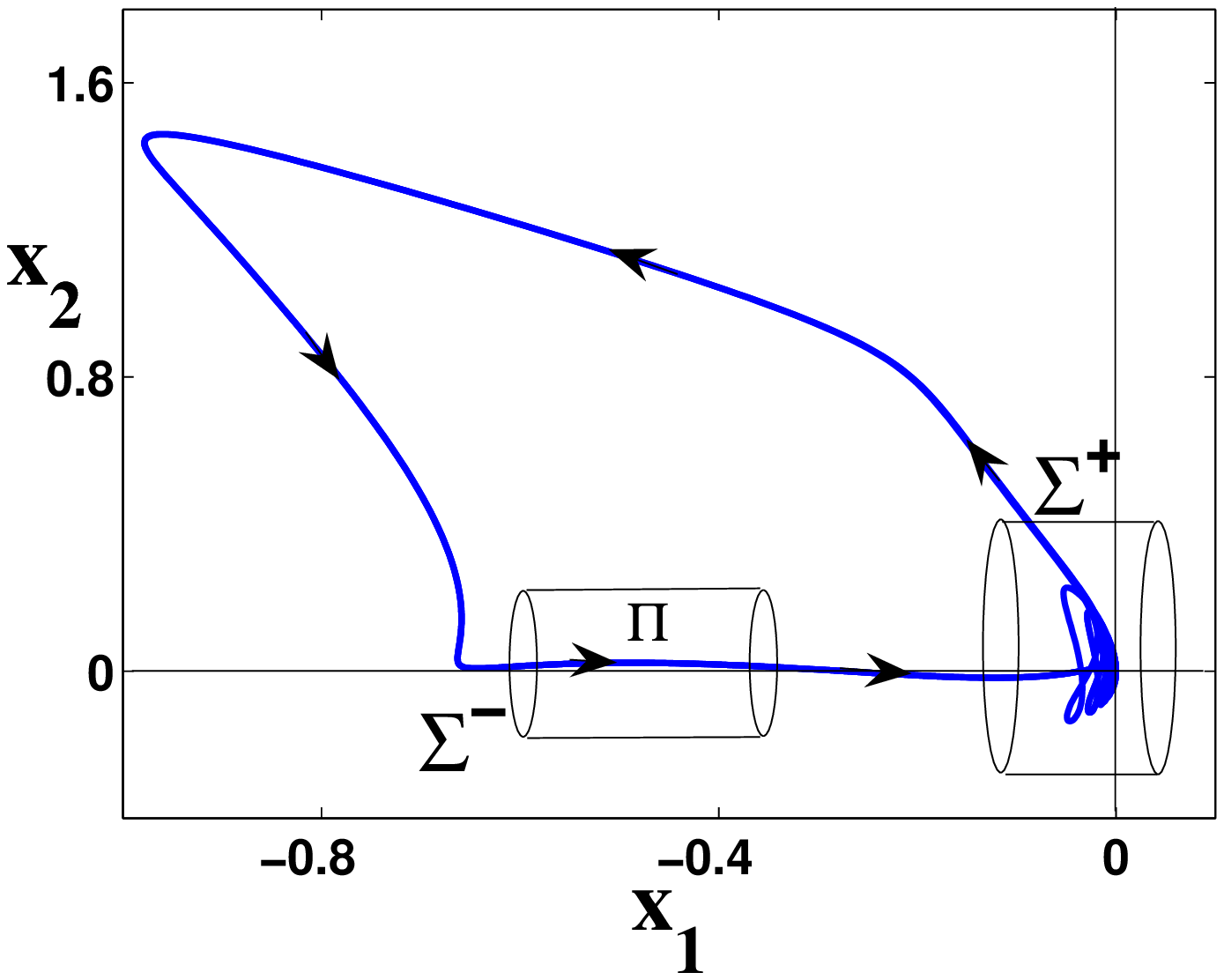, height=2.75in, width=3.0in, angle=0}
{\bf b}\epsfig{figure=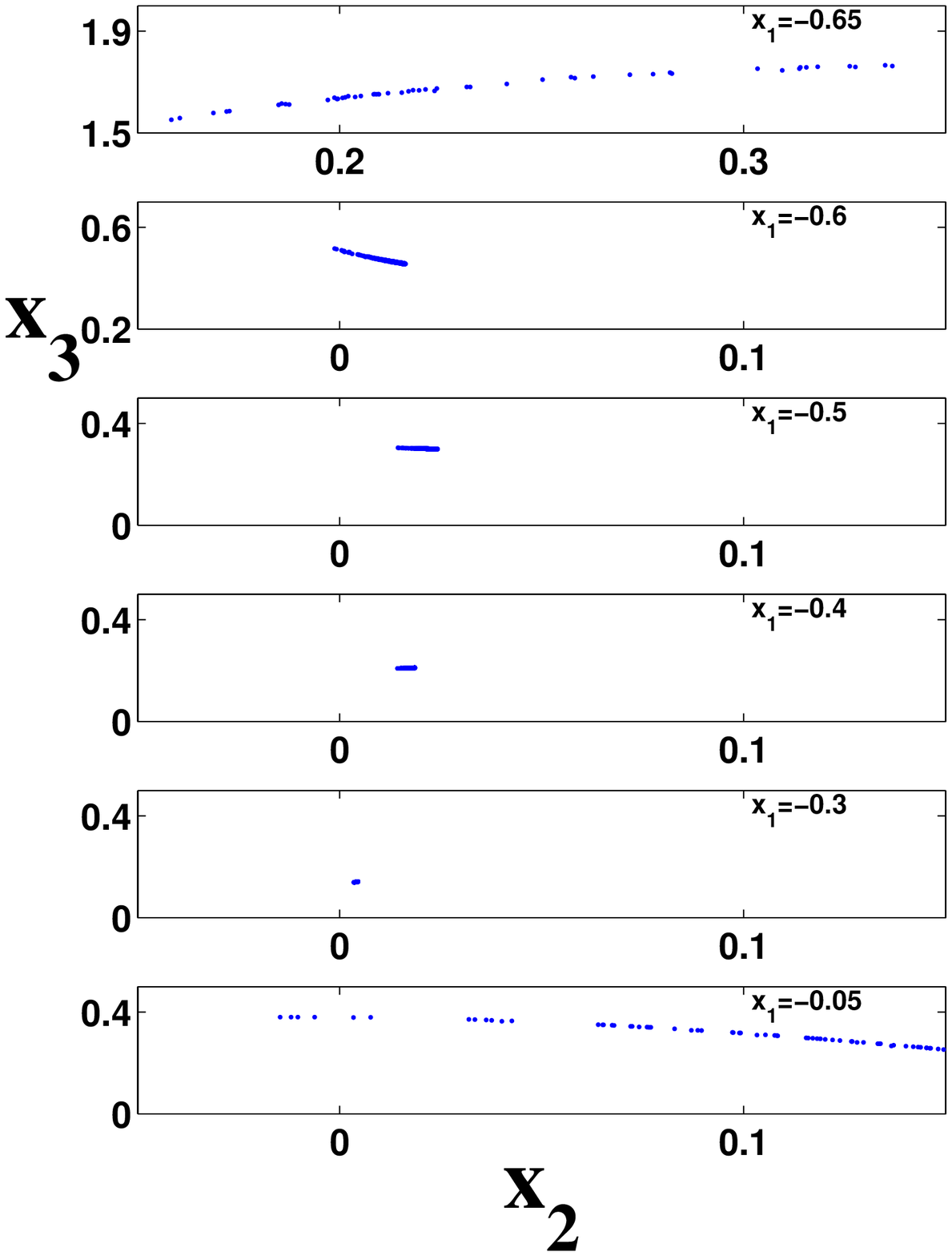, height=2.75in, width=3.0in, angle=0}
\end{center} 
\caption{{\bf a}. A periodic trajectory of the original system in the regime close to a subcritical AHB.
Cross-sections $\Sigma^+$ and $\Sigma^-$ are defined as follows:
$\Sigma^+=\left\{\left(x_1, \rho\cos\theta, \rho\sin\theta\right):\; \rho=0.1,\; x_1\in [-0.02, 0.02],\;
\theta\in [0,2\pi)\right\}$. $\Sigma^-$ is orthogonal to $x_1-$axis and is set at 
$x_1=-0.65$.
The numerical results shown in {\bf b} test the return mechanism and the strong contraction
property for (\ref{6.5}) (see {\bf (G1)} and {\bf (G2)} in the text). For this, we
covered $\Sigma^+$ with
a uniform mesh  $\Sigma_{N_i,N_j}^+=\left\{ \xi_{ij}\right\},$ $i=1,2,\dots,N_i,$ $j=1,2,\dots, N_j$. 
Using points $\xi_{ij}$ as initial conditions, we integrated
(\ref{6.5}) until the corresponding trajectories hit the second cross-section, $\Sigma^-$,
transverse to the stable manifold. The top plots in {\bf b} shows the image of $\Sigma_{N_i,N_j}$ 
in $\Sigma^-$ under the flow-defined map.
To test the strong contraction property, we repeated this experiment by taking 
$\Sigma^-$ progressively closer to the origin (the location of $\Sigma^-$ is 
indicated in each plot). The top five images in {\bf b} clearly indicate
that the projection of  the vector field onto $\Sigma^-$ is strongly contracting. The bottom plot
shows that the trajectories starting from $\Sigma^+$ reach a small neighborhood of the 
unstable equilibrium.
}\label{f.2}
\end{figure} 
It turns out that in the range of parameters of interest, the AHB of the equilibrium for $\alpha=0$ 
can be either subcritical, or 
supercritical, or degenerate depending on the value of $p$.  This gives rise to several 
qualitatively distinct oscillatory regimes 
generated by (\ref{6.5}) for small values of $\alpha> 0$.  The corresponding bifurcation scenarios 
for (\ref{6.1})-(\ref{6.3}) 
for fixed values of $p$ and varying $\nu$ were studied numerically in \cite{FKRT}. 
Below, we reproduce some of these numerics for the system in new coordinates and supplement 
them with a set of new numerical experiments relevant to the analysis 
of this paper. 
After that we state our results for each of the following cases: supercritical AHB,
subcritical AHB, and the transition layer between the regions of sub- and supercritical 
bifurcations. 

\begin{figure}
\begin{center}
{\bf a}\epsfig{figure=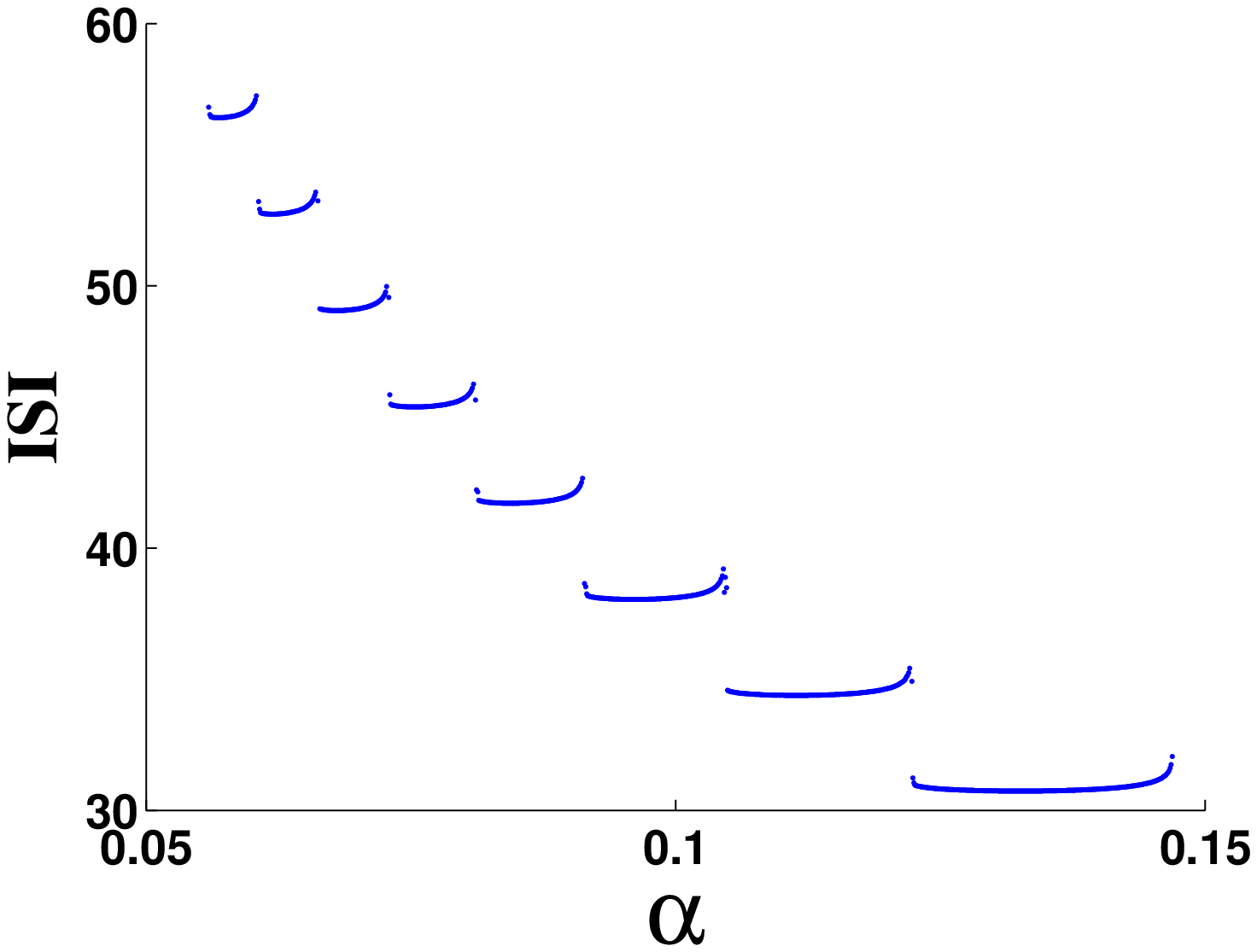, height=2.25in, width=3.0in, angle=0}
{\bf b}\epsfig{figure=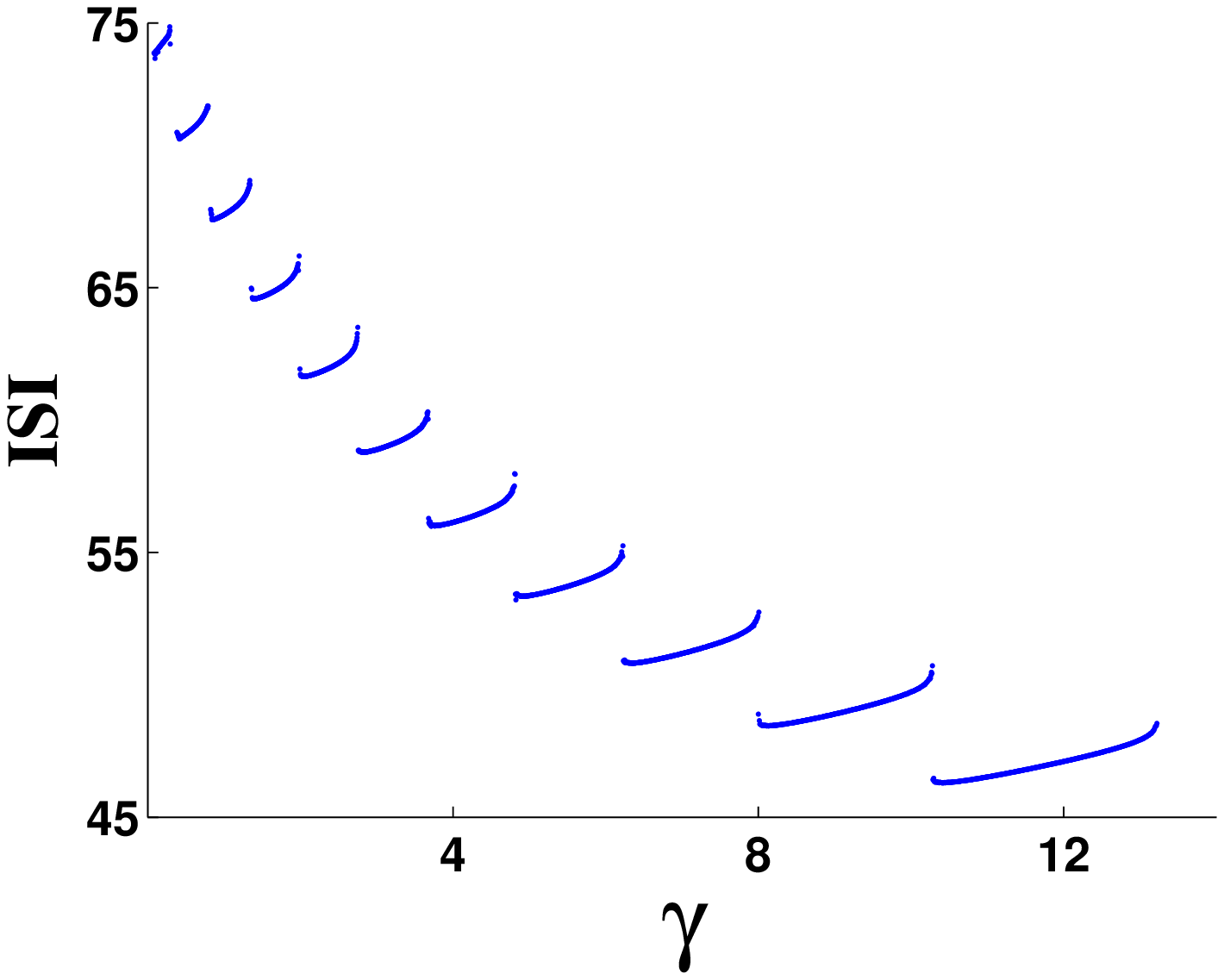, height=2.25in, width=3.0in, angle=0}
\end{center}
\caption{
The duration of the ISIs is plotted for:
({\bf a}) $\gamma=0.778$ and ({\bf b}) $\alpha=0.045$.
}\label{f.3}
\end{figure}

We start with discussing the supercitical case. It is well known that a supercritical AHB
produces a stable periodic orbit in a small neighborhood 
of the bifurcating equilibrium. The period of the nascent orbit is approximately 
equal to $2\pi\beta^{-1},$ where $\beta=b(0)$ is the imaginary part of the pair of complex conjugate  eigenvalues 
of the matrix of the linearized system at the bifurcation.  The numerical time series of $x_1$ and $x_{2,3}$ show 
that, while $x_{2,3}$ oscillate with the frequency prescribed by the Andronov-Hopf Bifurcation Theorem 
\cite{KUZ, MarsdenMcCracken}, the former oscillates 
with twice that frequency (see Figure \ref{f.1}a,b). 
The asymptotic analysis of Section 3 explains this counter-intuitive effect and shows that 
this is, in fact, a generic property of the periodic orbits born from a supercritical bifurcation in $R^n$, $n\ge3$.  
Note that the difference in frequencies of oscillations in $x_1$ and $x_{2,3}$ can not be understood from the 
topological normal form from the AHB \cite{GEO}, because the latter does not contain the information about 
the geometry of the periodic orbit. In Section 4.1,
we compute two geometric invariants of the periodic orbits as a curve in $\Re^3$: the curvature and the torsion.  
The latter shows that generically the orbit born from a supercritical AHB is not planar.  This together with certain 
symmetry properties of the orbit explains the frequency doubling of the oscillations in $x_1$.
As was noted in \cite{FKRT}, for increasing values of $\alpha > 0$, the small periodic orbit born from the supercritical 
AHB undergoes a cascade of period-doubling bifurcations, the first of which is shown in Figure \ref{f.1}c,d. Already 
at the first period-doubling bifurcation, the periodic orbit lies outside the 
region of validity of the local power series expansions, developed in the present paper.  
Therefore, our analysis does not explain the period-doubling
bifurcations for increasing values of $\alpha>0$.
In Section 4.3 we complement our analytical results with the numerical construction of the $1D$ first return map. 
The latter explains the mechanism for period-doubling cascade and the window of complex dynamics reported in 
\cite{FKRT}. 
\begin{figure}
\begin{center}
{\bf a}\epsfig{figure=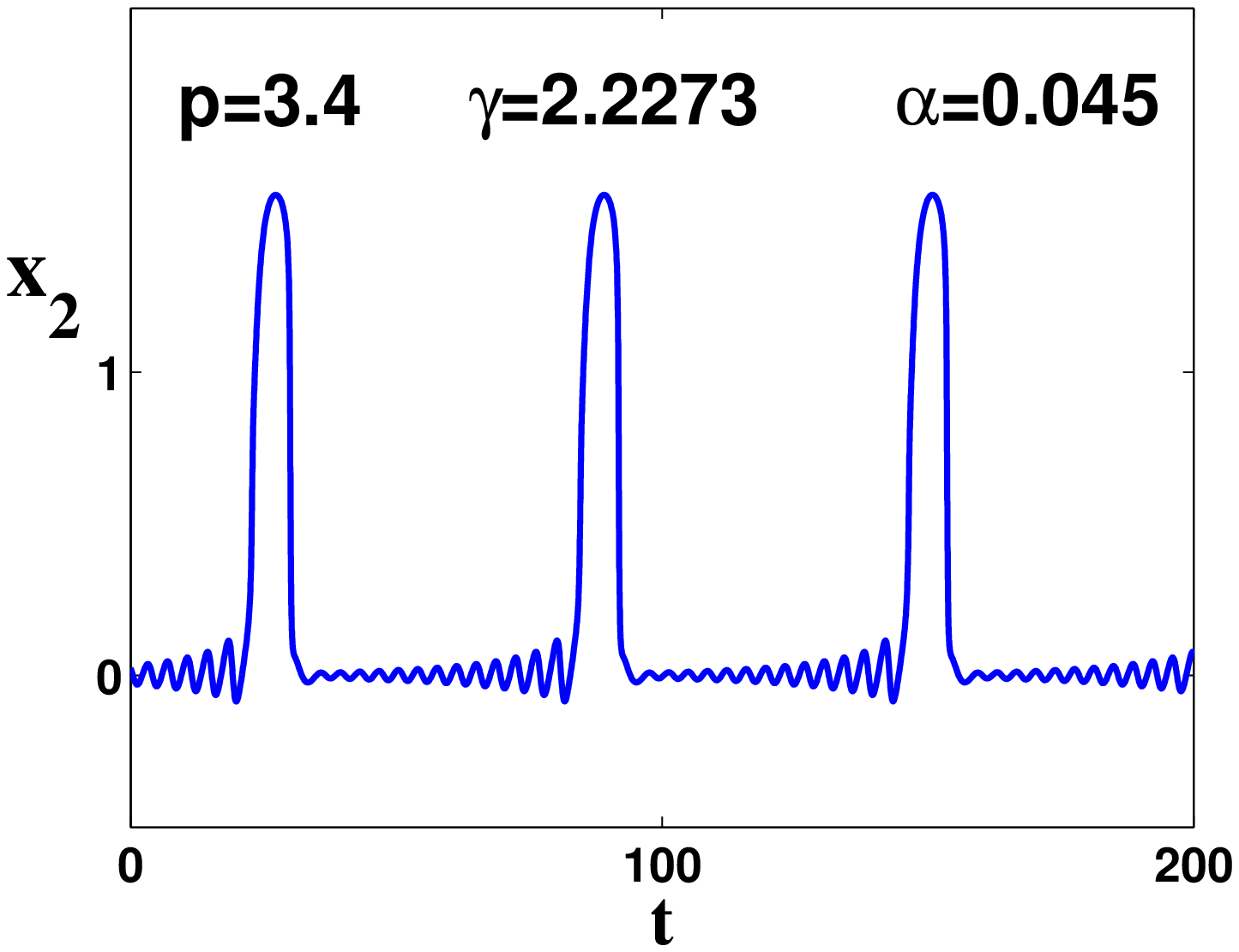, height=1.6in, width=3.0in, angle=0}\quad
{\bf b}\epsfig{figure=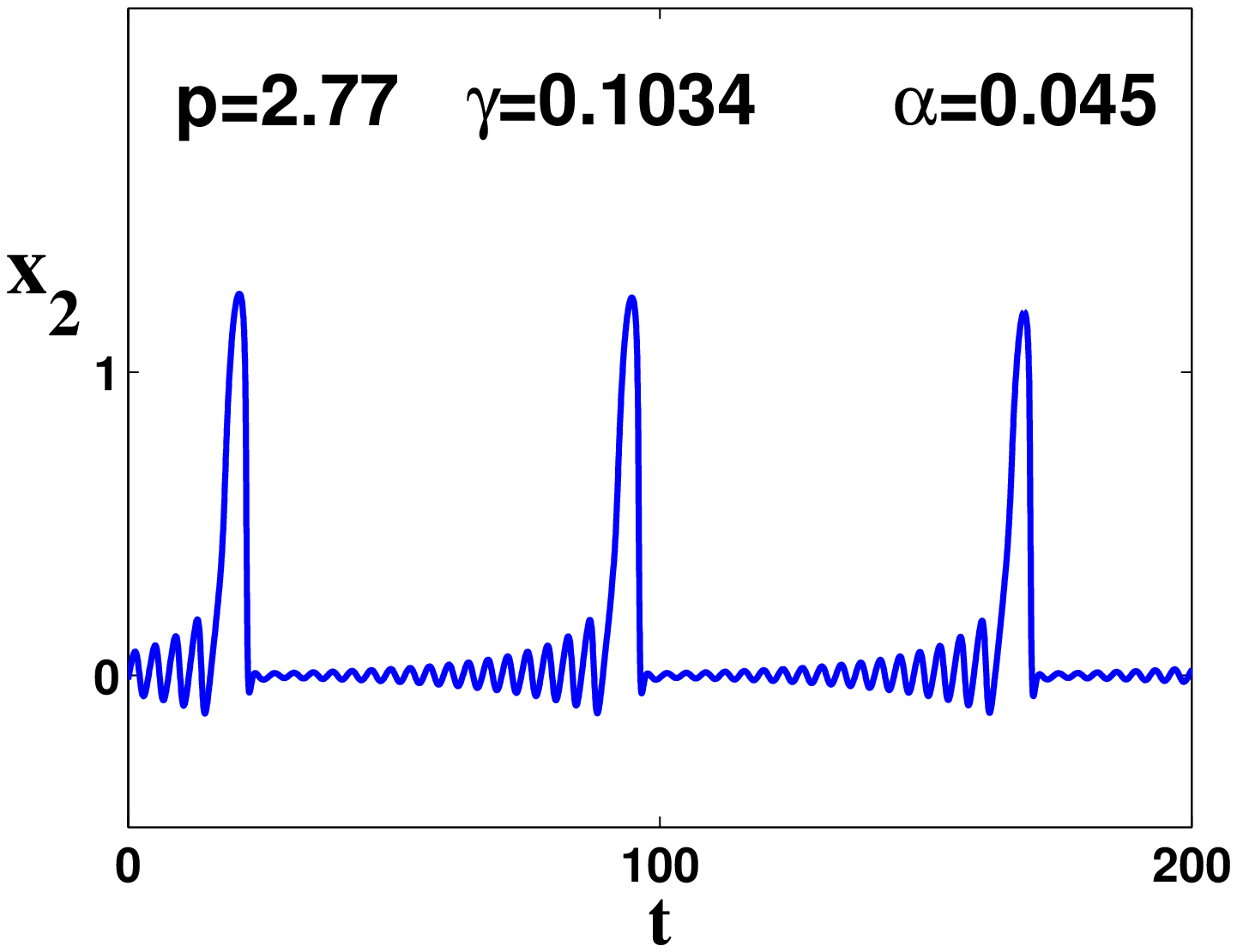, height=1.6in, width=3.0in, angle=0}
{\bf c}\epsfig{figure=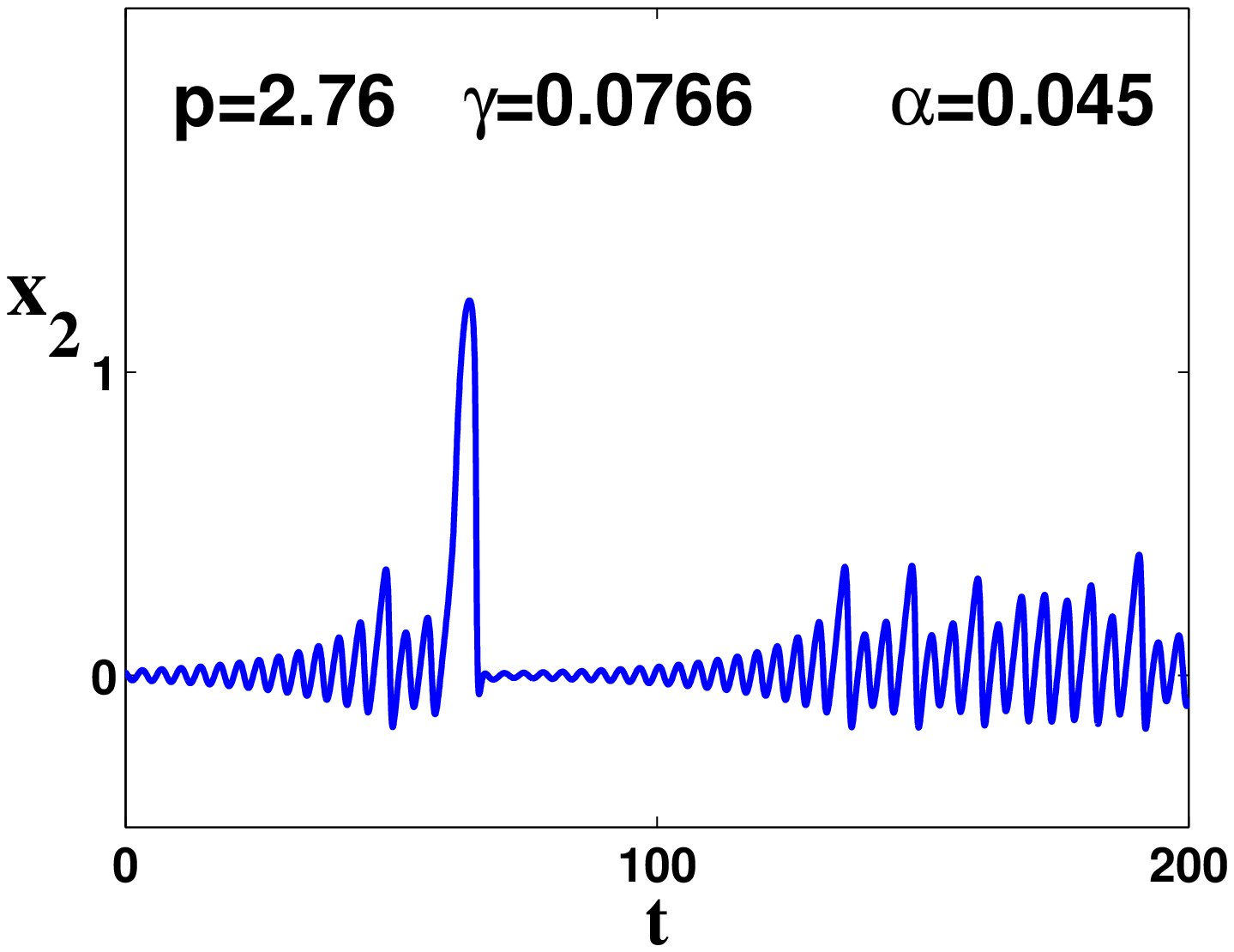, height=1.6in, width=3.0in, angle=0}
{\bf d}\epsfig{figure=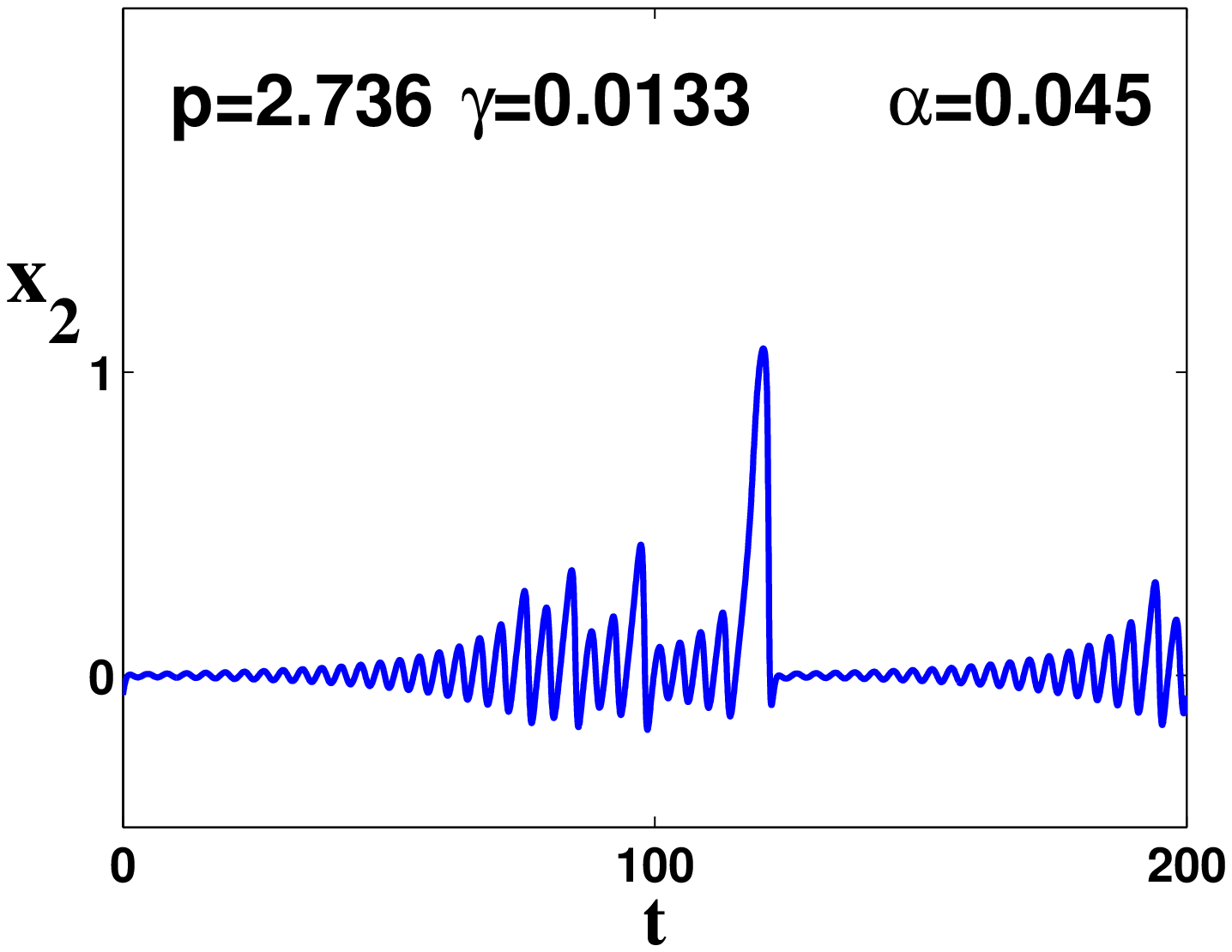, height=1.6in, width=3.0in, angle=0}\quad
\end{center}
\caption{The time series of $x_2$ are plotted for decreasing values of $\gamma>0$.
In plots {\bf a} and {\bf b}, the number of small amplitude oscillations separating the spikes increases,
while the qualitative form of the solutions remains the same.
For smaller values of $\gamma$, the system exhibits intermittency:
series of regular small amplitude oscillations are followed by those of irregular
oscillations of the intermediate amplitude ({\bf c} and {\bf d}). 
}\label{f.4}
\end{figure}
\begin{figure}
\begin{center}
{\bf a}\epsfig{figure=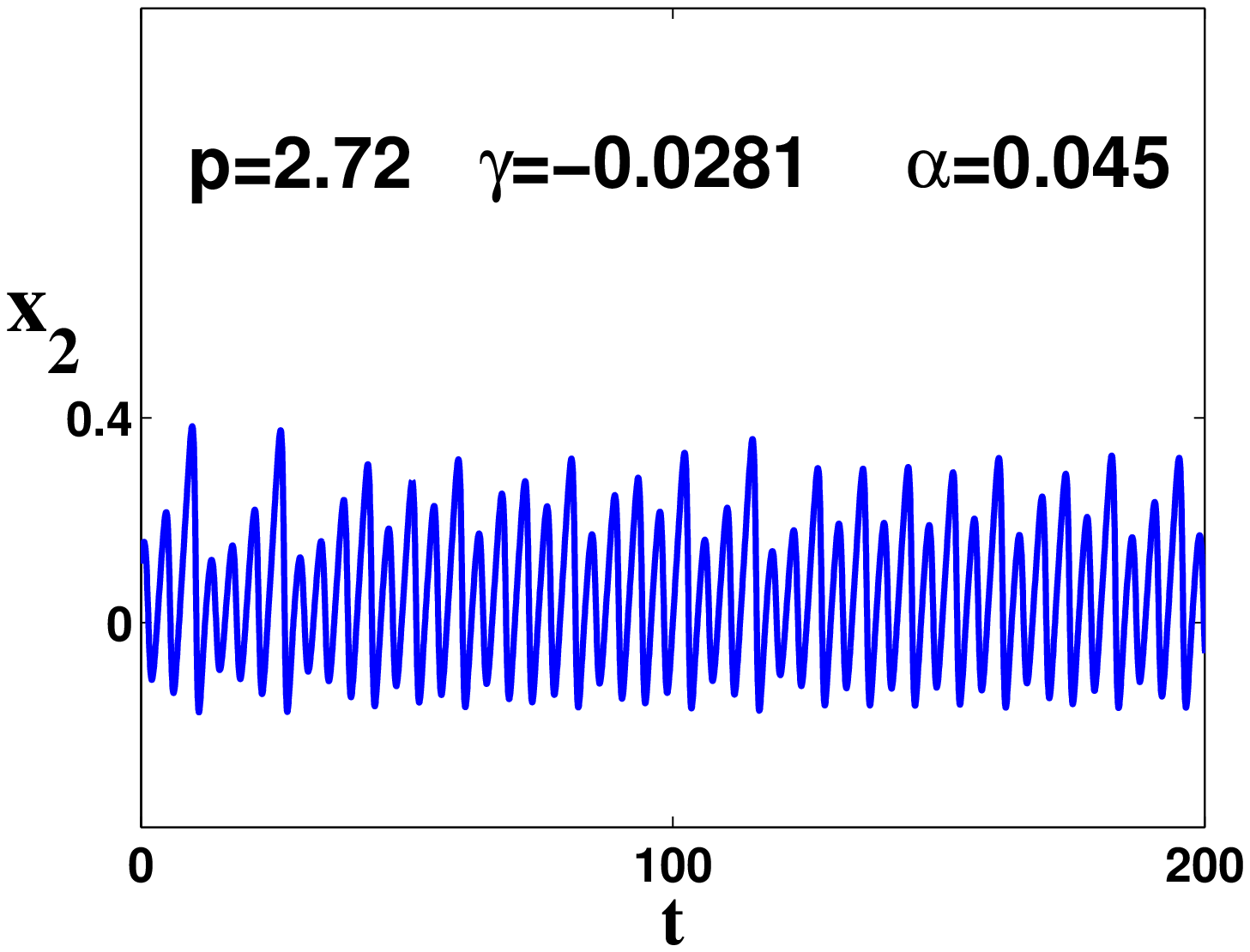, height=1.6in, width=3.0in, angle=0}\quad
{\bf b}\epsfig{figure=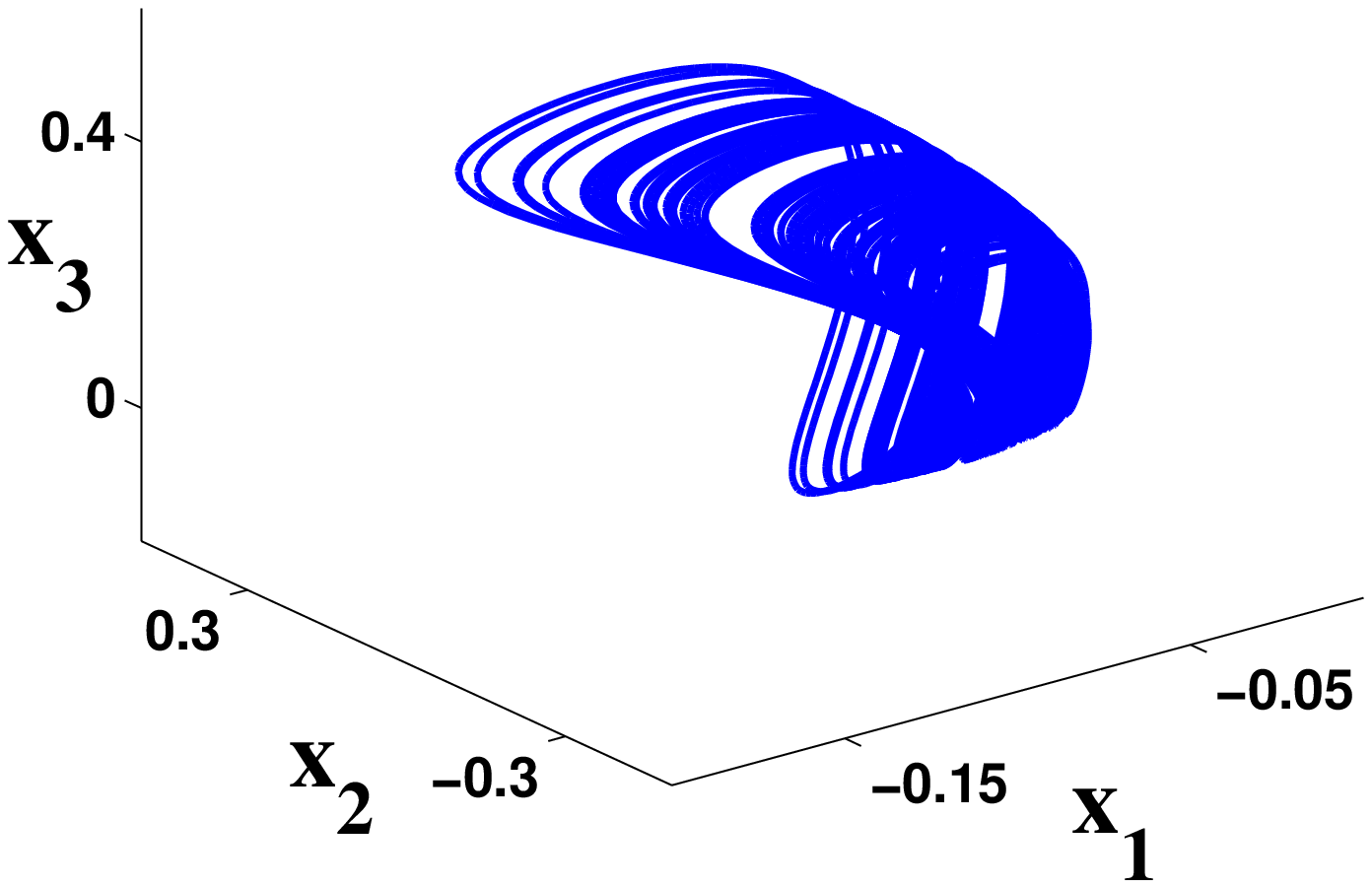, height=1.6in, width=3.0in, angle=0}
{\bf c}\epsfig{figure=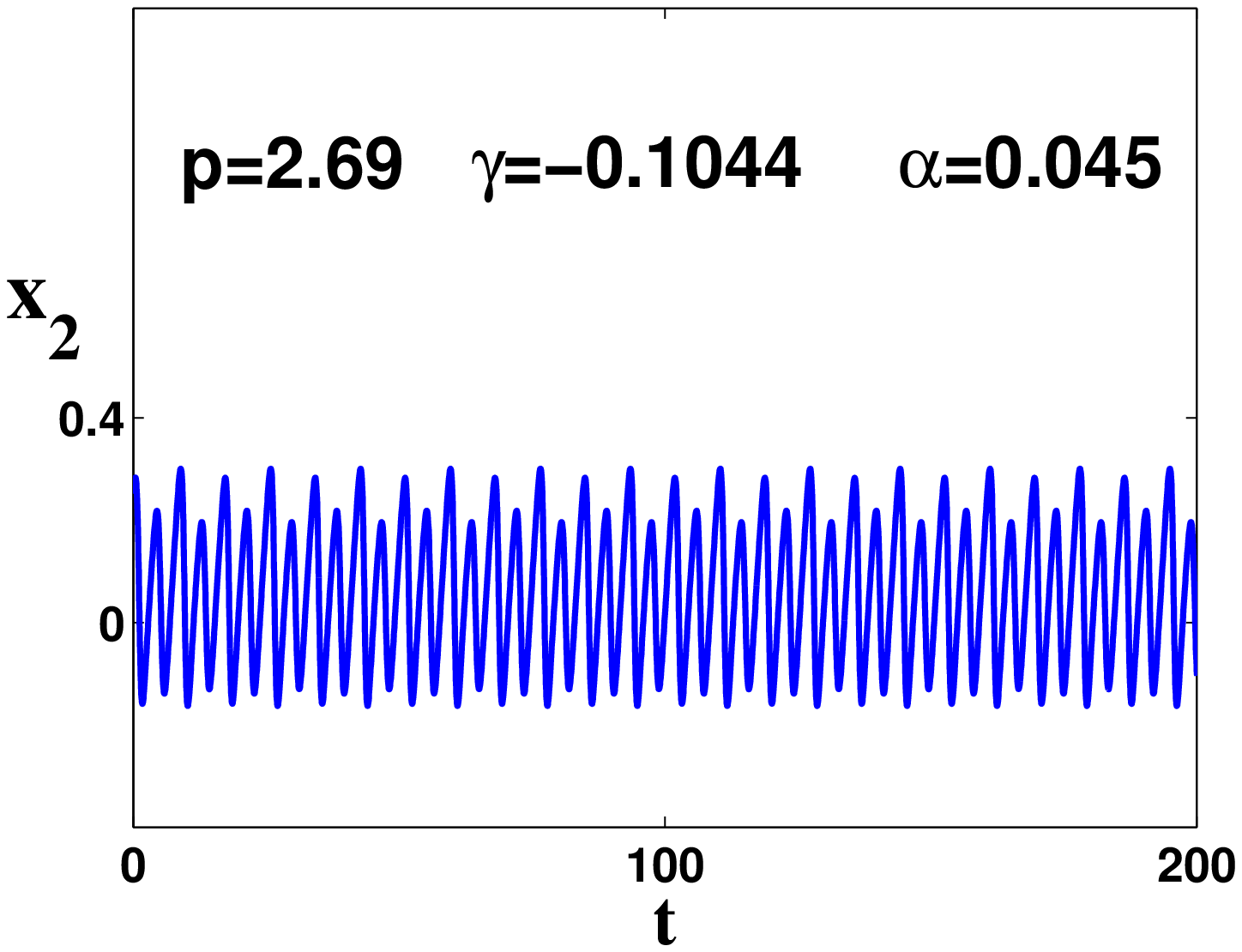, height=1.6in, width=3.0in, angle=0}\quad
{\bf d}\epsfig{figure=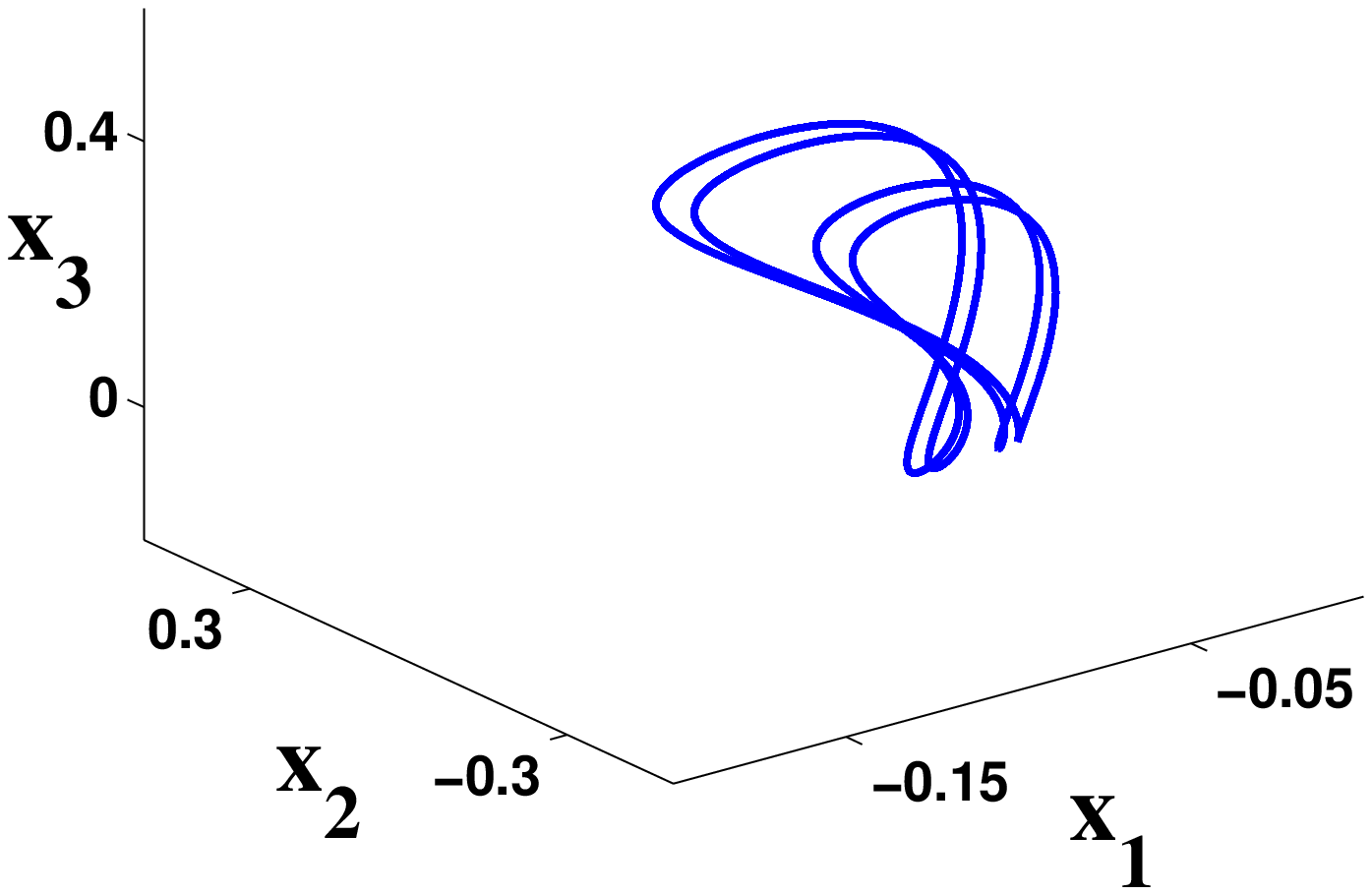, height=1.6in, width=3.0in, angle=0}
{\bf e}\epsfig{figure=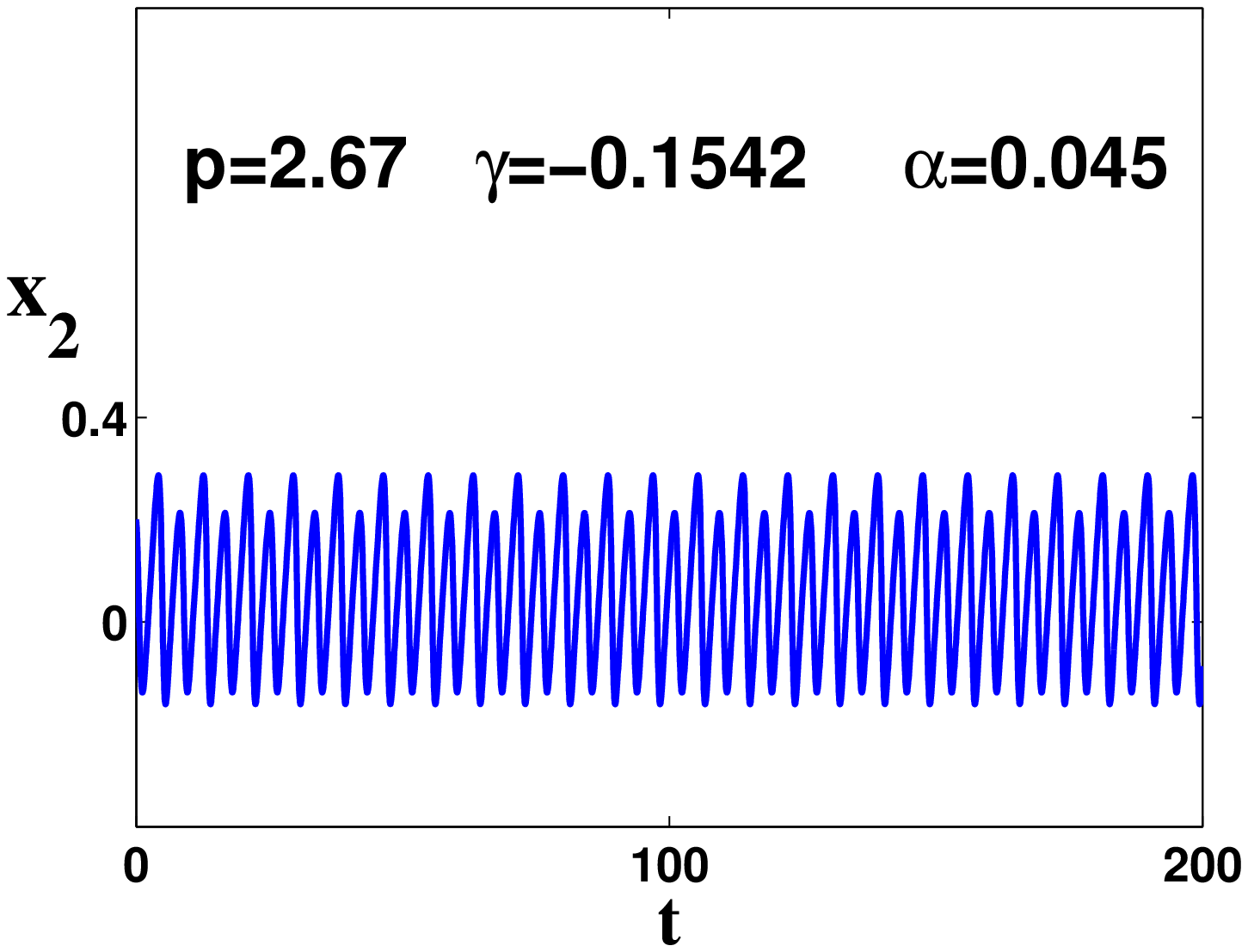, height=1.6in, width=3.0in, angle=0}\quad
{\bf f}\epsfig{figure=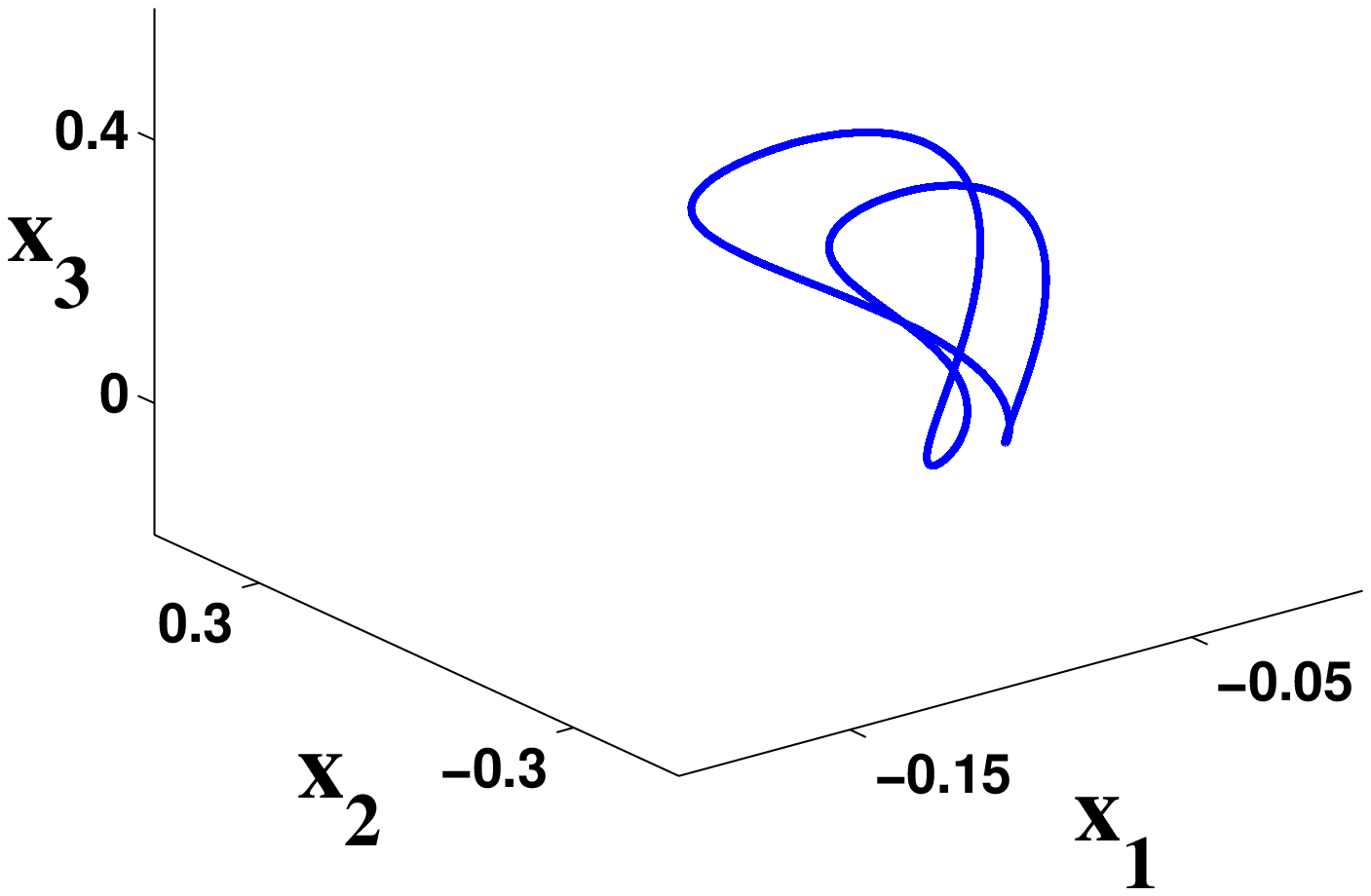, height=1.6in, width=3.0in, angle=0}
{\bf g}\epsfig{figure=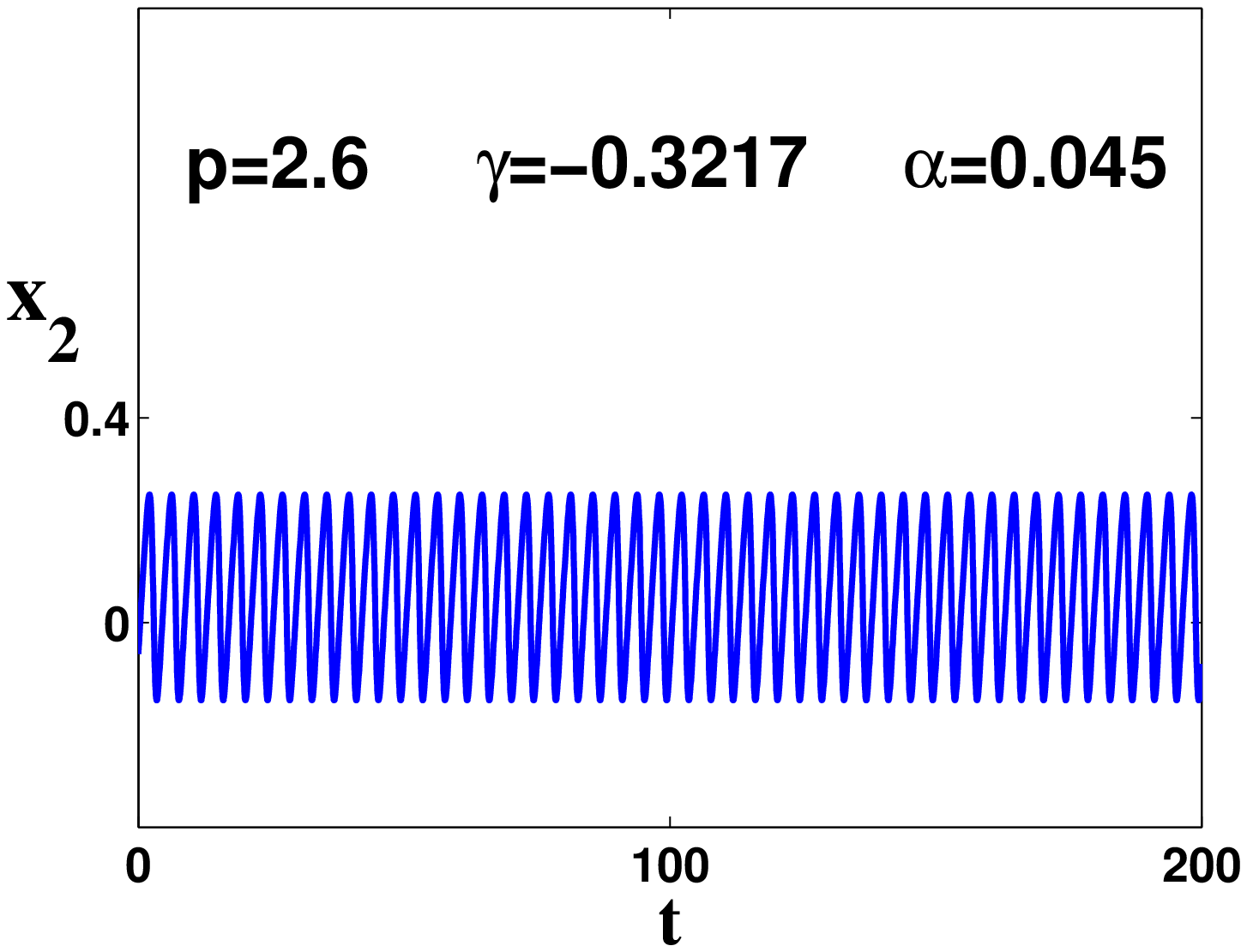, height=1.6in, width=3.0in, angle=0}\quad
{\bf h}\epsfig{figure=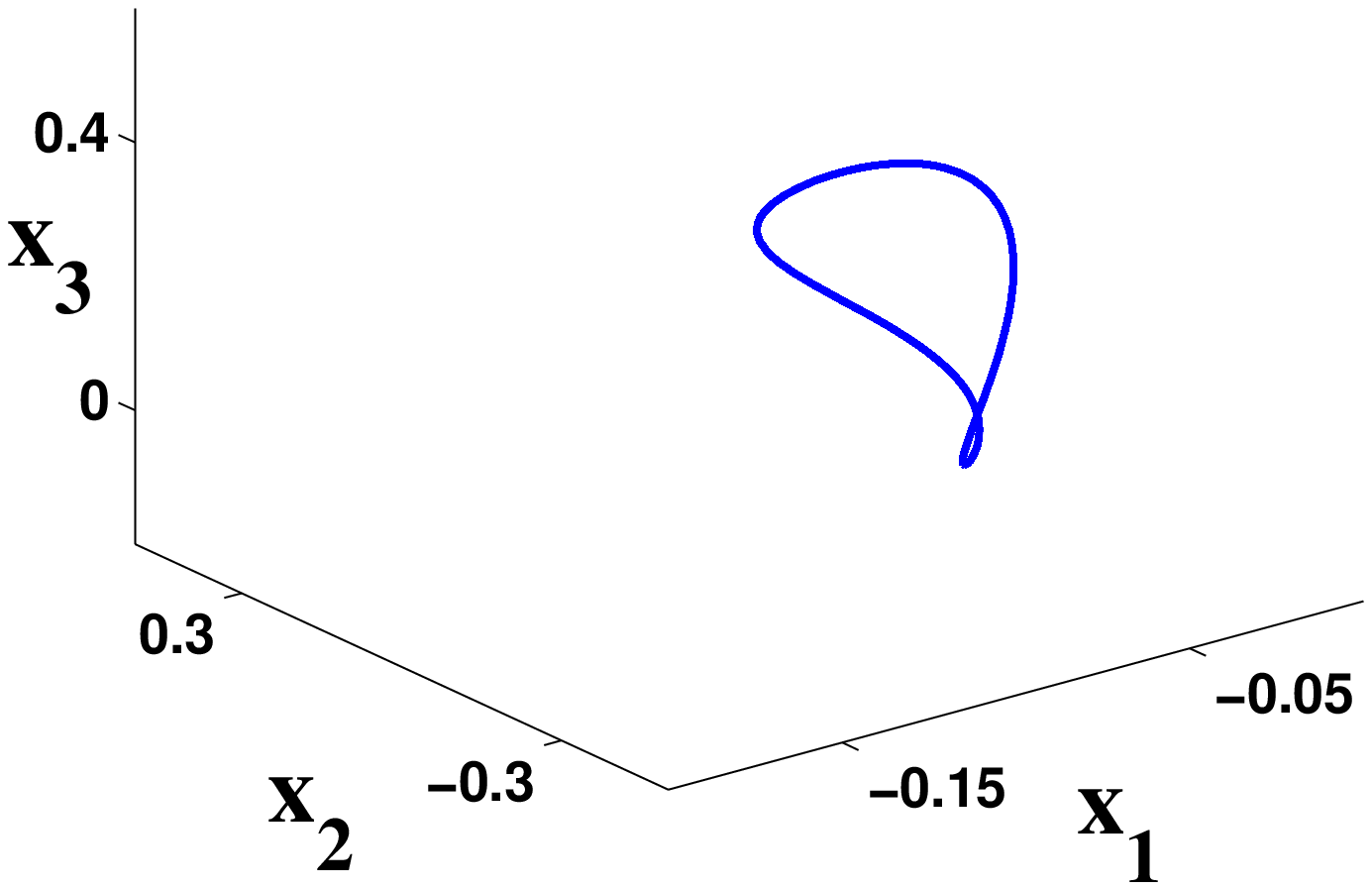, height=1.6in, width=3.0in, angle=0}
\end{center}
\caption{
The time series of $x_2$ (left column) and the corresponding trajectories of (\ref{6.1})-(\ref{6.3})
(right column) plotted for negative values of $\gamma$ close to $0$. The transition from the 
complex dynamics in {\bf a,b} to regular approximately harmonic oscillations in {\bf g,h}
contains a reverse cascade of period-doubling bifurcations ({\bf c}-{\bf f}).
}\label{f.5}
\end{figure}
 	
We next turn to the subcritical case.  
The dynamics resulting from the subcritical AHB depends on the properties of the vector field outside 
of a small neighborhood of the equilibrium. We conducted a series of numerical experiments to study
the global properties of the vector field (see caption of Figure \ref{f.2} for details). Based 
on these numerical observations, we identify two essential features of the vector field:
{\bf (G1)} the return mechanism and {\bf (G2)} the strong contraction property. Specifically,
\begin{description}
\item[(G1)]
For a suitably chosen cylindrical crossection, $\Sigma^+$, placed sufficiently close to the origin, and
another crossection, $\Sigma^-$, transverse to $W^s(O)$ (see Figure \ref{f.2}a), the flow-defined map
$Q:\Sigma^+\to\Sigma^-$ is well-defined for $\alpha\ge 0$ and depends smoothly on the parameters of the
system.
\item[(G2)] There is a region $\Pi$, adjacent to $\Sigma^-$ and containing a subset of $W^s(O)$
(see Figure \ref{f.2}a), in which the projection of the vector field in direction transverse to 
$W^s(O)$ is sufficiently stronger than that in the tangential direction.
\end{description} 
In Section 4.2, conditions {\bf (G1)} and {\bf (G2)} are made precise.
Properties {\bf (G1)} and {\bf (G2)} guarantee that for small values of $\alpha>0$, the trajectories leaving
a small neighborhood of the saddle-focus, after some relatively short time enter $\Pi$. In $\Pi$,
the trajectories approach $W^s(O)$ closely and then follow it to a sufficiently small neighborhood of the 
unstable equilibrium. For small $\alpha>0$, the trajectories starting from a sufficiently small neighborhood
of the saddle-focus, remain in some larger (but still small) neighborhood of the origin for a long time.
Eventually, they hit $\Sigma^+$ and the dynamics described above repeats. Therefore, under conditions
{\bf (G1)} and {\bf (G2)}, the
subcritical AHB bifurcation results in a 
sustained motion consisting of long intervals of time that the trajectory spends near a weakly unstable saddle-focus 
and relatively brief excursions outside of a small neighborhood of the origin (see Figure \ref{f.2}a).  
In terms of the time series, the dynamical variables undergo series of small amplitude oscillations 
alternating with large spikes (Figure \ref{f.1}e-h). For increasing values of $\alpha > 0$, the number of small 
amplitude oscillations decreases and so do the time intervals between consecutive spikes, so-called 
interspike intervals (ISIs).  
The multimodal solutions arising near the subcritical AHB are not necessarily periodic. Nonetheless, 
the ISIs may be characterized in terms of the control parameters present in the system, regardless of whether the
underlying trajectories 
are periodic or not.  Specifically, in Section 4.2, we derive an asymptotic relation for the 
duration of the ISIs: 
\be\lbl{6.6}
\tau \sim {1\over 2\alpha} \ln\left(1+ {\alpha\over\gamma(p)}O\left(\epsilon^{-\varsigma}\right)\right),\;\;
\varsigma\ge 4,
\ee
where the first Lyapunov coefficient $\gamma(p)$ reflects the distance of the system from the transition 
from sub- to supercritical AHB, with $\gamma(p)$ 
being positive when the AHB is subcritical and equal to zero at the transition point. 
Finally, $\epsilon>0$ is a small parameter and positive $\alpha=O(\epsilon^2)$.
The role of $\epsilon$ will become clear later.  We illustrate (\ref{6.6}) with 
numerically computed plots of the ISIs under the variation of control parameters in Figure \ref{f.3}.  
The ISIs provide a 
convenient and important characteristics of the oscillatory patterns involving pronounced spikes.  For example, it 
is widely used in both theoretical and experimental neuroscience for description of patterns of electric activity 
in neural cells.  An important aspect of the ISIs, is that the dependence of the ISIs on the control parameters 
in the system (which often can be directly established experimentally) reveals the bifurcation structure of the system.  
Therefore, the analytical characterization of the ISIs in terms of the bifurcation parameters, 
such as given in (\ref{6.6}), is important in applications. For the problem at hand, 
the ISIs depend on the interplay of the two parameters 
$\alpha$ and $\gamma(p)$, which reflect the proximity of the system to a codimension 2 bifurcation.  
The dependence of the ISIs on $\alpha$ was studied in \cite{GW} 
for a model problem near a Hopf-homoclinic bifurcation.  Our analysis extends the formula 
for the ISIs obtained in \cite{GW} to a wide class of systems and emphasizes the role 
of $\gamma(p)$, which is important for the present problem.  
We also note that the structure of the global vector field suggests that (\ref{6.5}) is 
 also close to  a homoclinic bifurcation, which can also influence the ISIs.  However, under 
the variation of the control parameters $\alpha$ and $p$ the system remains 
bounded away from the homoclinic bifurcation, so that the latter effectively does not affect the ISIs. 
In fact, (\ref{6.6}) can be easily extended to include the distance from the homoclinic bifurcation as well.  
As follows from (\ref{6.6}), the ISIs increase for decreasing values of $\gamma>0$ (see Figure \ref{f.3}b).  
This observation prompted our interest in investigating the transition in the system 
dynamics as $\gamma$ crosses $0$. Numerical simulations show that as $\gamma$ approaches $0$ from the positive side, the 
oscillations become very irregular and very likely to be chaotic (Figure \ref{f.4}).  
When $\alpha >0$ is fixed and $\gamma$ crosses 0, 
the irregular dynamics is followed by the reverse period-doubling cascade and terminates with the creation of a 
regular small amplitude periodic orbit, which can be then tracked 
down to a nondegenerate supercritical AHB (Figure \ref{f.5}a-h).
This scenario presents a certain interest as it suggests the formation of 
the chaotic attractor in a well-defined bifurcation setting.  Namely, the transition 
from the region of subcritical to supercritical AHB for small $\alpha >0$ leads 
to the formation of a chaotic attractor and a reverse period-doubling cascade.
As in the case of the complex dynamics appearing for 
increasing values of $\alpha$ for supercritical AHB, the irregular oscillations in the present 
case can not be understood using the local analysis alone.  Therefore, 
we complemented our analysis with the study of numerically constructed first return map.  
The latter gives a clear geometric picture of the origins of the chaotic dynamics and the 
period-doubling cascade of periodic orbits in this parameter regime.

\section{The local expansions}
\setcounter{equation}{0}
In the present section, we study trajectories of (\ref{6.5}) in a small neighborhood of the 
equilibrium at the origin. Assume that $h(x,\alpha)$ is a smooth function in a small neighborhood
of $(0,0)$, so that all power series expansions below are justified. 
We are interested in the dynamics of (\ref{6.5}) for values of $\alpha$ close to $0$. 
The role of $\alpha$ in our analysis is twofold. On the one hand, $\alpha$ is a control
parameter in this problem, on the other hand, the smallness of $\alpha$ is used
in the asymptotic analysis below. To separate these two roles, we use the following
rescaling
\be\lbl{5.1a}
\alpha=\mu\epsilon^2,
\ee
where $\epsilon>0$ is a fixed sufficiently small constant and $\mu$ varies in a certain interval
of size $O(1)$. We introduce cylindrical coordinates in $\Re^3$
\be\lbl{cl}
\left(x_1,x_2,x_3\right)\mapsto\left(x_1, \rho,\theta\right)\in\Re\times\Re^+\times S^1,\quad
x_2=\rho\cos\theta,\;x_3=\rho\sin\theta,
\ee
and define
\be\lbl{5.2}
D=\left\{ \left(x_1,\rho, \theta \right):\; \left| x_1 \right| \le M^2\epsilon^2, \rho\in(0,2M\epsilon]\right\},\;
D_0=\left\{ \left(x_1,\rho, \theta \right):\; \left| x_1 \right| \le M^2\epsilon^2, \rho\in(0,M\epsilon]\right\},
\ee
where $M>0$ is sufficiently large.

To characterize the trajectories of (\ref{6.5}) in $D$ we will use an exponentially stable slow manifold, $S$,
whose leading order approximation is given by
\be\lbl{5.3}
S_0= \left\{ 
\left( x_1,\rho,\theta\right):\; x_1=U(\theta)\rho^2,\; 0\le\rho\le 2\epsilon M\right\},
\ee
where
\be\lbl{5.4}
U(\theta)= a + A\cos\left(2\theta-\vartheta\right).
\ee
and $a, A$ and $\vartheta$ are computable constants (see Appendix A).

Below we show that  the trajectories of (\ref{6.5}) with initial data from $D_0$ approach an
$O(\epsilon)$  neighborhood of $S_0$ in time $O\left(\left| \ln \epsilon \right|\right)$
and stay in this neighborhood as long as they remain in $D$. The reduction of the system dynamics to $S$
yields a complete description of the trajectories of (\ref{6.5}) in $D$. The qualitative character of the
solution behavior in $D$ depends on the sign of {\it the first Lyapunov coefficient}
\be\lbl{5.5}
\gamma={1\over 2\pi}\int_0^{2\pi} \left( U\left(\theta\right) \bar Q_2\left(\theta\right)+
\bar Q_3\left(\theta\right)\right) d\theta,
\ee
where $2\pi-$periodic trigonometric polynomial $\bar Q_{2,3}$ are given in Appendix A.
Finally, by $\beta=b(0)$ we denote the absolute value of the imaginary parts of the
complex conjugate pair of eigenvalues at the AHB (see (\ref{6.5})).
The results of the analysis of this section are summarized in

\noindent
{\bf Theorem 3.1} 
{\sc
Suppose $\beta>0$ and $\gamma\ne 0$. Let $\epsilon>0$ denote a small parameter.
Then for $\alpha=O(\epsilon^2)$ and for sufficiently large (independent of $\epsilon$)
$M>0$, the trajectories
with initial conditions in $D_0$ enter an $O(\epsilon)$ neighborhood of $S_0$ in time
$O\left(\left|\ln\epsilon\right|\right)$ and remain there as long as they stay in $D$.
In the neighborhood of $S_0$, the trajectories can be uniformly approximated on any
interval of time $[t_0,\bar t]$ of length $O(1)$:
\be\lbl{5.7}
x_1(t)=\rho^2 U(\theta)+O(\epsilon^3), \quad x_2(t)=\rho\cos\theta +O(\epsilon^2),\quad\mbox{and}\quad
x_3(t)=\rho\sin\theta +O(\epsilon^2),
\ee
where $U(\theta)$ and $\gamma$ are defined in (\ref{5.4}) and (\ref{5.5}), respectively,  and
$$
\rho(t) = \left(\left(\rho(t_0)^{-2}+{\gamma\over\alpha}\right) e^{-2 \alpha(t-t_0)}-
{\gamma\over\alpha}\right)^{-1\over 2}+O(\epsilon^2),\quad
\dot\theta=\beta,\quad
\tan\left(\theta(0)\right)= {x_3(t_0)\over x_2(t_0)}.
$$
} 

\noindent {\bf Remark 3.1}
\begin{description}
\item[a)] From Theorem 3.1, one can easily deduce the existence of a periodic orbit $O_\alpha$,
when $\alpha\gamma<0$. $O_\alpha$ is stable  (unstable) if $\alpha>0$ ($\alpha<0$). The leading order
approximation for $O_\alpha$ follows from (\ref{5.7}):
\be\lbl{5.8}
x_1=\bar\rho^2 U(\theta)+O(\epsilon^3),\;\; x_2=\bar\rho\cos\theta +O(\epsilon^2),\;\;
x_3=\bar\rho\sin\theta +O(\epsilon^2),\;\;
\mbox{and}\;\;
\theta\in [0, 2\pi),
\ee
where $\bar\rho=\sqrt{\alpha\over -\gamma}$. Equations (\ref{5.4}) and (\ref{5.8})
imply that the frequency of oscillations in $x_1$ is twice
that of oscillations in $x_{2,3}$, provided $A\ne 0$ in (\ref{5.4}) 
(see Figure \ref{f.1}a,b). In Section 4.1, we give a geometric interpretation of this frequency doubling 
effect.
\item[b)] If $\alpha>0$ and $\gamma>0$, Theorem 3.1 describes the trajectories with initial conditions
near a weakly unstable equilibrium. It shows that along these trajectories, $x_{2,3}$ undergo  approximately
harmonic oscillations, whose amplitude grows at the rate $O(\epsilon^2)$. In addition, they satisfy the
following scaling relation:
\be\lbl{5.9}
\bar x_1(t) \sim a\rho^2(t),
\ee
where $\bar x_1(t)$ stands for the average value of $x_1(t)$ over one cycle of oscillations.
\item[c)]
The domain of validity of the asymptotic expansions in (\ref{5.7}) extends much farther than
$D$. As will follow from the proof of Theorem 3.1, $M$ in the definition of $D$ may be taken 
up to $o\left(\epsilon^{-1\over 3}\right)$. This
implies that, with the error term $o(1)$, (\ref{5.7}) remains valid for the ranges of $x_1$ and $x_{2,3}$ up to 
$o(\epsilon^{4\over 3})$ and $o(\epsilon^{2\over 3})$, respectively.
\end{description}

In the remainder of this section, we prove Theorem 3.1. By expanding
$h(x,\alpha)$ into finite Taylor sum with the reminder term, we have  
\begin{eqnarray}\nonumber
\dot x_1 &=& -x_1+\sum a_{ij}(\alpha)x_i x_j +\sum a_{ijk}(\alpha) x_ix_jx_k + O(4),\\ \lbl{5.1}
\dot x_2 &=& a(\alpha) x_2 -b(\alpha) x_3 +\sum b_{ij}(\alpha)x_i x_j +\sum b_{ijk}(\alpha) x_ix_jx_k +O(4),\\ \nonumber
\dot x_3 &=& b(\alpha) x_2 +a(\alpha) x_3 +\sum c_{ij}(\alpha)x_i x_j +\sum c_{ijk}(\alpha) x_ix_jx_k +O(4),
\end{eqnarray}
where
\begin{eqnarray*}
\sum \sigma_{ij}(\alpha) x_ix_j &=& \sum_{i=1}^3 \sum_{j=i}^3 \sigma_{ij}(\alpha) x_ix_j,\\
\sum \sigma_{ijk}(\alpha) x_ix_jx_k &=& \sum_{i=1}^3 \sum_{j=i}^3 \sum_{k=j}^3 \sigma_{ijk}(\alpha) x_ix_jx_k,\quad
\sigma\in\left\{a,b,c\right\}.
\end{eqnarray*}
To study system of equations (\ref{5.1}) in a small neighborhood of the origin, 
we rewrite it in cylindrical coordinates (\ref{cl}) and rescale variables
\be\lbl{5.10}
 x_1=\xi\epsilon^2\;\mbox{and}\; \rho=r\epsilon.
\ee
In new coordinates, we have
\begin{eqnarray}
\dot \xi &=&-\xi+ r^2 T_1(\theta) +O\left(\epsilon\right),\nonumber\\ \lbl{5.11}
\dot r &=& \epsilon r^2 R_1(\theta)+
\epsilon^2 r\left(\mu+\xi R_2(\theta)+r^2R_3(\theta)\right)+O(\epsilon^3),\\
\nonumber
\dot\theta &=& \beta+\epsilon r L_1(\theta)+ O\left(\epsilon^2\right),
\end{eqnarray}
where $T_1, L_1$ and $R_{1,2,3}$ are homogeneous trigonometric polynomials.

Using $\phi=\beta^{-1}\theta$ as a new variable, from (\ref{5.11}), we obtain
\begin{eqnarray}\lbl{5.11a}
\frac{d\xi}{d\phi}&=&-\xi+r^2 P_1(\phi)+\epsilon \Xi\left(\xi,r,\phi,\epsilon\right),\\ 
\lbl{5.12}
\frac{dr}{d\phi}&=&\epsilon Q_1(\phi) r^2 +\epsilon^2 r\left(\mu+\xi Q_2(\phi)+r^2 Q_3(\phi)\right)+
\epsilon^3 \Psi\left(\xi,r,\phi,\epsilon\right),
\end{eqnarray}
where $Q_{1,2,3}(\phi)$ and $P_1(\phi)$ are trigonometric polynomials of period $\omega=2\pi\beta^{-1}$.
Functions $\Xi$ and $\Psi$ are bounded and continuous in $D'\times\Re\times [0, 1]$.
Here $D'$ stands for the domain of the rescaled variable $(\xi,r)$ when the original variable
$x\in D\subset \Re^3$:
\be\lbl{Dprime} 
D^\prime= \left\{ \left( \xi,r \right):\; \left| \xi \right| \le 4M^2, \; 0<r\le 2 M \right\}.
\ee
Similarly, we define
\be\lbl{D0prime}
D^\prime_0=\left\{ \left(\xi,r\right):\; \left| \xi \right| \le M^2, \; 0<r\le M\right\}.
\ee

To determine the slow manifold for (\ref{5.1}), we introduce the following notation. 
By $\bar p_1$ and $\tilde P_1(\phi)$ we denote the mean value and the oscillating part 
of the trigonometric polynomial $P_1(\phi)$ in (\ref{5.11a}), respectively:
\be\lbl{5.13}
\bar p_1={\beta\over 2\pi}\int_0^{2\pi\over\beta} P_1(\phi)d\phi
\quad\mbox{and}\quad \tilde P_1(\phi)= P_1(\phi)-\bar p_1,
\ee
and define
\be\lbl{5.16}
\xi_0(\phi)=\bar p_1 +p_1(\phi) \quad\mbox{and}\quad
p_1(\phi)=
{\tilde P_1(\phi)-{\tilde P_1}^\prime (\phi)\over 1+4\beta^2}.
\ee

We first prove an auxiliary lemma, which characterizes the trajectories of (\ref{5.11a}) and (\ref{5.12})
in $D^\prime$.
In the following two lemmas, we will keep track of how the remainder terms in the asymptotic expansions
depend on the size of $D^\prime$. This will be used later to determine the domain of validity for the 
asymptotic approximation of the slow manifold. Below $\kappa=O\left(M\epsilon\right)$ means that there
exist positive constants $\epsilon_0$ and $C$ (independent of $\epsilon$) such that
$\left|\kappa\right| \le CM\epsilon$ for $\epsilon\in\left[0, \epsilon_0\right]$.
In the remainder of this section, all estimates are valid for sufficiently large $M>0$ independent 
of $\epsilon$ (as stated in the Theorem 3.1), but also allow the possibility that $M>1$ grows as 
$\epsilon\to 0$.

\noindent
{\bf Lemma 3.1} 
{\sc Suppose 
$\left(\xi, r\right)$ remains in $D^\prime$ for $\phi\in\left[\phi_0,\bar\phi\right]$. Then
\be\lbl{5.15}
\xi(\phi)= \xi_0(\phi)r^2(\phi) +
e^{-(\phi-\phi_0)}\left\{ \xi(\phi_0)-\xi_0(\phi_0) r^2(\phi_0)\right\} +O\left(M^3\epsilon\right),
\quad\phi\in\left[\phi_0, \bar\phi\right],
\ee
where $\xi_0(\phi)$ is defined in (\ref{5.16}).
}

\noindent{\bf Proof:}
Using (\ref{5.13}), we integrate (\ref{5.11a}) over $\left[\phi_0,\phi\right]$
\be\lbl{5.17}
e^\phi\xi(\phi)=e^{\phi_0}\xi(\phi_0)+
\bar p_1 \int_{\phi_0}^\phi e^s r^2(s)ds  + \int_{\phi_0}^\phi e^s r^2(s) \tilde P_1(s) ds+ 
\epsilon\int_{\phi_0}^\phi e^s \Xi\left(\xi(s), r(s),s,\epsilon\right) ds.
\ee
Using integration by parts, we have
$$
\int_{\phi_0}^\phi e^s r^2(s)ds=e^\phi r^2(\phi)-e^{\phi_0} r^2(\phi_0)
-2\int_{\phi_0}^\phi e^s r(s)\dot r(s)ds.
$$
From (\ref{5.12}) and (\ref{Dprime}), we find that for $\left(\xi,r\right)\in D^\prime$
$$
\left|\int_{\phi_0}^\phi e^s r(s)\dot r(s)ds\right| \le C_1 e^{\phi-\phi_0}M^3\epsilon,
$$ 
where $C_1>0$ is independent of $\epsilon$.
Therefore,
\be\lbl{5.18}
\int_{\phi_0}^\phi e^s r^2(s)ds=e^\phi r^2(\phi)-e^{\phi_0} r^2(\phi_0)
+e^{\phi-\phi_0}O\left(M^3\epsilon\right).
\ee
Similarly, using integration by parts in the second integral on the right hand
side of (\ref{5.17}), we obtain
\begin{eqnarray}\nonumber
\int_{\phi_0}^\phi e^s r^2(s) \tilde P_1(s) ds &=& e^\phi r^2(\phi) \tilde P_1(\phi)-e^{\phi_0} r^2(\phi_0) \tilde P_1(\phi_0)\\
   & &-2\int_{\phi_0}^\phi e^s r(s) \dot r(s) \tilde P_1(s)ds - \int_{\phi_0}^\phi e^s r^2(s) \tilde P_1^\prime(s)ds.\lbl{A.2}
\end{eqnarray}
Applying integration by parts in the last integral on the right hand side of (\ref{A.2}),
we have
\begin{eqnarray}\nonumber
\int_{\phi_0}^\phi e^s r^2(s) \tilde P_1^\prime(s) ds &=& e^\phi r^2(\phi) \tilde P_1^\prime(\phi)-e^{\phi_0}r^2(\phi_0) \tilde P_1^\prime(\phi_0)\\
&& -2\int_{\phi_0}^\phi e^s r(s) \dot r(s) \tilde P_1^\prime (s)ds - \int_{\phi_0}^\phi e^s r^2(s) \tilde P_1^{\prime\prime}(s)ds.
\lbl{A.3}
\end{eqnarray}
By direct verification, one finds that 
$
\tilde P_1^{\prime\prime}=\left(1+4\beta^2\right)\tilde P_1.
$
The integrals involving $\dot r$ on the right hand sides of (\ref{A.2}) and
(\ref{A.3}) are bounded by $e^{\phi-\phi_0}O(M^3\epsilon)$ since, by (\ref{5.12}), $\dot r=O(\epsilon M^2)$.
Using these observations, from (\ref{A.2}) and (\ref{A.3}) we obtain
\be\lbl{A.4}
\int_{\phi_0}^\phi e^s r^2(s) \tilde P_1(s)ds = e^\phi r^2(\phi)p_1(\phi)-e^{\phi_0} r^2(\phi_0)p_1(\phi_0)+ e^{\phi-\phi_0}O(M^3\epsilon),
\ee
where $p_1(\phi)$ is given by (\ref{5.16}).
The combination of (\ref{5.17}), (\ref{5.18}), and  (\ref{A.4}) yields (\ref{5.15}).
\\
$\Box$

The following lemma shows that the trajectories of (\ref{5.11a}) and (\ref{5.12}),
which stay in $D^\prime$ for sufficiently long time 
enter a small neighborhood of $S$ 
no later than in time $O(\left|\ln\epsilon\right|)$
and remain in this neighborhood as long as they stay in $D^\prime$ (see (\ref{Dprime})).

\noindent
{\bf Lemma 3.2} 
{\sc 
Let $\left(\xi(\phi),r(\phi)\right)$ denote a trajectory of (\ref{5.11a}) and (\ref{5.12}),
which stays in $D^\prime$ for $\phi\in \left[\phi_0, \bar\phi\right]$, $\bar \phi>\phi_1$,
$\phi_1=\phi_0+\left|\ln\epsilon\right|$. Then
\be\lbl{5.20}
\xi(\phi)=\xi_0(\phi) r^2(\phi) +O\left(M^3\epsilon\right),
\ee
for $\phi\ge\phi_1$ and as long as $\left(\xi(\phi),r(\phi)\right)\in D^\prime$.
}

\noindent{\bf Proof:} Denote $\phi_1=\phi_0+\left|\ln\epsilon\right|$.
By plugging in $\phi=\phi_1$ into (\ref{5.15}), we have
\be\lbl{5.21}
\left|\xi(\phi_1) - \xi_0(\phi_1)r^2(\phi_1)\right|\le C_2\epsilon,
\ee
for some $C_2>0$ independent of $\epsilon>0$.
Using Lemma 3.1 again with $\phi_0:=\phi_1$, we obtain
$$
\xi(\phi)= \xi_0(\phi)r^2(\phi) +
e^{-(\phi-\phi_1)}\left\{\xi(\phi_1)-\xi_0(\phi_1)r^2(\phi_1)\right\}
+O(M^3\epsilon),\;\phi\ge \phi_1.
$$
The expression in the curly brackets is $O(\epsilon)$ by (\ref{5.21}).
\\
$\Box$

\noindent {\bf Remark 3.2}
Lemmas 3.1 and 3.2 show that the trajectories of (\ref{5.11a}) and (\ref{5.12})
converge to an exponentially stable manifold $S$, whose  
leading order approximation is given in the definition of $S_0$. The method, 
which we used in Lemmas 3.1 and 3.2 to obtain the
leading order approximation of the slow manifold, can be extended to calculate the
higher order terms in the expansion for $S$.
 
Having shown that the trajectories approach an $O(\epsilon)$ neighborhood of $S_0$ in time 
$O\left(\left|\ln\epsilon\right|\right)$, next we
reduce the dynamics of (\ref{5.11a}) and (\ref{5.12}) to the slow manifold. 
For this, we define
\be\lbl{5.22}
I(\phi)=\left({1\over r(\phi)}+\epsilon q_1(\phi)\right)^2,
\ee
where $q_1^\prime(\phi)=Q_1(\phi)$ and $Q_1(\phi)$ is a trigonometric polynomial on 
the right hand side of (\ref{5.12}). For the sake of definiteness, 
we choose $q_1^\prime(\phi)$ such that 
$$
\int_0^{2\pi\over\beta}q_1(s)ds=0.
$$
The following lemma provides the desired reduction.

\noindent
{\bf Lemma 3.3} 
{\sc For $\phi\ge\phi_1$ and as long as $(\xi,r)\in D^\prime$, $(I,\phi)$ satisfies the 
following system of equations
\begin{eqnarray}\lbl{5.24}
\dot I &=& -2\epsilon^2\left(\mu I+\gamma+Q(\phi)\right)+O(M^3\epsilon^3),\\ \lbl{5.25}
\dot \phi &=& 1+O(M\epsilon),
\end{eqnarray}
where $Q(\phi)$ is a trigonometric polynomial with period $\omega={2\pi\beta^{-1}}$ and zero mean 
$\int_0^\omega Q(\phi) d\phi =0$. The expression for $\gamma$ is given in
(\ref{5.5}).
}

\noindent{\bf Proof:}
The change of variables
\be\lbl{5.26}
r(\phi)={1\over J(\phi)-\epsilon q_1\left(\phi\right)},\quad \mbox{where}\quad q_1^\prime(\phi)=Q_1(\phi),
\ee
in Equation (\ref{5.12}) yields
\be\lbl{5.27}
\dot J ={-\epsilon^2\over J}\left(\left(\mu+\xi Q_2(\phi) \right) J^2 +Q_3(\phi)\right)+O(M^2\epsilon^3).
\ee
After another change of variables, $I=J^2$, we obtain
\be\lbl{5.28}
{d I\over d\phi} = -2\epsilon^2 \left( \left(\mu+ \xi Q_2(\phi)\right) I +Q_3(\phi)\right)
+ O(M^3\epsilon^3),
\ee
By Lemma 3.2, for $\phi\ge\phi_1$
\be\lbl{5.29}
\xi(\phi)=\xi_0(\phi)r^2(\phi) +O(M^3\epsilon).
\ee
By plugging in (\ref{5.29}) into (\ref{5.28}), we obtain
\be\lbl{5.30}
{d I\over d\phi} = -2\epsilon^2\left(
\mu I +\xi_0(\phi)Q_2(\phi)+Q_3(\phi)\right)+
O(M^3\epsilon^3).
\ee
Equation (\ref{5.24}) follows by rewriting  (\ref{5.30}) 
\be\lbl{5.31}
{d I\over d\phi} = -2\epsilon^2\left(
\mu I +\gamma+ Q(\phi)\right)+
O(M^3\epsilon^3),
\ee
where
\be\lbl{5.32}
\gamma = {1\over\omega}\int_0^\omega\left(\xi_0(\phi)Q_2(\phi)+Q_3(\phi)\right)d\phi\;\;\mbox{and}\;\; 
\int_0^{\omega} Q(\phi)d\phi =0,
\ee
and $\xi_0(\phi)$ is given in (\ref{5.16}).
Equation (\ref{5.25})  follows from the last equation in (\ref{5.11}) and the definition of $\phi$.\\
$\Box$

The statements in Theorem 3.1 can now be deduced from Equations (\ref{5.24}) and (\ref{5.25}).
We only need to show that a trajectory with initial condition in $D_0^\prime$ remains in $D^\prime$
for times longer than $O\left(\left|\ln\epsilon\right|\right)$. This follows from the fact that
$\dot I=O\left(\epsilon^2\right)$ in $D^\prime$. Consequently, it takes time $O(\epsilon^{-2})$ for
$I$ to undergo $O(1)$ change necessary for leaving $D^\prime$ from $D^\prime_0$.
We omit any further details. In conclusion, we  note that the domain of validity of the asymptotic 
analysis of this section extends
much further beyond $D^\prime$. Indeed, from (\ref{5.20}) we observe that the remainder term 
tends to $0$ with $\epsilon\to 0$ provided $M=o\left(\epsilon^{-1\over 3}\right)$. This means that
the expansions for $r$ and $\xi$ can be controlled in regions of size up to $O\left(\epsilon^{-1\over 3}\right)$
and $O\left(\epsilon^{-2\over 3}\right)$ respectively. For the original variables $\rho,$
$x_1$, these estimates translate into  $O\left(\epsilon^{2\over 3}\right)$
and $O\left(\epsilon^{4\over 3}\right)$ respectively.

We end this section by deriving a useful estimate for the time of flight of trajectories
passing close to the saddle-focus. For this, consider an initial value problem for system of equations
(\ref{5.11a}) and (\ref{5.12}). Suppose that the initial condition implies that $I\left(\phi_0\right)=I_0>0$.
We would like to know how long it takes for $I$ to reach a given value $\bar I=I\left(\phi_0+\Delta\right)$,
$\Delta>0$.
We assume that $I_0$ and $\bar I$ are sufficiently separated, e.g., $2\bar I\le I_0$.
To estimate $\Delta$, note that
\be\lbl{5.51}
\Delta=\Delta_1+\Delta_2,
\ee
where $\Delta_1=O\left(\left|\ln\epsilon\right|\right)$ is time necessary to reach 
an $O(\epsilon)$ neighborhood of $S_0$ from $I_0$. Denote 
\be\lbl{5.52}
  I_1:= I\left(\phi_1\right)=I_0+O\left(\epsilon^2\left|\ln\epsilon\right|\right),\;
  \mbox{where}\; \phi_1=\phi_0+\Delta_1.
\ee
The second term on the right hand side of (\ref{5.51}) can be estimated by integrating (\ref{5.24})
over $\left[ \phi_1, \phi_1+\Delta_2\right]:$
\be\lbl{5.53}
\bar I= I_1 e^{-2\alpha\Delta_2}+{\gamma\over\mu}\left(1- e^{-2\alpha\Delta_2}\right)
-2\epsilon^2\int_0^{\Delta_2} e^{-2\alpha\left(\Delta_2-s\right)} Q\left(\phi_0+s\right)ds +
O\left(\epsilon^3\right).
\ee
The integral on the right hand side of (\ref{5.53}) is bounded uniformly for $\Delta_2>0$.
This observation combined with (\ref{5.52}) implies 
\be\lbl{5.54}
\bar I= I_0 e^{-2\alpha\Delta_2}-{\gamma\over\mu}\left(1- e^{-2\alpha\Delta_2}\right)
+O\left(\epsilon^2\right).
\ee
Note that the contribution of $O\left(\epsilon^2 \left| \ln\epsilon \right|\right)$ term in (\ref{5.52})
to (\ref{5.54}) is negligible provided $\Delta_2$ is sufficiently large.
From  (\ref{5.54}), we obtain the desired estimate
\be\lbl{5.55}
\Delta={1\over 2\mu\epsilon^2}\ln \left({I_0 +{\gamma\over \mu}\over \bar I +{\gamma\over\mu}+O(\epsilon^2)}\right).
\ee

\section{The oscillations}
\setcounter{equation}{0}

In the present section, we use the local analysis of Section 3 to study certain oscillatory patterns
arising in the model of solid fuel combustion (\ref{6.1})-(\ref{6.3}). By plugging in the values of the 
parameters of (\ref{6.1})-(\ref{6.3}) into the expression for $\gamma$ (\ref{5.5}), we find that $\gamma(p)$
is a quadratic function with two zeros at $p_1\approx 0.34$ and $p_2\approx 2.73$. 
These values are in a good agreement with those obtained by numerical bifurcation analysis in \cite{FKRT}.
The quadratic character of $\gamma(p)$ is explained by the fact that $\gamma$ is determined by the
second order terms in the Taylor expansions of $h_{2,3}(0,0)$.
We concentrate on the parameter region around $p=p_2$, where $\gamma(p)$ changes its sign.
For small positive $\alpha$, there are two dynamical regimes: for values of $p$ lying 
to the left and to the right of some $O(\epsilon)$ neighborhood
of $p_2$ corresponding to the supercritical and subcritical AHBs. The transition region between these two
parameter regimes adds to the repertoire of qualitatively distinct dynamical behaviors. Below, we study
these three cases in more detail.

\subsection{The supercritcal AHB: geometry of the periodic orbit}
\begin{figure}
\begin{center}
{\bf a}\epsfig{figure=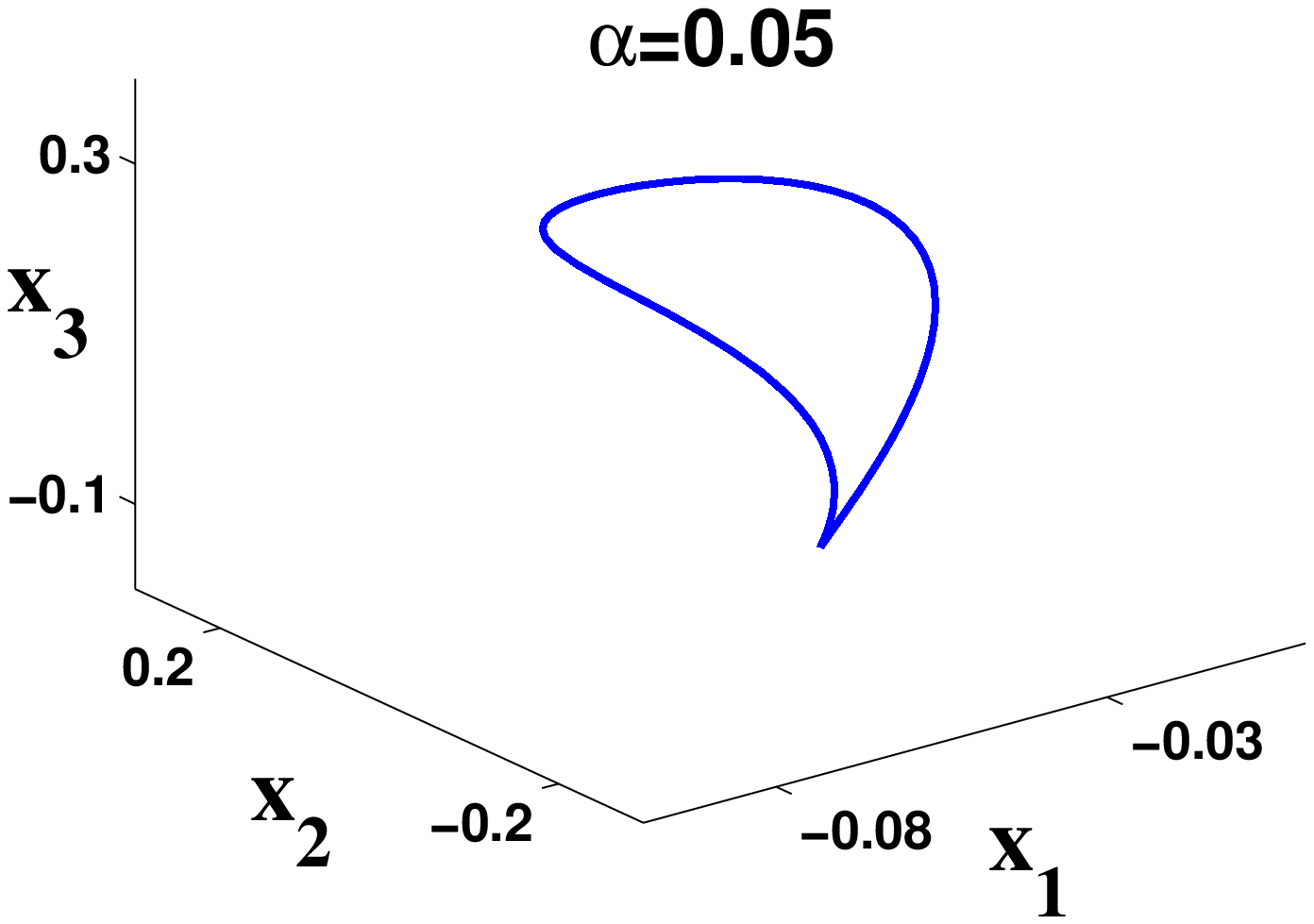, height=1.8in, width=2.0in, angle=0}
{\bf b}\epsfig{figure=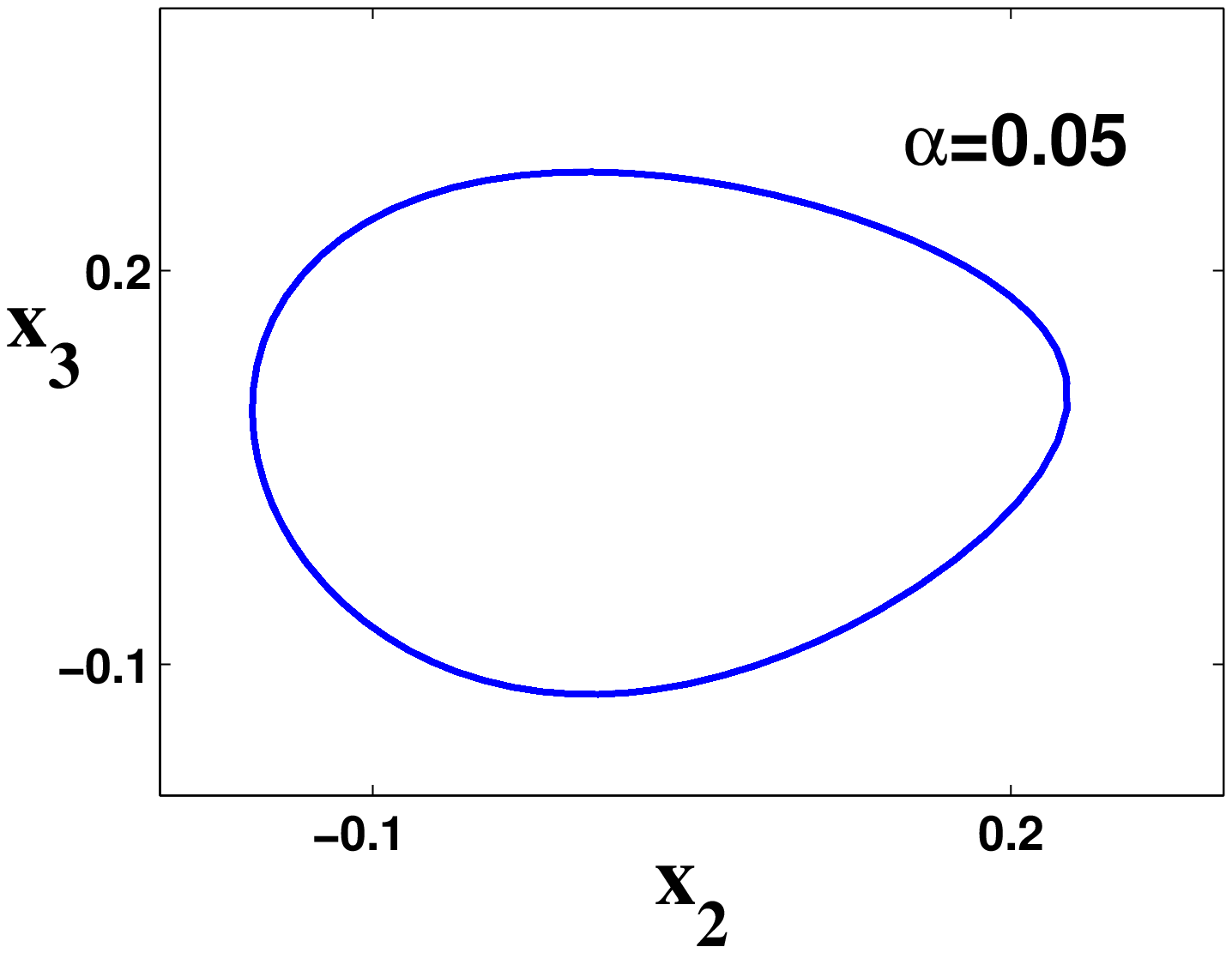, height=1.8in, width=2.0in, angle=0}
{\bf c}\epsfig{figure=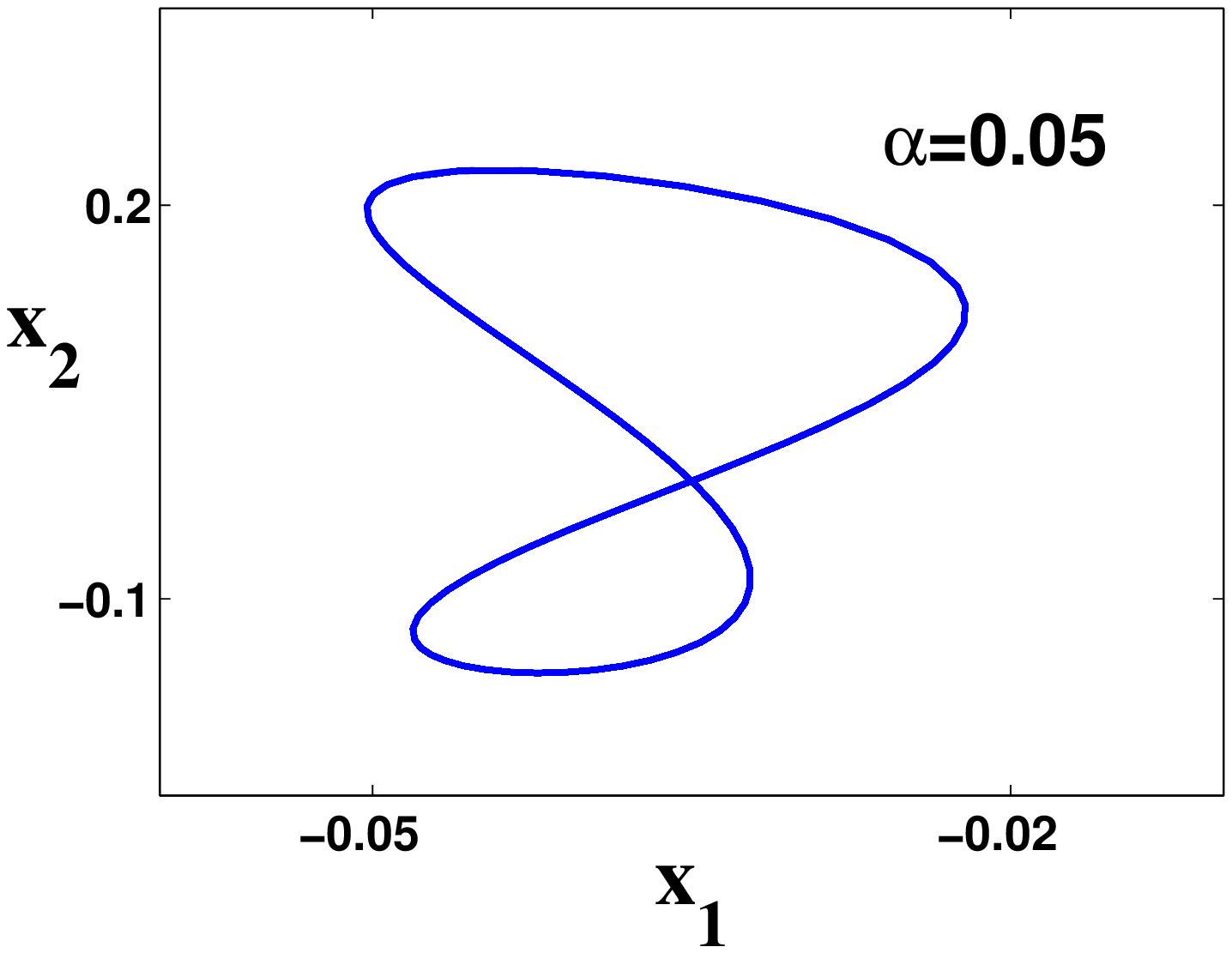, height=1.8in, width=2.0in, angle=0}
\end{center}
\caption{ A periodic orbit born from the supercritical AHB ({\bf a}), and
its projections onto $x_2-x_3$({\bf b}) and $x_1-x_2$ planes ({\bf c}).
Note that the periodic orbit shown in ({\bf a}) is not planar. This 
accounts for the presence of the self-intersection in the its projection
in ({\bf c}) (see text for details).
}\label{f.6}
\end{figure}

\setcounter{equation}{0}
It follows from Theorem 3.1 that for $\gamma(p)<0$ and sufficiently small $\alpha>0$,
(\ref{6.5}) has a stable limit cycle $O_\alpha$, whose leading order approximation is given by
\be\lbl{4.1}
x(\theta)=\left(\bar\rho^2\left(a+A \cos 2\theta\right), \bar\rho\cos\theta, \bar\rho\sin\theta\right),
\quad\bar\rho=\sqrt{\alpha\over -\gamma},\quad \theta\in[0, 2\pi)
\ee
(see Figure \ref{f.6}a).
Moreover, the analysis of Section 3 shows that all trajectories starting from a sufficiently small 
neighborhood of the origin and not belonging to $W^s(O)$, converge to $O_\alpha$. The leading
order approximation of $O_\alpha$ in (\ref{4.1}) reveals a remarkable property of the oscillations 
generated by the limit cycle born from the supercritical AHB: the frequency of oscillations in $x_1$
is twice as large as that of the oscillations in $x_{2,3}$. To explain the frequency doubling effect,
we recall that 
\be\lbl{4.2}
\dot\theta=\beta+O(\epsilon).
\ee
Therefore, (\ref{4.1}) implies that, unless $A=0$,
the frequency of oscillations in $x_1$ is $\beta\pi^{-1}$, while that of oscillations in $x_{2,3}$
is ${\beta\over 2\pi}$. The latter coincides with the frequency of the bifurcating periodic solution.
Below, we complement the analytical explanation of the frequency doubling with the geometric interpretation.

The oscillations in $x_{2,3}$ can be well understood using the topological normal form of the AHB
\cite{GEO,GH, KUZ}. Indeed, since at the bifurcation the center manifold at the origin is tangent to the $x_2-x_3$
plane, a standard treatment of the AHB using the center manifold reduction \cite{Carr, GH} shows that the projection
of the bifurcating limit cycle onto $x_2-x_3$ plane to leading order is a circle (Figure \ref{f.6}b)
and the projection
of the vector field is given by the equation of the angular variable (\ref{4.2}). Therefore, the oscillations in
$x_{2,3}$ are approximately harmonic with the period equal approximately to $2\pi\beta^{-1}$. The topological
normal from, however, does not describe the oscillations in $x_1$. For this, one needs to take into account the
geometry of the bifurcating periodic orbit as a curve in $\Re^3$. The geometry of $O_\alpha$ is fully determined by 
the two geometric invariants: the curvature, $\mbox{k}(\theta)$, and the torsion, $\kappa(\theta)$. For the purposes
of the present discussion, we need only the latter, but we compute both invariants for completeness.
After some algebra, the parametric equation for the periodic orbit (\ref{4.1}) yields  
\begin{eqnarray}\nonumber 
\mbox{k}(\theta)&=& {\left|\left[\dot x, \ddot x\right]\right|\over \left|\dot x\right|^3}=
\sqrt{1+2A^2\left(3\cos 4\theta +5\right)}+\mbox{h.o.t},\\
\lbl{4.3} 
\kappa(\theta) &=& {\left(\dot x, \ddot x, \dddot x\right)\over\left|\left[\dot x, \ddot x\right]\right|^2}=
{-\gamma\over\alpha} {6A\sin 2\theta\over 1+2A^2\left(3\cos 4\theta +5\right)}
+\mbox{h.o.t},\quad \theta\in[0, 2\pi).
\end{eqnarray}
The geometry of the leading order approximation of the periodic orbit is determined by the three
parameters $\alpha,$ $\gamma,$ and $A$. The former two parameters are the same as the parameters in the  
topological normal form of the nondegenerate AHB \cite{KUZ}; while the latter captures the geometry of the slow manifold (or unstable
manifold) near the origin. From the geometrical viewpoint the bifurcation is degenerate if either $\gamma$ or $A$ is equal
to zero. The latter condition holds if and only if the slow manifold is either a plane or a circular paraboloid
near the origin. In this case, Equation (\ref{4.3})
implies that (to leading order) the bifurcating orbit lies in a plane. Generically, $O_\alpha$ is not planar 
(see Figure \ref{f.6}a). 
The fact that the periodic orbit is generically not planar combined with the
symmetry of the orbit about the $x_1-$axis implies that its projection onto any plane containing $x_1-$axis has to have
a self-intersection (see Figure \ref{f.6}c). From  Figure \ref{f.6}c, it is clear that as the phase point goes 
around $O_\alpha$ once, $x_1$
has to trace its range at least twice. Therefore, the frequency of oscillations in $x_1$ has to be at least twice
as high as that of the oscillations in $x_{2,3}$. The above discussion implies that the frequency doubling of oscillations in
$x_1$ is a generic  geometric property of the AHB.  

\subsection{The subcritcal AHB: multimodal oscillations}

If the AHB is subcritical ($\gamma>0$) the loss of stability of the equilibrium at the origin of (\ref{6.5}) 
results in the creation of multimodal trajectories which spend a considerable amount of time near a weakly 
unstable equilibrium. To describe the resultant dynamics we give the following definition.

\noindent
{\bf Definition 4.1} 
{\sc We say that a trajectory of (\ref{6.5}) undergoes multimodal oscillations for $t\ge 0$ if
there exist positive constants $r_1$ and $r_2$ independent of $\epsilon$ and an unbounded sequence of times
$$
0\le t_1<t_2<\dots,\quad \lim_{i\to\infty} t_i=\infty,
$$ 
such that 
\begin{description}
\item[a.] $\left|x\left(t_{2i}\right)\right|\le r_1\epsilon^2\;\;\&\;\;\left|x\left(t_{2i-1}\right)\right|> r_2,\quad i=1,2,\dots$,
\item[b.] $\left(\exists t^\prime> t_1:\;\left| x(t^\prime)\right| \le r_1\epsilon^2\right)\;\Rightarrow\; 
\exists i\in\Na\;:\;
\left|x(t)\right|\le r_1\epsilon^2, 
t\in\left[ \min(t^\prime, t_{2i}), \max(t^\prime, t_{2i})\right]\;\;\forall t^\prime \in \left(t_{2i-1}, t_{2i+1}\right).$
\end{description}
The time intervals $\tau_i=t_{2i+1}-t_{2i-1},\; i=1,2,\dots$ are called ISIs.
}

The proximity of (\ref{6.5}) to the AHB alone is clearly not sufficient to account for  the appearance of 
the multimodal oscillations in (\ref{6.5}).
Below, we formulate two additional assumptions on the vector field outside of the small neighborhood of the equilibrium, which
are relevant to (\ref{6.1})-(\ref{6.3}).
Under these conditions, we show that the subcritical AHB results in sustained multimodal
oscillations. In addition, we determine the asymptotics of the ISIs for positive $\alpha=O(\epsilon^2)$. 
For (\ref{6.5}) to generate multimodal oscillations, it is necessary that the trajectories leaving an $O(\epsilon^2)$
neighborhood of the origin reenter it after some interval of time. Therefore, the vector field in $O(1)$ neighborhood of the
origin must provide a return mechanism. Our second assumption on the global vector field of (\ref{6.5})  
is that the trajectories approaching the origin
along a $1D$ stable manifold, $W^s(O)$, are subject to a strong contraction toward $W^s(O)$, i.e., in a small neighborhood
of an $O(1)$ segment of $W^s(O)$, the projection of the vector field onto a plane transversal to $W^s(O)$ is sufficiently stronger than
that along $W^s(O)$ (Figure \ref{f.2}a). This guarantees that the trajectories entering such 
region of strong contraction approach $W^s(O)$ very 
closely and follow it to an $O(\epsilon^2)$ neighborhood of the unstable equilibrium, from where they are propelled away along 
the unstable manifold, $W^u(O)$. Due to the proximity of (\ref{6.5}) to the AHB, the motion away from the origin is very slow.
This results in the pronounced intervals of the small amplitude oscillations. Below, we  summarize these observations into 
two formal assumptions on the global vector field {\bf (G1)} and {\bf(G2)}. For this, we first need to introduce some
auxiliary notation.
For analytical convenience, we assume that in an $O(1)$ neighborhood of the origin, $W^s(O)$ can be and has been straightened by  
a smooth change of coordinates. More specifically, the nonlinear terms in (\ref{6.5}) satisfy
\be\lbl{4.3a}
h_{2,3} \left(x_1,0,0,\alpha\right)=0,\quad w(x_1,\alpha)=-x_1+h_1\left(x_1,0,0,\alpha\right)>0,\; 
x_1\in [d_1,0),
\ee
for some $d_1<0$. Note that $\left(d_1, 0,0\right)$ is assumed to be sufficiently far away from the origin (see Figure \ref{f.8}).  
To describe the mechanism of return, we introduce two crossections:
\be\lbl{4.4}
 \Sigma^+=[-c_1\delta^2, c_1\delta^2]\times D_\delta\quad\mbox{and}\quad\Sigma^-=\left\{d_1\right\}\times D_{\delta_1},
\ee
where $D_\delta \subset\Re^2$ denotes a disk of radius $\delta$ centered at the origin:
\be\lbl{4.4a}
D_\delta=\left\{ y\in\Re^2:\;\left|y\right|\le \delta\right\}.
\ee
Positive constants $\delta$ and $\delta_1$ are sufficiently small (see Figure \ref{f.2}a) and
$c_1>0$ is chosen so that $\Sigma^+$ intersects the slow manifold $S$ transversally.
In addition, we require that $\delta=o(\epsilon^{2/3})$ to guarantee that $\Sigma^+$
belongs to the region of validity of the local analysis of Section 3 (see Remark 3.1c).
Let $x_0\in \Sigma^+$ and  consider a trajectory of (\ref{6.5}) 
starting from $x_0$.
We assume that for sufficiently small $\epsilon>0$ and $\alpha=O(\epsilon^2)$,
every such trajectory  intersects $\Sigma^-$ from the left. Denote the point of the first 
intersection by $Q(x_0)\in\Sigma^-$.
We assume that 

\noindent
{\bf(G1)} the first return map $Q:\Sigma^-\to\Sigma^+$ depends smoothly on $\epsilon$ and 
$\min_{x\in Q\left(\Sigma^+\right)} \left|x\right|\ge\zeta>0$.

To measure the rate of contraction toward $W^s(O)$, we consider a $2\times 2$ matrix
\be\lbl{4.5}
A\left(x_1,\alpha\right)=\left(\begin{array}{cc}
{\partial f_2\over \partial x_2} &{\partial f_2\over \partial x_3} \\
 {\partial f_3\over \partial x_2} &{\partial f_3\over \partial x_3} 
\end{array}
\right)_{(x_1,0,0,\alpha)}.
\ee
Let $\lambda_{1}(x_1,\alpha)\le \lambda_2(x_1,\alpha)$ denote the eigenvalues of the symmetric matrix 
$$
A^s\left(x_1,\alpha\right)={1\over 2}\left(A\left(x_1,\alpha\right)+A^T\left(x_1,\alpha\right)\right).
$$
Denote $\underline\lambda(x_1,\alpha)=-\lambda_1(x_1,\alpha)$ and 
$\bar\lambda(x_1,\alpha)=-\lambda_2(x_1,\alpha)$. 
We assume that for sufficiently small $\epsilon>0$ and $\alpha=O(\epsilon^2)$,

\noindent
{\bf(G2)} 
\begin{eqnarray}\lbl{4.6}
{\p \over \p x_1} \bar\lambda(x_1,\alpha) <0, & x_1\in[d_1,0];&\\
\lbl{4.7}
\exists\; d_2 \in (d_1,0):\; \bar\lambda(x_1,\alpha)>0,& x_1\in[d_1,d_2]&\quad \& 
\quad \min_{x_1\in[d_1, d_2]} {\bar\lambda(x_1,\alpha)\over w(x_1,\alpha)} = 
O\left(\left|\ln\epsilon\right|\right),
\end{eqnarray}

\noindent
{\bf(G3)} 
$$
\max_{x_1\in[d_1, 0)} {\underline\lambda(x_1,\alpha)\over w(x_1,\alpha)} = 
           O\left(\left|\ln\epsilon\right|\right).
$$

Under these conditions, we have

\noindent
{\bf Theorem 4.1}
{\sc
Let $\gamma>0$, conditions in (\ref{4.3a}), {\bf (G1)} and {\bf (G2)} hold. 
Then for sufficiently small $\epsilon>0$ and $\alpha=O\left(\epsilon^2\right),$ a trajectory
of (\ref{6.5}) with initial condition from an $O(\epsilon^2)$ neighborhood of the
origin and not belonging to $W^s(O)$ undergoes multimodal oscillations.
The ISIs are uniformly bounded from below
\be\lbl{4.8b}
\tau_i\ge \tau^-={1\over 2\alpha}\ln\left(1+{\alpha\over\gamma\epsilon^4}C^-\right),\;\; i=2,3, \dots.
\ee
If, in addition, {\bf (G3)} holds then the ISIs satisfy two-sided bounds
\be\lbl{4.8a}
\tau^-\le \tau_i \le \tau^+={1\over 2\alpha}\ln\left(1+{\alpha\over\gamma\epsilon^{4+\chi}}C^+\right),\;\; i=2,3, \dots,
\ee
for some $\chi>0$.
Positive constants $C^\pm$ do not depend on $\epsilon,$ $\alpha,$ and $\gamma$.
}
\newpage

\noindent
{\bf Remark 4.1}
\begin{description}
\item[(a)]
The principal assumptions on the global vector field are formulated in {\bf (G1)} and
(\ref{4.7}). The condition in (\ref{4.6}) makes the derivation of certain  estimates
in the proof of the theorem easier and is used for analytical convenience. Likewise,
{\bf (G3)} is not essential for the proposed mechanism. However, without this condition
obtaining the two-sided estimates for the ISIs requires an additional argument in the 
proof. Condition {\bf (G3)} means that the two eigenvalues of $A$ have the same order
of magnitude. This condition is not restrictive.
\item[(b)]
For fixed $\gamma>0$, the estimates in (\ref{4.8a}) and (\ref{4.8b}) can be rewrittten as 
\be\lbl{4.8c}
{\ln(1+\tilde C^-\alpha)\over 2\alpha} \le \tau_i\le {\ln(1+\tilde C^+\alpha)\over 2\alpha},
\ee
with constants $\tilde C^\pm$ independent from $\alpha$. 
Similarly, for fixed $\alpha$ and varying $\gamma>0$ 
inequalities in (\ref{4.8a}) and (\ref{4.8b}) imply
\be\lbl{4.8d}
\bar C^- -{1\over 2\alpha}\ln\gamma \le \tau_i \le \bar C^+ -{1\over 2\alpha}\ln\gamma,
\ee
where constants $\bar C^\pm=O\left(\left|\ln\epsilon\right|\right)$ do not depend on $\gamma>0$.
\end{description}

\begin{figure}
\begin{center}
\epsfig{figure=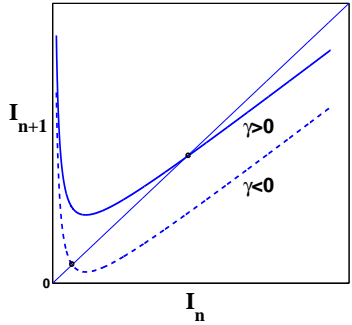, height=1.8in, width=2.0in, angle=0}
\end{center}
\caption{ 
The interpretation of the bifurcation scenario arising in (\ref{6.1})-(\ref{6.3}) for
fixed positive value of $\alpha=O(\epsilon^2)$ and varying $\gamma$ near $0$ using
a one-parameter family of $1D$ maps. For negative values of $\gamma=O(1)$, the map
in (\ref{9.4}) has a stable fixed point (upper graph), which corresponds to a 
stable limit cycle born from the supercritical AHB. For increasing values of $\gamma$
the fixed point looses stability via a period-doubling bifurcation (lower graph),
which initiates a period-doubling cascade leading to the formation of the chaotic attractor.
}\label{f.7}
\end{figure}

\subsection{The transition from subcritical to supercritical AHB}

The numerical experiments presented in Section 2 show that the transition from subcritical to
supercritical AHB contains a distinct bifurcation scenario involving the formation
of chaotic attractor via a period-doubling cascade. The analytical explanation of these 
phenomena is outside the scope of the present paper.
Below we comment on the difficulties arising in the analytical treatment of
this problem. In the present section, we use a combination of the analytic and numerical techniques
to elucidate the origins of the complex dynamics arising in the parameter regime near
the border between regions of sub- and supercritical AHB.
The principal features of this bifurcation scenario are 
summarized in Figures \ref{f.4} and \ref{f.5}. In these numerical experiments, we
kept $\alpha$ fixed at a small positive value and varied $\gamma$.
By taking progressively smaller values of $\gamma>0$, one first observes that
the multimodal patterns exhibit an increase in the ISIs (Figure \ref{f.4}a,b; see also Figure \ref{f.3}b).
We consider these oscillatory patterns regular (even if they are not periodic),
because the timings of the spikes remain within narrow bounds in accord with (\ref{4.8a}).
For smaller values of $\gamma>0$, the oscillatory patterns become irregular and
are characterized by long intervals of oscillations of small and intermediate amplitudes
between successive spikes (Figure \ref{f.4} c,d).
As $\gamma$ becomes negative, the trajectories lose spikes and
consist of irregular oscillations 
(Figure \ref{f.5} a,b). Further decrease of $\gamma$ leads the system through
the reverse cascade of period-doubling bifurcations (Figure \ref{f.5} c-h).
The period-doubling cascade terminates with the creation of the limit cycle,
which can be followed to a nondegenerate supercritical AHB by letting $\alpha\to 0$
(Figure \ref{f.5} g,h).
To account for the bifurcation scenario described above, we construct the first
return map. Since near the origin the trajectories spend most of the time in the vicinity of
the $2D$ slow manifold, the first return map is effectively one-dimensional. The distinct
unimodal structure of the $1D$ first return map affords a lucid geometric interpretation
for the bifurcation scenario in the $3D$ systems of differential equations near the
transition from sub- to supercritical AHB. Specifically, we show that the mechanism
for generating complex dynamics in the continuous system in this parameter regime is
the same as in the classical scenario of the period-doubling transition to chaos in
the one-parameter families of unimodal maps \cite{GH}. 
\begin{figure}
\begin{center}
{\bf a}\epsfig{figure=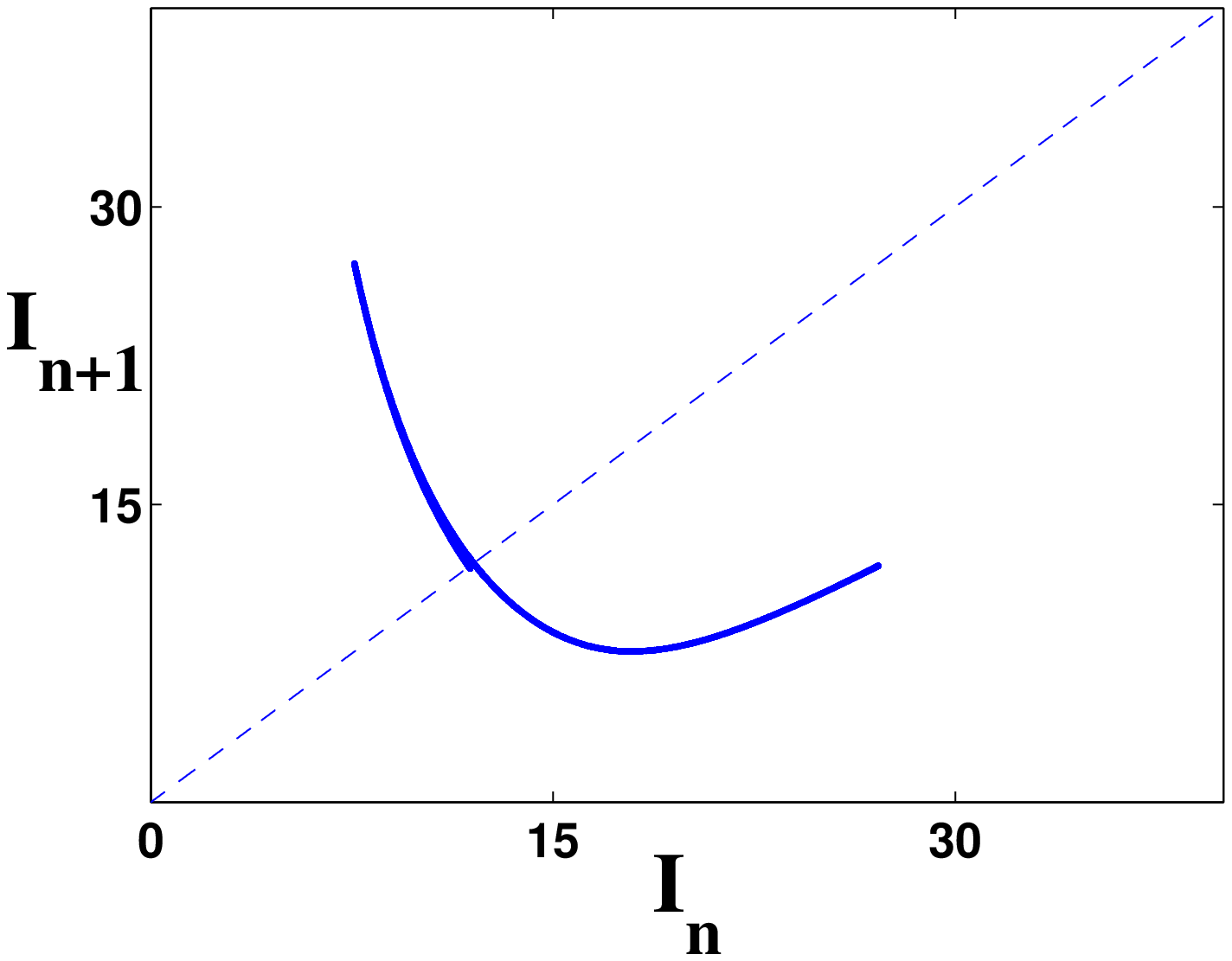, height=1.8in, width=1.8in, angle=0}\quad
{\bf b}\epsfig{figure=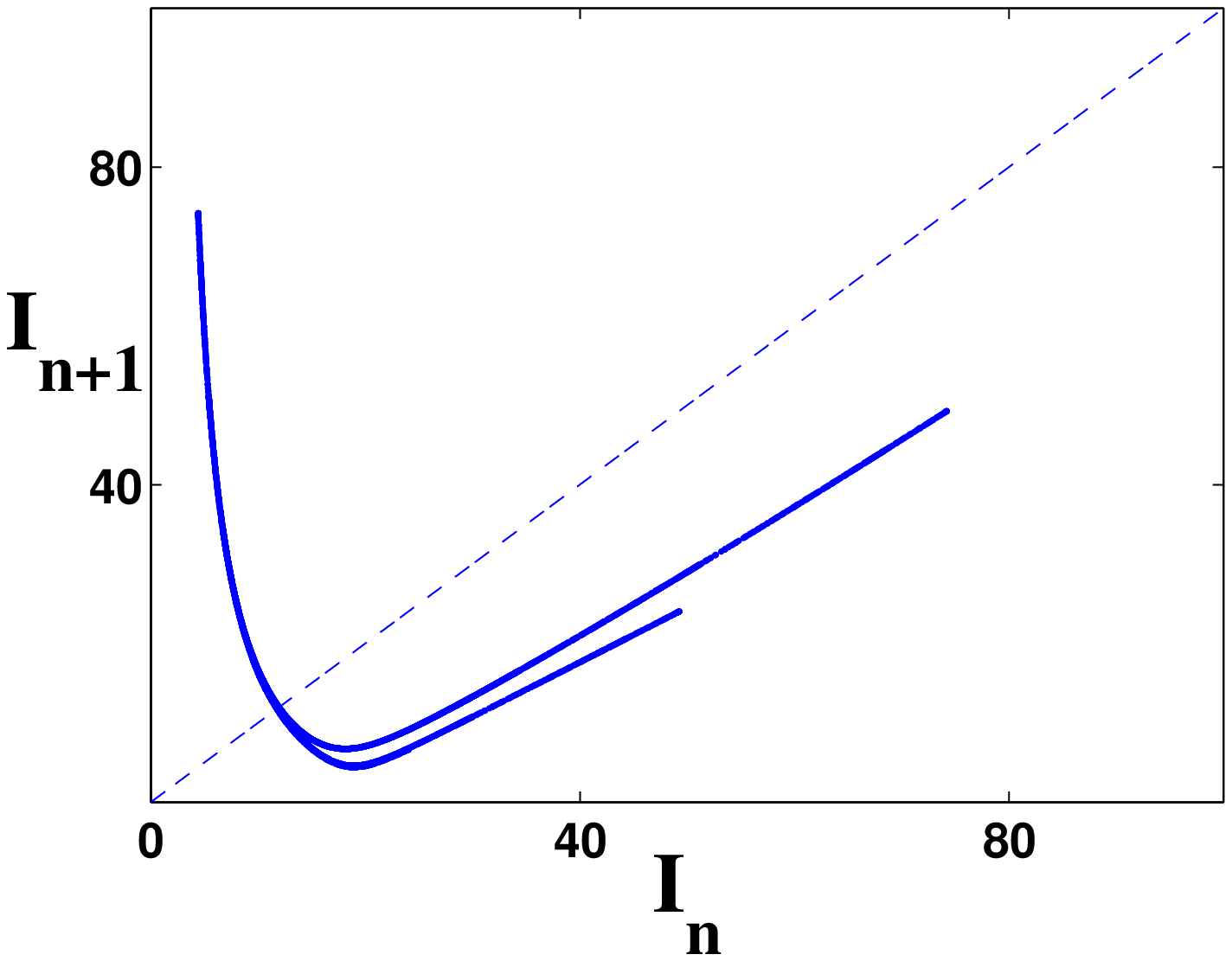, height=1.8in, width=1.8in, angle=0}\\
{\bf c}\epsfig{figure=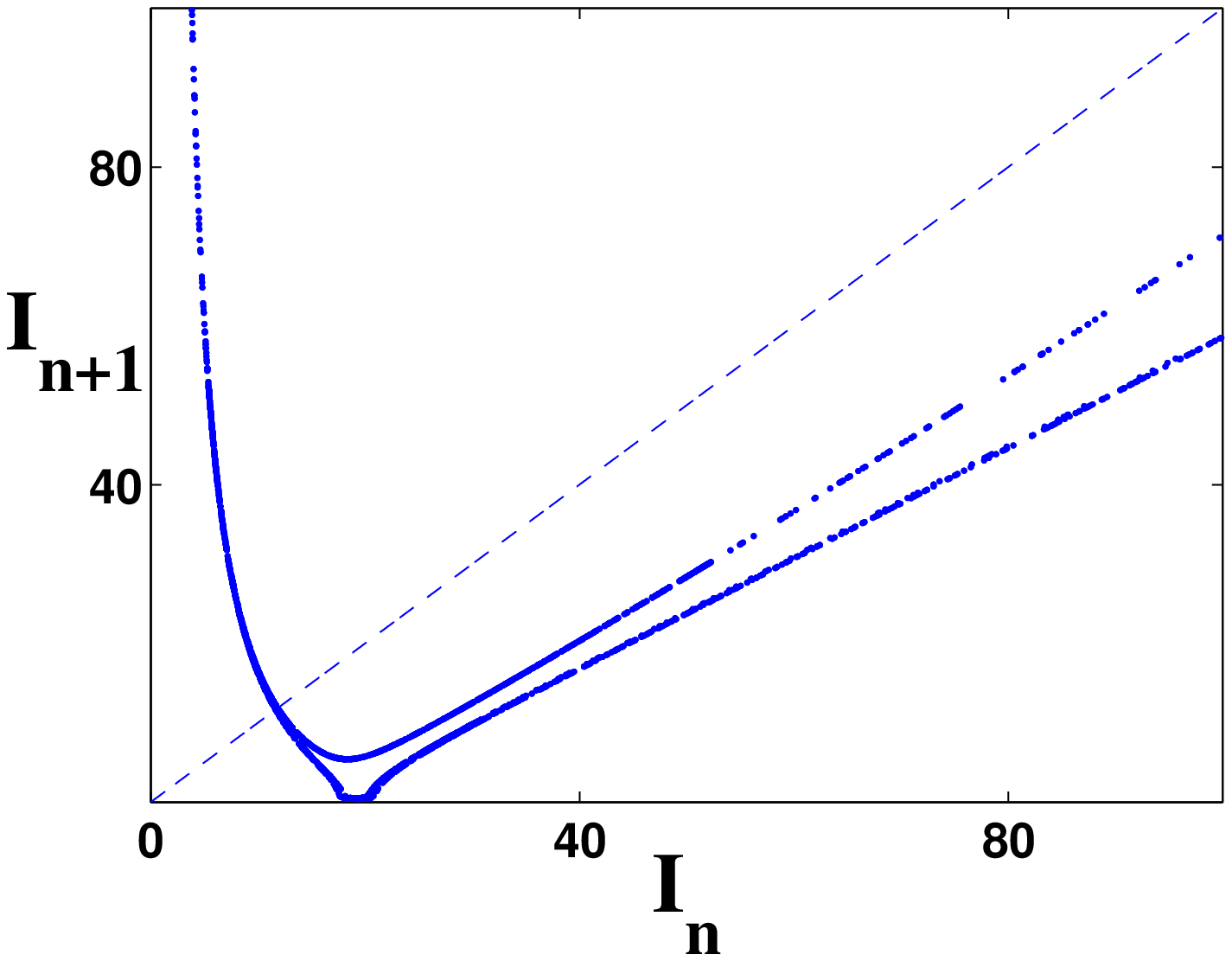, height=1.8in, width=1.8in, angle=0}\quad
{\bf d}\epsfig{figure=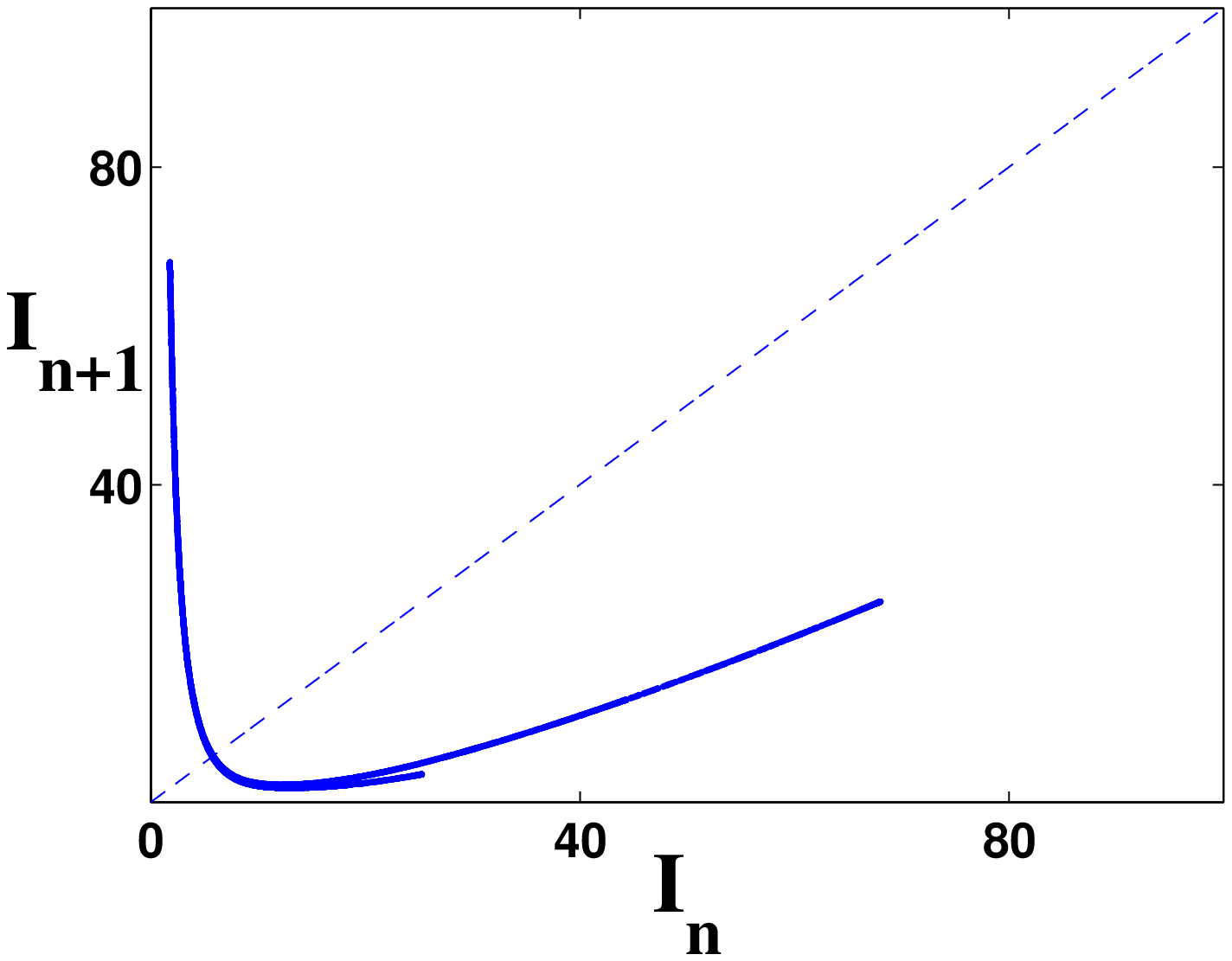, height=1.8in, width=1.8in, angle=0}
\end{center}
\caption{ The plots of first return map for $\mI$ (\ref{9.6}) computed for
({\bf a}) $\alpha= 0.0445$, $\gamma=-0.0486$,
({\bf b}) $\alpha= 0.0445$, $\gamma=-0.0126$,
({\bf c}) $\alpha= 0.0445$, $\gamma=0.0290$,
({\bf d}) $\alpha=0.12 $, $\gamma=-1.0741$.
The plots in ({\bf a}-{\bf c}) illustrate the bifurcations in the family of the first return maps
for the values of $\gamma$ near zero. The map shown in ({\bf d}) corresponds to the 
irregular oscillations following the period-doubling cascade for increasing values
of $\alpha$ in the case of supercritical AHB (see Remark 4.2). 
}\label{f.8}
\end{figure}

We next turn to the derivation of the first return map.
We start with extending the asymptotic analysis of Section 3 to cover the case of $\gamma=o(1)$. 
This requires computing additional terms on the right hand side of the reduced equation
(\ref{5.24}), because already for 
$\left|\gamma\right|=O(\epsilon)$, the term involving $\gamma$ on the right hand side of (\ref{5.24})
is comparable with the $O(\epsilon^3)$ remainder term.
We then use the more accurate approximation for the equation for $I$ to compute a
$1D$ mapping:
\be\lbl{9.0}
G:\; I(\phi)\mapsto I\left(\phi + 2\pi\beta^{-1}\right),
\ee
which describes how  $I$ changes after one cycle of oscillations. 
As noted above, the reduced equation (\ref{5.24}) is not suitable for analyzing
the case of small $\left|\gamma\right|$ and one needs to include more terms in the expansion on the right hand side of (\ref{5.24}).  
In analogy with the topological normal form for the degenerate AHB 
\cite{KUZ, GL}, we expect that terms up to $O(\epsilon^5)$) are needed in (\ref{5.24}) to resolve the dynamics 
for small $\left|\gamma\right|$.  This can be achieved by a straightforward albeit tedious calculation.  Below we describe the 
formal procedure of obtaining the required expansions.  The justification of these expansions can be done in complete analogy 
with the analysis of Section 3.  First, we compute terms on the right hand sides of (\ref{5.11a}) and (\ref{5.12}) up to $O(\epsilon^3)$ and $O(\epsilon^6)$ respectively.  Then we look for a solution of (\ref{5.11a}) in the following form:
\be\lbl{9.1}
\xi(\phi)=\xi_0(r,\phi)+\epsilon\xi_1(r,\phi)+\epsilon^2\xi_2(r,\phi)+O(\epsilon^3).
\ee
By taking an initial  condition from $D_0$ and using Anzats (\ref{9.1}) in (\ref{5.11a}), we recover $\xi_0(r,\phi)$ 
(see (\ref{5.16})) and find the next 
two terms in the expansion of $\xi(r,\phi)$: $\xi_1(r,\phi)$ and $\xi_2(r,\phi)$.  Thus, we obtain the approximation of the slow manifold 
with the accuracy $O(\epsilon^3)$.  Next, we plug in (\ref{9.1}) into (\ref{5.12}) and collect terms multiplying equal powers of $\epsilon$
to obtain 
\be\lbl{9.2}
\dot r = \epsilon R_1(r,\phi)+\dots+\epsilon^5 R_5(r,\phi)+O(\epsilon^6),
\ee
where $R_k(r,\phi)$ are $2\pi\beta^{-1}$– periodic functions of the second argument.
Finally, by rewriting (\ref{9.2}) in terms of $I$ (see (3.25)), integrating it over one period of oscillations, 
$\omega$, and disregarding $O(\epsilon^5)$ terms, we obtain a map for the change in $I$ after one cycle of oscillations
\be\lbl{9.4}
I_{n+1}=G(I_n),\;\; G(I)\equiv I-2\epsilon^2\omega\left(\mu I +\gamma +{\epsilon^2c\over I}\right),\; \omega=2\pi\beta^{-1}.
\ee
Except for the last term in the definition of $G$, the map in (\ref{9.4}) follows from the reduced system (\ref{5.24}) and (\ref{5.25}).  
The last term is only needed if the value of $\left|\gamma\right|$ does not exceed $O(\epsilon)$.  Our calculations show that 
the second Lyapunov coefficient, $c$, is negative for the values of parameters used in (\ref{6.1})-(\ref{6.3}).
Therefore, for sufficiently small $\epsilon>0$, (\ref{9.4}) defines a unimodal map.  Away from an $O(\epsilon^2)$ neighborhood of 
$0$, the graph of  $\tilde I=G(I)$ is almost linear with a weakly attracting slope $1-2\omega\mu\epsilon^2$ for $\mu>0$ (Figure \ref{f.7}).
For $\mu>0$, $G$ has a unique fixed point:
\be\lbl{9.5}
\bar I=\left\{
\begin{array}{cc}
{-\gamma\over\mu} +{c\over\gamma}\epsilon^2+O(\epsilon^3),& \gamma<0,\\
\sqrt{-c\over\mu}\epsilon +O(\epsilon^2),&\gamma=0,\\
{-c\over\gamma}\epsilon^2+O(\epsilon^3),& \gamma>0.
\end{array}
\right.
\ee
For negative values of $\gamma=O(1)$, $I=\bar I$ is a stable fixed point of $G$, which corresponds to a stable limit cycle  born
from a  supercritical AHB (see Figure \ref{f.5} g,h). It follows from (\ref{9.4}) that for increasing values of $\gamma$, 
the graph of the map moves down. Consequently, 
the fixed point moves to the left and eventually looses stability via a period-doubling bifurcation (Figure \ref{f.5} e,f).  
Further increase in $\gamma$, yields more period-doubling bifurcations (Figure \ref{f.5} c,d).  Given the unimodal character of 
$G$, one expects that this sequence of period-doubling bifurcations eventually leads to the formation of a chaotic attractor 
(Figure \ref{f.5} a,b). The unimodal character of the map combined with the manner of its
dependence on $\gamma$ suggest a clear geometric mechanism for the formation of the chaotic attractor and the subsequent
period-doubling cascade arising near the transition from sub- to supercritical AHB (Figure \ref{f.5}a-h). 
We now outline the  limitations of the asymptotic analysis of this section and difficulties
arising in the justification of (\ref{9.4}).
According to (\ref{9.5}) already at the moment of the first period-doubling bifurcation, fixed point $\bar I$
belongs to an $O(\epsilon^2)$ neighborhood of $0$. In this neighborhood,
$\rho=O(1)$ (see (\ref{5.22}))  and, therefore, $\bar I$ lies outside of the region of validity of
the asymptotic analysis. In this region, (\ref{9.4}) may only be considered as a formal 
asymptotic expression. A rigorous justification of the from of the map in the 
boundary layer meets substantial analytical difficulties: it requires introducing an additional set 
of intermediate asymptotic expansions and matching them with those obtained in Section 3. 
We do not address this problem in the present work. 
Below, we resort to using numerical techniques to verify the principal features
of the first return map suggested by the asymptotic analysis: the unimodality of $G$ and its
dependence on $\gamma$. The numerics confirmed our predictions about the form
of the first return map and it also revealed certain additional features.

For numerical construction of the first-return map, we fix angle $\bar \theta \in [0, 2\pi)$
and define a crossection in cylindrical coordinates:
$
\Sigma=\{(\rho, \xi, \theta): \theta=\bar \theta\}.
$
Let $(\rho(t), \xi(t), \theta(t))$ be a trajectory of (\ref{6.5}) and $0<t_1<t_2<t_3<\dots<t_k<\dots$
denote a sequence of times, at which $(\rho(t_k),\xi(t_k),\theta(t_k))\in \Sigma,\;k=1,2,3,\dots$.
Then, we define 
\be\lbl{9.6}
\mG: \mI(t_k) \mapsto \mI(t_{k+1})\quad\mbox{and}\quad \mI(t) =\frac{1}{\rho^2(t)}.
\ee
By a suitable choice of  $\bar\theta$, we can achieve that a trajectory of
(\ref{6.5}), intersects $\Sigma$ once during one cycle of oscillations. Therefore,
for a trajectory lying in a $2D$ slow manifold, (\ref{9.6}) defines a 
$1D$ first-return map. Note that maps (\ref{9.4}) and (\ref{9.6}) are 
related via rescaling of the coefficients, since to leading order $ \mI=\epsilon^{-2}I$.
Therefore, the graphs of $G(I)$ and $\mG(\mI)$ are
similar in the domain, where (\ref{9.4}) is valid. We computed the first return maps
for fixed $\alpha=0.0445$ and for several values of $\gamma$. 
The representative plots are shown in Figure \ref{f.8} a-c.
The numerically computed maps confirmed  our expectations about the graph of
$\mG$ in the boundary layer near $0$: in this region, the map is decreasing
with a strongly expanding slope. 
Away from $0$, the graph of $\mG$ contains an almost linear branch,
whose slope is very close to $1-2\alpha\omega,$ as predicted by (\ref{9.4})
(see lower branches in Figure \ref{f.8}b,c).
These numerics also confirm that under the variation of $\gamma$, the graph of the first return
map is translated in vertical direction (Figure \ref{f.8} a-c). Therefore, the family of the first return maps 
possesses two principal ingredients (the unimodality and the additive dependence on $\gamma$),
which are necessary for the qualitative explanation of the bifurcation scenario 
given in Figure \ref{f.7}. In addition, our numerical experiments reveal
a new feature of the first return map: for small $|\gamma|$, the graph of the map has 
another almost linear branch in the outer region away from the origin (see Figure \ref{f.8}b).
This branch of the graph of the map also has an attracting slope.
The presence of this branch in the first return map indicates the existence of another (branch of the) 
slow manifold different from that described in Section 3. 
A possible explanation for the appearance
of the upper branch in the first return map is due to the unstable manifold of the periodic orbit.
For small $|\gamma|$, the first-return map (\ref{9.6}) is multivalued 
in the outer region. However, since both the upper and the lower branches have positive slopes less than $1$, 
the qualitative dynamics for the map shown in Figure \ref{f.8}b does not depend on the exact mechanism for
the selection between the branches in (\ref{9.6}). Although explaining the multivaluedness of 
the first return maps shown in Figures \ref{f.8}b,c presents an interesting problem, it is not 
critical for the qualitative explanation of the bifurcation scenario arising during the transition from 
the subcritical to supercritical AHB. The latter was the main goal of the present subsection.

\noindent
{\bf Remark 4.2 }
In addition, to the region in the parameter space containing the border between the regions of
sub- and supercritical AHB, there is another parameter regime resulting in
the complex dynamics.
In \cite{FKRT}, it was shown numerically that the limit cycle born from the supercritical
AHB in (\ref{6.5}) undergoes a period-doubling cascade leading to the formation of the chaotic attractor
for increasing values of $\alpha$. The first period-doubling bifurcation in this cascade is shown in Figure \ref{f.1}c,d.
This bifurcation scenario is consistent with the form of the first-return map constructed 
in this subsection. Indeed, it is easy to see from (\ref{9.5}) that for increasing values of $\mu$ the fixed
point, $\bar I$, moves to the left. Therefore, the explanation given above for the period-doubling cascade
resulting from the variation of $\gamma$ near $0$ also applies to the case of increasing $\alpha$ and
fixed $\gamma <0$ (see Figure \ref{f.8}d).

\section{The proof of Theorem 4.1}
\setcounter{equation}{0}

In the present section, we show that under the assumptions of Theorem 4.1
the trajectories of (\ref{2.5}) with initial conditions in $\Sigma^-$ enter $D_0$.
The local analysis in Section 3 describes the behavior of trajectories from
the moment they reach $D_0$ and until they leave $D$. In particular,
it shows that the dynamics near the origin has two phases: the fast approach to 
the slow manifold and the slow drift away from the origin along the slow manifold.
Upon leaving $D$, the trajectories are reinjected back to $\Sigma^-$ by the  
return mechanism postulated in {\bf (G1)}. This scenario implies that the 
system undergoes multimodal oscillations as stated in Theorem 4.1.
\begin{figure}
\begin{center}
\epsfig{figure=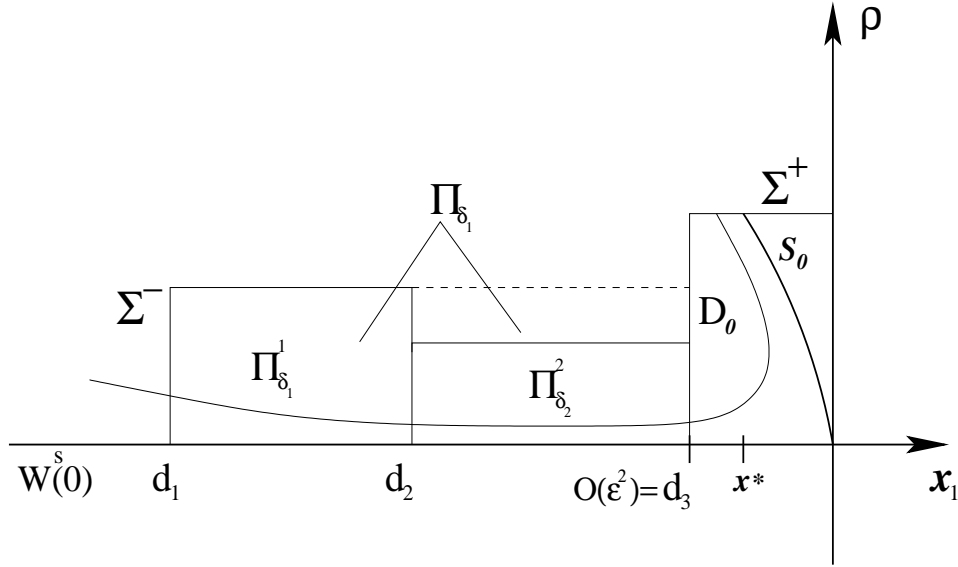, height=3.0in, width=5.0in, angle=0}
\end{center}
\caption{ A schematic representation of the trajectory of (\ref{6.1})-(\ref{6.3})  approaching
a weakly unstable saddle-focus. The trajectory crossing $\Sigma^-$ is subject to strong vector
field transverse to $W^s(O)$ in $\Pi^1_{\delta_1}$. This guarantees that it remains close to $W^s(O)$
until it enters an $O(\epsilon^2)$ neighborhood of the origin. From that moment and until the trajectory
reaches $\Sigma^+$ its behavior is described by the local analysis of Section 3.
}
\label{f.9}
\end{figure}

We start with presenting several auxiliary estimates, which will be needed for the
proof. By (\ref{4.3a}), one can choose $C_1>0$ such that
\be\lbl{7.1}
w(x_1, \alpha)\ge C_1 \left|x_1\right|,\quad x_1\in\left[ d_1,0\right]
\ee 
and for sufficiently small $\alpha\ge 0$.
Next, we note that $\bar\lambda(0,0)=0$. This equation and (\ref{4.6}), by  the Implicit
Function Theorem, imply that there exists $0>x^*(\alpha)=O(\epsilon^2)$ such that
$$
\bar\lambda(x_1,\alpha)\left\{\begin{array}{ll}
<0,& x_1\in[d_1, x_1^*(\alpha)),\\
>0,& x_1\in(x_1^*(\alpha),0].
\end{array}
\right.
$$

Using (\ref{4.6}) and (\ref{7.1}), we have 
$$
{\bar\lambda(x_1,\alpha)\over w(x_1,\alpha)}\ge {\bar\lambda(x_1^*,\alpha)\over w(x_1,\alpha)}=0,\quad x_1\in [d_1, x^*_1(\alpha)],
$$
where $w(x_1,\alpha)$ is defined in (\ref{4.3a}).

Let $\bar\mu>0$ be such that ({\bf G1}) and ({\bf G2}) hold for $\alpha\in [0,\bar\mu\epsilon^2]$. 
In the remainder of this section, unless stated otherwise,  it is assumed that $\alpha\in 
\left[0,\bar\mu\epsilon^2 \right]$. To simplify notation, we will often omit the dependence 
of various functions on $\alpha$.
Using (\ref{4.3a}), we rewrite (\ref{6.5}) in the following form:
\begin{eqnarray}\lbl{7.2}
\dot x_1 &=& w(x_1)+y_1\phi_1(x_1,y)+y_2\phi_2(x_1,y),\\ 
\lbl{7.3}
\dot y &=& A(x_1)y+\psi(x_1,y),
\end{eqnarray}
where $w(x_1)$ and $A(x_1)$ are defined in (\ref{4.3a}) and (\ref{4.5}) respectively;
$y=(y_1,y_2):=(x_2,x_3)$, and 
\be\lbl{7.4}
w(0)=0,\; \psi(x_1,y)=O\left(\left|y\right|^2\right),\; \phi_{1,2}(x_1,y)=\phi_{1,2}(x_1)+O(\left|y|\right),\quad
x\in \left[d_1, 0\right].
\ee
Let 
\be\lbl{7.4a}
d_3=-M^2\epsilon^2
\ee
denote the $x_1-$coordinate on the left lateral boundary of $D_0$ (see (\ref{5.2})) and $\Pi_{\delta_1}=[d_1,d_3]\times D_{\delta_1}$
(see Figure \ref{f.9}).
From (\ref{6.5}) and (\ref{4.3a}), one can see that  for ${\delta_1}>0$ sufficiently small (independent of $\epsilon$),
the right hand side of (\ref{7.2})
$$
w(x_1)+y_1\phi_1(x_1,y)+y_2\phi_2(x_1,y)>0,\quad (x_1,y)\in\Pi_{\delta_1}.
$$
Therefore, in $\Pi_{\delta_1}$, (\ref{7.2}) and
(\ref{7.3}) may be rewritten as follows
\be\lbl{7.5}
{d y(x_1)\over d x_1}=\tilde A(x_1) y +\tilde\psi(x_1,y),
\ee
where 
$$
\tilde A(x_1)={1\over w(x_1)} A(x_1)\quad\mbox{ and}\quad \tilde\psi(x_1,y)={\psi(x_1,y)\over w(x_1,y)}\left(1+O\left(\left|y\right|\right)\right).
$$
By (\ref{7.1}) and (\ref{7.4}), we have
\be\lbl{7.6}
\left|\tilde\psi(x_1,y)\right|\le {C_2\left|y\right|^2\over \left|x_1\right|},\quad (x_1,y)\in\Pi_{\delta_1},
\ee
for some $C_2>0$ independent of $\epsilon$.

To follow the trajectories from $\Sigma^-$ to $D_0$, we introduce two regions:
$$
\Pi_{\delta_1}^1=[d_1,d_2]\times D_{\delta_1}\quad\mbox{and}\quad\Pi_{\delta_2}^2=[d_2,d_3]\times D_{\delta_2}\;\;\mbox{(see Figure \ref{f.9})}.
$$
Recall that $D_{\delta_{1,2}}$ denotes the disk of radius $\delta_{1,2}$ (see (\ref{4.4a})).
Positive constant $\delta_1$ is the same as in (\ref{4.4}) and $0<\delta_2\ll \delta_1$ will be specified later (Figure \ref{f.9}). 
Recall that $d_3=-M\epsilon^2$ denotes the value of the $x_1-$coordinate on the the left lateral boundary of $D_0$ (see (\ref{5.2})).
By taking $M>0$ large enough, we can arrange  $d_3<x^*(\alpha)<0$ for $\alpha\in \left[0,\bar\mu\epsilon^2 \right]$.
By {\bf(G1)}, the vector field in $\Pi^1_{\delta_1}$ is sufficiently strong so that the trajectories entering $\Pi_{\delta_1}^1$
through $\Sigma^-$ get into a narrow domain $\Pi_{\delta_2}^2$ and remain there until they reach $D_0$.
The following lemma allows to control the trajectories in $\Pi_{\delta_1}=[d_1,d_3]\times D_{\delta_1}$.

\noindent{\bf Lemma 5.1}
{\sc Let $[\ell_1,\ell_2]\subseteq[d_1,d_3],\; 0<\bar\delta\le {\delta_1},$ and
\be\lbl{7.7}
-\upsilon=\max_{x_1\in[\ell_1,\ell_2]}\sup_{\left|y\right|=1}\left(\tilde A(x_1)y,y \right)<0.
\ee
Then
\be\lbl{7.8}
\left|y(\ell_1)\right|\le\bar\delta\quad\Rightarrow\quad \left|y(x_1)\right|\le 2\bar\delta e^{-\upsilon(x_1-\ell_1)},\; x_1\in[\ell_1,\ell_2],
\ee
provided
\be\lbl{7.9}
{4 C_2\bar\delta \left(1- e^{-\upsilon(\ell_2-\ell_1})\right)\over\left|\ell_2\right|\upsilon}\le 1.
\ee
}
\noindent
{\bf Proof:}
We may assume that $|y(x_1)|\not=0,x \in [\ell_1,\ell_2]$, since, otherwise, 
by the uniqueness of solution of the initial value problem for (\ref{7.5}), $y(x_1)=0,$ $x_1\in[\ell_1,\ell_2]$,
and (\ref{7.8}) holds. From (\ref{7.5}), (\ref{7.6}), and (\ref{7.7}), we have
\be\lbl{7.10} 
\frac{d|y(x_1)|}{dx_1} \le -\upsilon \left|y(x_1)\right| + \bar\psi\left(x_1,\left|y\right|\right),\;\;
\bar\psi\left(x_1,\left|y\right|\right)={ C_2 \left|y \right|^2 \over \left|x_1 \right|},\;\; x_1\in[\ell_1,\ell_2],\;
\left|y\right|\le\left|\delta_1\right|.
\ee
where $\upsilon$ is defined in (\ref{7.7}), and $\left|y\right|$ denotes the Euclidean norm of $y\in\Re^2$.
Let $\bar y(x_1)$ denote the solution of the initial value problem 
\be\lbl{7.11}
\bar y^\prime=-\upsilon \bar y +\bar\psi(x_1,\bar y),\;
\bar y(\ell_1)=\bar\delta.
\ee
It is sufficient to show that (\ref{7.8}) holds for
$\bar y(x_1)$. We represent $\bar y(x_1)$ as the limit of the sequence
of successive approximations
\be\lbl{7.12}
\bar y^{(n+1)}(x_1)=\bar\delta e^{-\upsilon(x_1-\ell_1)}+\int_{\ell_1}^{x_1} e^{-\upsilon(x_1-s)}
\bar\psi\left(s,\bar y^{(n)}(s)\right) ds,\quad n=0,1,2,\dots
\ee
and $\bar y^{(0)}(x_1)=\bar \delta e^{-\upsilon(x_1-\ell_1)}.$
We use induction to show 
\be\lbl{7.13}
\left|\bar y^{(n)}(x_1)\right| \le 2 \bar\delta e^{-\upsilon(x_1-\ell_1)},\quad x \in [\ell_1,\ell_2],\quad
n=0,1,2,\dots .
\ee 
Inequality (\ref{7.13}) holds for $n=0$. We show that $(\ref{7.13})_{n=k} \Rightarrow (\ref{7.13})_{n=k+1}$.
Using the definition of $\bar\psi$ in (\ref{7.10}) and $(\ref{7.13})_{n=k}$, we have
\begin{eqnarray}\nonumber
\int_{\ell_1}^{x_1} e^{-\upsilon(x_1-s)}\bar\psi(s,\bar y^{(k)}(s) ds &\le&
\frac{C_2}{|\ell_2|}\int_{\ell_1}^{x_1} e^{-\upsilon(x_1-s)}\left|\bar y^{(k)}(s)\right|^2 ds\le
\frac{C_2}{|\ell_2|}\int_{\ell_1}^{x_1}e^{-\upsilon (x_1-s)} 4 {\bar\delta}^2 e^{-2\upsilon(s-\ell_1)} ds\\ 
\lbl{7.14}
&=&\bar\delta e^{-\upsilon(x_1-\ell_1)} {4 C_2\bar\delta \left( 1- e^{-\upsilon(\ell_2-\ell_1)}\right)\over\left|\ell_2\right|\upsilon}.
\end{eqnarray}
Using (\ref{7.14}) and (\ref{7.9}),
from $(\ref{7.13})_{n=k}$ we obtain  $(\ref{7.13})_{n=k+1}$. By induction, (\ref{7.13}) holds. By taking $n\to\infty$ 
in (\ref{7.13}), we have 
\be\lbl{7.15}
\left|\bar y(x_1)\right|\le 2\bar\delta e^{-\upsilon(x_1-\ell_1)},\; x_1\in[\ell_1,\ell_2].
\ee
The statement in (\ref{7.8}) then follows from (\ref{7.15}) and the standard theory for differential inequalities
\cite{hartman}.\\
$\Box$

In the following lemma, we determine the size of $\Pi^2_{\delta_2}$.

\noindent{\bf Lemma 5.2}
{\sc There exists positive $\delta_2=O(\epsilon^2)$, such that the trajectories of (\ref{7.2}) and (\ref{7.3}) 
entering $\Pi^2_{\delta_2}$ from $\Pi^1_{{\delta_1}}$ remain in $\Pi^2_{\delta_2}$ until they reach $D_0$ (see Figure \ref{f.9}.)
}

\noindent {\bf Proof:} 
Denote
$$
-\upsilon_2=\max_{x_1\in[d_2,d_3]}\sup_{\left|y\right|=1}\left(\tilde A(x_1)y, y\right).
$$
Recall 
$$
\bar\lambda (x^*)=0 \quad \mbox{and}\quad O(\epsilon^2)= d_3< x^*=O(\epsilon^2).
$$
Therefore, by (\ref{4.6}), $\bar\lambda(d_3)=O(\epsilon^2)$ is positive.
In addition,
\begin{eqnarray*}
\upsilon_2&=&-\max_{x_1\in[d_2,d_3]} \sup_{\left|y\right|=1}\left(\tilde A(x_1)y, y\right)=
\min_{x_1\in[d_2,d_3]}{-\sup_{\left|y\right|=1}\left( A(x_1)y, y\right)\over w(x_1)}\\
&=& {\min_{x_1\in[d_2,d_3]} \bar\lambda(x_1)\over\max_{x_1\in[d_2,d_3]} w(x_1)}= 
{\bar\lambda(d_3)\over\max_{x_1\in[d_2,d_3]} w(x_1)}= O(\epsilon^2). 
\end{eqnarray*}
and $\upsilon_2>0$. 
Next we apply Lemma 5.1 with $\bar\delta:=\delta_2,$ $\upsilon:=\upsilon_2$, and $\ell_{1,2}:=d_{2,3}$. Note that for small $\upsilon_2$
the inequality (\ref{7.9}) can be rewritten as 
\be\lbl{7.16}
{4 C_2\delta_2 \left(d_3-d_2 +O(\upsilon_2^2)\right)\over\left|d_3\right|}\le 1.
\ee
Since $d_3=O(\epsilon^2)$  (see (\ref{7.4a})) and $\upsilon_2=O(\epsilon^2)$, one can choose $\delta_2=O(\epsilon^2)$ so 
that (\ref{7.16}) holds.\\
$\Box$

Having found the size of $\Pi^2_{\delta_2}$, we now determine the rate of contraction in $\Pi^1_{\delta_1}$ sufficient
to funnel the trajectories entering $\Pi^1_{\delta_1}$ through $\Sigma^-$ to $\Pi^2_{\delta_2}$.

\noindent{\bf Lemma 5.3}
{\sc The trajectories of (\ref{7.2}) and (\ref{7.3}) 
entering $\Pi^1_{\delta_1}$ through $\Sigma^-$ remain in $\Pi^1_{\delta_1}$ until they reach 
$\Pi^2_{\delta_2}$.
}

\noindent {\bf Proof:} 
Denote
$$
-\upsilon_1=\max_{x_1\in[d_1,d_2]}\sup_{\left|y\right|=1}\left(\tilde A(x_1)y, y\right).
$$
From {\bf (G2)}, one finds that $\upsilon_1=O(\left|\ln \epsilon\right|)$ is positive.
Lemma 5.3 now follows from Lemma 5.1 with  $\bar\delta:={\delta_1},$ $\upsilon:=\upsilon_1$, and $\ell_{1,2}:=d_{1,2}$.  
Indeed, for $\upsilon_1=O\left(\left|\ln\epsilon\right|\right)$, we have 
\be\lbl{7.17}
\upsilon_1\ge {4C_5{\delta_1}\over\left|d_2\right|}.
\ee
Inequality (\ref{7.17}) is sufficient for (\ref{7.9}) to hold. By Lemma 5.1, we have
\be\lbl{7.19}
\left|y(d_2)\right|\le 2{\delta_1} e^{-\upsilon_1(d_2-d_1)}.
\ee
With $\upsilon_1=O\left(\left|\ln\epsilon\right|\right)$, by (\ref{7.19}), we can achieve 
$\left|y(d_2)\right|<\delta_2=O(\epsilon^2)$.\\
$\Box$

Lemmas 5.2 and 5.3 imply that the trajectories entering $\Pi^1_{\delta_1}$ through $\Sigma^-$
stay in $\Pi^1_{\delta_1}\bigcup\Pi^2_{\delta_2}$ until they reach $D_0$. Moreover, the
inequality in {\bf (G1)} guarantees that such trajectories are bounded away from $W^s(O)$.
Thus, we can use the analysis of Section 3 to describe the evolution of trajectories 
from the moment they reach $D_0$ until they leave $D$. By Remark 3.1c, this description extends
to any region where $x_{2,3}=o\left(\epsilon^{2/3}\right),$ i.e., the trajectories can be controlled 
until they hit $\Sigma^+$. This is followed by the return to $\Pi_{\delta_1}^1$ according to {\bf (G1)},
and the next cycle of the multimodal oscillations begins. The analysis of this section applies to 
any trajectory starting from $D_0$ and not belonging to $W^s(O)$.

It remains to estimate the ISIs. For this, we compute the time needed for the trajectory starting
in an $O(\epsilon^2)$ neighborhood to return back to this neighborhood after making one global 
excursion. Since the time of flight of the trajectory outside a small neighborhood of the 
origin depends regularly on the control parameters, the duration of the very long ISIs is determined
by the time spent in the neighborhood of the origin. To estimate the latter, 
we note that at the moment a trajectory enters $D_0$, we have  
\be\lbl{7.20a}
\left|y(d_3)\right|\le C_3 \epsilon^2.
\ee
This follows from Lemma 5.3, since $\left|y(d_3)\right|\le \delta_2$ and $\delta_2=O\left(\epsilon^2\right).$
To obtain the lower bound on $\left|y(d_3)\right|$, recall that by {\bf (G1)}, we have 
\be\lbl{7.20b}
\left|y(d_1)\right|\ge \zeta>0.
\ee
As follows from {\bf (G3)}, the maximal rate of contraction in $\Pi_{\delta_1}$ does not exceed 
$O\left(\left|\ln\epsilon\right|\right)$ in absolute value. This combined with (\ref{7.20b}) implies that
\be\lbl{7.20c}
\left|y(d_3)\right|\ge C_4\epsilon^{2+{\chi\over 2}}
\ee
for some $\chi\ge 0$ and $C_4$ independent of $\epsilon$. We omit the proof of (\ref{7.20c}), because it is 
completely analogous to that of Lemma 5.1. 

Let $t=t^-$ denote the moment 
of time when a trajectory 
of (\ref{6.5}) enters $D_0$ from $\Pi^2_{\delta_2}$. After switching back to the original parametrization 
of $y$ by time, $t$, we rewrite (\ref{7.20a}) and (\ref{7.20c}):
\be\lbl{7.20}
C_3\epsilon^{2+{\chi \over 2}}\le\left|y(t^-)\right|\le C_4\epsilon^2.
\ee
For $t>t^-$ the trajectory approaches and remains close to the slow manifold as long as 
$\left|y\right|=o(\epsilon^{2/3})$
(see Remark 3.1c). Let ${2\over3}<j<1$ and denote 
\be\lbl{7.21} 
t^+=\min\left\{t>t^-:\;\left|y(t)\right|=\epsilon^j\right\}.
\ee
From (\ref{7.20}) and (\ref{7.21}), we have  the following bounds for $I_0:=I(t^-)$ 
and $I_1:=I(t^+)$:
\be\lbl{7.22}
C_6\epsilon^{-2-\chi}\le I_0\le C_7\epsilon^{-2}
\quad\mbox{and}\quad
I_1=O(\epsilon^{2(1-j)}),\; j\in\left({2\over 3}, 1\right).
\ee
For these ranges of values of $I_0$ and $I_1$, from (\ref{5.55}) we have 
\be\lbl{7.23}
\tau_{in}=t^+-t^-={1\over 2\mu\epsilon^2}\ln\left(\left(1+{\mu\over\gamma}I_0\right)
\left(1+{\mu\over\gamma}o(1)\right)\right).
\ee
The combination of (\ref{7.22}) and (\ref{7.23}) yields two-sided bounds for $\tau_{in}$:
\be\lbl{7.24}
\tau^-_{in}\le \tau_{in}\le \tau_{in}^+,
\ee
where
$$
\tau_{in}^- ={1\over 2\mu\epsilon^2}
\ln\left(1+{\mu\over\gamma\epsilon^2}C_8^-\right)\quad\mbox{and}\quad
\tau_{in}^+ ={1\over 2\mu\epsilon^2}
\ln\left(1+{\mu\over\gamma\epsilon^{2+\chi}}C_8^+\right),
$$
where positive constants $C^\mp_8$ can be chosen independent of $\epsilon$.
Inequalities (\ref{7.24}) provide bounds for the time that a multimodal trajectory spends
near the unstable equilibrium. On the other hand, the time of flight outside a small neighborhood 
of the origin, $\tau_{out}$, depends regularly on $\epsilon$, by {\bf (G2)}.
Therefore, by adjusting constants $C_8^\mp$ in (\ref{7.24}) if necessary, one can obtain 
uniform bounds on the ISIs $\tau=\tau_{in}+\tau_{out}$ for sufficiently small $\epsilon>0$,
positive $\gamma$ and $\mu$ from bounded intervals, as stated in Theorem 4.1.

\section{Discussion}
\setcounter{equation}{0}

In the present paper, we investigated a mechanism for generation of multimodal oscillations
in a class of systems of differential equations close to an AHB. 
Our analysis covers both cases of sub- and supercritical AHB. 
For the supercritical case, we identified a novel geometric feature of the bifurcating limit cycle, 
the frequency doubling effect.
It turns out that generically in the normal system of coordinates the oscillations 
in one of the variables are twice faster than in the remaining two variables.
Therefore, the leading order approximation of the limit cycle bifurcating from the supercritical AHB
requires two harmonics.
The asymptotic analysis of the present paper explains the frequency-doubling. In addition, we provide
a complementary geometric interpretation to this counterintuitive effect. In particular,
we showed that it is a consequence of the geometry of the limit cycle. 
The latter is not captured by the topological normal form of the AHB.
The analysis of the multimodal oscillations arising from the subcritical AHB requires
additional assumptions on the global behavior of trajectories. 
We identified two principal properties of the global vector field:
the mechanism of return and the strong contraction property.
In the presence of this global structure, the subcritical AHB produces
sustained multimodal oscillations combining the small amplitude oscillations near
the unstable equilibrium with large amplitude spikes. 
The resultant motion is recurrent in a weak sense:
it may not be periodic but nevertheless the timings of the spikes possess certain
regularity. We have shown that the ISIs have well-defined asymptotics near the AHB
and comply to the two-sided bounds, which depend on the principal bifurcation
parameters. Our estimates show that near the AHB, the ISIs can be extremely long
and can change greatly under relatively small variation of the bifurcation parameters.
The ability of the system to exhibit such extreme variability in the ISI duration is
important in many applications, in particular, in the context of neuronal dynamics.
Previous studies investigated different possible mechanisms for generating multimodal
patterns with very long ISIs \cite{BER,DK05,DRSE,MC,RE89}.
For the finite dimensional approximation of the model
of solid-fuel combustion (\ref{6.1})-(\ref{6.3}), the main motivating example for our work, 
the proximity to the homoclinic bifurcation was suggested in \cite{FKRT} as a possible mechanism
for prolonged ISIs. 
Our conclusions confirm the importance of the proximity to the homoclinic bifurcation for explaining
the oscillatory patterns in (\ref{6.1})-(\ref{6.3}). The proximity to the homoclinic bifurcation
is implicitly reflected in our assumptions
on the global vector field.
However, we emphasize the critical role of the AHB:
the duration of the ISIs can be effectively controlled by the parameters associated with the AHB
without changing the distance of the system to the homoclinic bifurcation.
In all our numerical experiments, the system remained bounded away from the homoclinic bifurcation,
nevertheless it exhibited patterns with very long ISIs whose duration was amenable to control.
Our analysis extends the estimate for the ISIs obtained in \cite{GW} to a wide class of problems. 
It also emphasizes the proximity of the system to the border between 
sub- and supercritical AHB, as another factor in creating oscillatory patterns with long ISIs.
We show that as this border is approached from the subcritical side
the ISI grow logarithmically. This observation is important for explaining the oscillatory patterns
generated by (\ref{6.1})-(\ref{6.3}) since in this model the AHB can change its type under the variation
of the second control parameter $p$. This situation is not specific to the model of 
solid fuel combustion. Recent studies suggest 
that there is a class of neuronal models close to the AHB whose type may change with the values of parameters 
\cite{DK05, OP}.
Therefore, it is important to understand the principal features of the transition from sub- to 
supercritical AHB. We found that when the border between the regions of sub- and supercritical
bifurcation is approached from the subcritical side (while the distance from the AHB remains fixed), 
the oscillations become chaotic. The regime of irregular oscillations is then followed by the 
reverse period doubling cascade. To understand the nature of this bifurcation scenario 
we used a combination of analytic and numerical techniques.
Using the insights gained from the asymptotic analysis, we constructed a $1D$ first-return map.
The map provides a clear geometric interpretation for the bifurcation scenario near the 
transition from sub- to supercritical AHB. Our study suggests that the formation of the chaotic attractor
via a period-doubling cascade is a universal feature of this transition. For example, the bifurcation
scenarios reported for the Hodgkin-Huxley model in \cite{DK05} are very similar to those studied in the
present paper and are likely to share the same mechanism.

Mixed-mode oscillations similar to those studied in the present paper, have been studied in for a class
of the slow-fast systems in $\Re^3$ near the AHB. Although, the work toward developing a complete mathematical 
theory for such oscillations is still in progress, the general mechanism for their generation and the bifurcation
structure of the problem have been greatly elucidated  recently \cite{KOP, MC, MSLG, WESCH}.
The present paper shows the relation between the mechanisms for the mixed-mode oscillations in the model
in \cite{FKRT} and for those in the slow-fast systems. The latter possess a well-defined structure of the
global vector field due to the presence of the disparate timescales in the governing equations \cite{JON}. The analyses
of the mixed-mode oscillations in \cite{KOP, MC, MSLG, WESCH} use in an essential way the relaxation structure 
of the problem. The model in \cite{FKRT} is an example of the mixed-mode generating system, which does not 
possess an explicit relaxation structure. In fact, it is hard to expect such structure in a system obtained via
projecting an infinite-dimensional system onto a finite-dimensional subspace. In formulating the assumptions on the
global vector-field ({\bf G1}) and ({\bf G2}), we were looking for the minimal requirements on the system near
an AHB that guarantee the existence of the mixed-mode solutions. Due to the lack of the information about the
global vector field of (\ref{6.1})-(\ref{6.3}), it appears impossible to verify these conditions analytically.
However, the numerical simulations clearly show that system of equations (\ref{6.1})-(\ref{6.3}) possesses
the qualitative structure required by ({\bf G1}) and ({\bf G2}) (Figure \ref{f.2}b). On the other hand, we expect that conditions 
({\bf G1}) and ({\bf G2}) should be possible to verify analytically for a wide class of slow-fast systems.
Therefore, we believe that our results will be useful for understanding mixed-mode oscillations in such systems. 
In particular, it would be interesting to apply this approach to the modified Hodgkin-Huxley system \cite{DK05}.
The numerical results reported in \cite{DK05} strongly suggest that the mechanism proposed in the present paper
is responsible for the generation of the very slow rhythms and chaotic dynamics in the Hodgkin-Huxley model.

\noindent
{\bf Acknowledgments.}
We thank Victor Roytburd for introducing us to this problem and to Michael Frankel and 
Victor Roytburd for helpful conversations. This work  was
partially supported through National Science Foundation Award No 0417624.

\renewcommand{\theequation}{A.\arabic{equation}}
\section*{Appendix A.} 
\setcounter{equation}{0}
\label{sec:C}
In this appendix, we list 
explicit expressions of various constants and trigonometric polynomials, which appear
in the definitions in of the slow manifold, $S$, and the first Lyapunov coefficient, $\gamma$.
All expressions are given in terms
of the coefficients of the power expansions on the right  hand side of (\ref{5.1}).
By $\sigma_{ijk}$ we denote $\sigma_{ijk}(0),$ $\sigma\in\{a,b,c\}$.
The following constants are used to define the leading order approximation of the slow
manifold in (\ref{5.3}) and (\ref{5.4}):
$$
a={a_{22}+a_{33}\over 2},\;
A=\sqrt{a^2_{23}+\left(a_{22}-a_{33}\right)^2\over 4\left(1+ 4\beta^2\right)},\;\mbox{and}\;
\vartheta=\arctan {2\beta a_{23}+\left(a_{22}-a_{33}\right)\over  a_{23}+2\beta \left(a_{22}-a_{33}\right)} .
$$
The following trigonometric polynomials enter the right hand sides of 
(\ref{5.11a}) and (\ref{5.12}): 
\begin{eqnarray*}
P_{1}(\phi)&=& a_{22}\cos^2\beta\phi+a_{23}\cos\beta\phi\sin\beta\phi+a_{33}\sin^2\beta\phi,\\
Q_{1}(\phi)&=&b_{22}\cos^3\beta\phi+(b_{23}+c_{22})\cos^2\beta\phi\sin\beta\phi
+(b_{33}+c_{23})\cos\beta\phi\sin^2\beta\phi+c_{33}\sin^3\beta\phi,\\
Q_{2}(\phi)&=&b_{12}\cos^2\beta\phi+(b_{13}+c_{12})\cos\beta\phi \sin\beta\phi
+c_{13}\sin^2\beta\phi,\\
Q_{3}(\phi)&=&b_{222}\cos^4\beta\phi+(b_{223}+c_{222})\cos^3\beta\phi \sin\beta\phi
                  +(b_{233}+c_{223})\cos^2\beta\phi \sin^2\beta\phi\\
                 &&+(b_{333}+c_{233})\cos\beta\phi \sin^3\beta\phi+c_{333}\sin^4\beta\phi,\\
\end{eqnarray*}
Functions $\bar Q_{2,3}(\theta)=Q_{2,3}(\beta^{-1}\theta)$ are used in the calculation
of the first Lyapunov coefficient $\gamma$ (see (\ref{5.5})).

\end{document}